\documentclass[11pt,a4paper]{amsart}
\usepackage[T1]{fontenc}
\usepackage{lmodern}
\usepackage[a4paper,textwidth=150mm,textheight=228mm,centering]{geometry}
\usepackage{amssymb,mathtools,microtype,booktabs,tabularx,array}
\usepackage[unicode,hidelinks,bookmarksnumbered=true]{hyperref}
\hypersetup{pdfauthor={Masoud Bayrami},
  pdftitle={Concentration in Hardy equations and singular heat flows}}
\numberwithin{equation}{section}
\newtheorem{theorem}{Theorem}[section]
\newtheorem{proposition}[theorem]{Proposition}
\newtheorem{lemma}[theorem]{Lemma}
\newtheorem{corollary}[theorem]{Corollary}
\newtheorem{suppprop}{Proposition}[theorem]

\newtheorem{suppcor}[suppprop]{Corollary}
\newtheorem*{condenserdistance}{Lemma}
\theoremstyle{remark}
\newtheorem{remark}[theorem]{Remark}
\newcommand{\R}{\mathbb R}
\newcommand{\T}{\mathbb T}
\newcommand{\Cyl}{\mathcal C}
\newcommand{\HH}{\mathcal H}
\newcommand{\Hrad}{\mathcal H_{\mathrm{rad}}}
\newcommand{\B}{\mathcal B}
\newcommand{\eps}{\varepsilon}
\newcommand{\sph}{|\mathbb S^{N-1}|}
\newcommand{\ecap}{\operatorname{Cap}}
\newcommand{\supp}{\operatorname{supp}}
\newcommand{\TV}{\mathrm{TV}}
\newcommand{\doi}[1]{\href{https://doi.org/#1}{\nolinkurl{#1}}}
\allowdisplaybreaks[1]
\begin{document}
\title[Concentration in Hardy equations and singular heat flows]
  {Concentration in Hardy equations\\and singular heat flows}
\author[Masoud Bayrami]{Masoud Bayrami}
\address{Department of Mathematical Sciences, Sharif University of Technology,
  Azadi Avenue, P.O. Box 11365-9415, Tehran, Iran}
\email{masoud.bayrami1990@sharif.edu}
\date{}

\begin{abstract}
Concentration of source measures can destroy continuity of nonlinear observables even when the
underlying data converge weakly. We study three manifestations of this phenomenon. For the local
critical Hardy equation we characterize energy admissibility and compactness by radial dyadic sums
under a bounded angular-density hypothesis, with one-sided admissibility results without angular
control, and obtain capacity criteria for general diffuse data at inverse power one. Radial endpoint
traces obey sharp moment bounds and have tightness criteria distinct from energy compactness. At the
sharp fractional $L^p$ Hardy coefficient, a localized power identity separates the Sobolev--Bregman
remainder from a pole flux. We show that a bounded profile has finite remainder precisely when its
nonnegative source pairing is finite, without requiring a pole trace for this equivalence. When the
trace exists, its coefficient is explicit, and finite remainder excludes nonconstant strong angular
traces; positive data with unbounded profiles demonstrate the failure of the equivalence in general.
For singular heat equations with a positive nonincreasing factor vanishing at infinity we identify
the exact range of reaction limits for a single time impulse and for inputs exhausted by parabolic
Kato measures together with specified invisible components, and control the effect of accumulating
impulses. These results distinguish energy defects, endpoint fluxes and loss of effective source
mass.

\end{abstract}
\maketitle

\section{Introduction}\label{sec-intro}

Weak convergence of positive source measures need not control the nonlinear
observables driven by those measures. Concentration may occur near a
singularity, near a boundary, or on small sets in the interior. We study
three distinct ways in which the observable can fail to be determined
continuously by the data: loss of energy compactness, a boundary contribution
to a nonlinear power balance, and loss of effective reaction mass. In each
case we identify the additional information needed to describe the defect.

\subsection*{Energy and traces at the critical Hardy constant}
Let \(B_R\subset\R^N\) be a ball with \(N\ge3\). We first consider the
problem
\begin{equation}\label{physical}
 \left(-\Delta-\frac{(N-2)^2}{4|x|^2}\right)u
   =u^{-\gamma}\mu,\qquad
 u|_{\partial B_R}=0,\qquad \gamma>0 .
\end{equation}
Here \(\mu\) is a positive measure on the punctured ball. We divide the
solution by the Hardy ground state \(|x|^{-(N-2)/2}\) and use the change of
variable \(t=\log(R/|x|)\). The critical Hardy form then becomes the
Dirichlet form on a half-cylinder. The outer boundary corresponds to \(t=0\)
and the pole to \(t\to\infty\). The difference between energy and trace
convergence can be seen for a single spherical source. Indeed, a transformed
measure of mass \(s^{-1}\) at height \(s\) gives the profile
\(\min(t/s,1)\). Its energy is \(s^{-1}\), while its trace at the pole is
one.

Boccardo and Orsina~\cite{BO} studied energy regularity and the loss of
concentrated data for equations with singular nonlinearities. See also the
results for measure data in Orsina and Petitta~\cite{OrsinaPetitta} and De
Cave, Durastanti and Oliva~\cite{DDO}. Hara related solvability to trace
inequalities and proved stability in energy in a metric on measures for
inverse powers at most one~\cite{Hara22,Hara23}. Klimsiak~\cite{Klimsiak}
studied regular symmetric Dirichlet forms and used quasi-uniform convergence
of linear potentials to prove stability. We refer to \cite{OlivaPetitta} for
a survey.

We seek conditions on the data that are equivalent to strong convergence for
the fixed critical form. For radial data and \(0<\gamma<1\), the
admissibility condition is the dyadic form of the weighted Hardy inequality
for the operator \(f\mapsto\int_0^t f(s)\,ds\) from \(L^2\) to
\(L^{1-\gamma}(d\zeta)\), see \cite{HardyMeasure}. For every \(\gamma>0\),
Theorem~\ref{angular} gives a necessary and sufficient condition on the
dyadic tails when the conditional angular densities are uniformly bounded.
In Theorem~\ref{endpoint} we give sharp moment conditions for the two radial
traces. For general diffuse data and \(\gamma=1\), Theorem~\ref{condenser}
characterizes strong convergence in terms of total mass and the difference
of the masses on the two plates of condensers of bounded capacity. For
\(\gamma\ne1\), we use three minima of the energy functional. Each is given
by a formula involving nested condensers, see Theorems~\ref{gap}
and~\ref{levels}. These criteria use the energy capacity of the cylinder.
Hara~\cite{Hara26} and Chen, Frank, Xiao, Yang, Yuan and Zhang~\cite{CFX26}
use capacity to construct Hardy weights and characterize boundary behavior
on general domains. The radial positive-power trace inequality is thus an
application of the established weighted Hardy theory. The additional results
here are the explicit two-ended compactness criterion for the critical
cylinder for every inverse power, the sharp and independent endpoint tests,
and the mass-and-condenser criterion for general diffuse data at inverse
power one. Proposition~\ref{prop:angular} identifies the one-sided
admissibility statements that survive without angular control and gives
counterexamples to the remaining implications. The minimum-action criterion
uses the convexity and lattice properties of the energy form; it is a
complementary formulation rather than a new general theory of trace
inequalities.

\subsection*{The pole contribution in a fractional power balance}
We next consider the linear equation
\[
 (-\Delta)^s u-\lambda |x|^{-2s}u=f
\]
with zero exterior values. At the sharp \(L^p\) Hardy coefficient for
\(p>2\), the power-energy term and the Hardy term may both be infinite. For
the local quadratic form, V\'azquez and Zographopoulos~\cite{VZ,VZ2}
identified the contribution at the pole. Frank and
Seiringer~\cite{FrankSeiringer} developed ground-state representations for
the fractional operator. Bogdan, Jakubowski, Lenczewska and Pietruska-Pa\l
uba~\cite{BogdanLp} proved the sharp \(L^p\) inequality with a
Sobolev--Bregman remainder. Chen and Weth~\cite{CW} proved a weighted
identity at the pole for the linear fundamental solution, and Chen, Gkikas
and Nguyen~\cite{CGN} used it to study measure data. Bogdan and
Merz~\cite{BM24,BM26} obtained angular-channel representations for the
quadratic forms and heat-kernel estimates for their semigroups.

In Theorem~\ref{angular-flux} we give an identity for the nonlinear source
pairing \(\int u^{p-1}f\). We assume that the profile obtained by dividing
by the ground state is bounded and that the spherical mean of its power has
a limit. The coefficient of the trace term is the weighted linear pole
normalization divided by \(p\). In logarithmic coordinates this
normalization is \(\sph\kappa'(\beta)\), where \(\kappa\) is the Hardy
coefficient corresponding to the power \(|x|^{-\beta}\). The remainder also
controls the angular oscillation of the profile. The sharp Hardy constant,
the Sobolev--Bregman representation, and the linear pole normalization are
established results; the new issue here is their combination in a nonlinear
power pairing for a possibly nonradial bounded profile.
Corollary~\ref{cor:notrace} shows that boundedness suffices for the equivalence
of finiteness of the remainder and of the nonnegative source pairing, even
without a limiting trace. In Theorem~\ref{fractional-counter} we construct
positive examples for each \(0<s<1\) with finite remainder and infinite
source pairing, showing why boundedness matters.

\subsection*{Attainable reactions in a heat equation}
Finally, on a flat torus we consider smooth nonnegative data \(a_n\,dt\,dx\)
that converge narrowly to a finite measure \(\Lambda\) and the corresponding
solutions of
\[
 \partial_t u_n-\Delta u_n=a_nF(u_n).
\]
We assume that \(F\) is positive and nonincreasing and tends to zero at
infinity. It may be singular at zero, as in the example
\(F(r)=r^{-\gamma}\). For a pulse with a given spatial density, the scalar
jump is obtained by solving \(dU/dc=F(U)\). This is the usual change of time
used for impulsive controls~\cite{BR,BMazzola}. Narrow convergence, however,
allows the approximating data to concentrate in space during the pulse.

Oliva and Petitta~\cite{OPheat} considered a singular reaction with an
integrable coefficient and an additive measure. Bal and Das~\cite{BD26}
studied mixed local and nonlocal diffusion with variable inverse powers and
integrable data. Das~\cite{Das25} allowed a measure coefficient in the
singular term and proved existence for a class of diffuse measures. For the
regularized local problem with zero initial data and constant inverse power
\(\gamma\ge1\), Das also studied the disappearance of polar data. For
admissible initial and boundary measures in heat equations with absorption,
see \cite{AJV}. Our aim is to determine all limits of \(a_nF(u_n)\,dt\,dx\)
over smooth approximations that converge narrowly to the same measure. For a
single time impulse, Theorem~\ref{atom-range} gives the full range of
densities between zero and the scalar jump. For atomless data exhausted by
parabolic Kato measures, Theorem~\ref{evo2-range} shows that the possible
limits are precisely the canonical reactions of submeasures of the datum.
This applies in particular to all space--time \(L^1\) data. The proof also
specifies which polar parts always disappear. Corollary~\ref{mixed-range}
uses a total variation estimate to combine the two cases when the impulse
times may accumulate. Unlike reduced-measure results that identify a
particular limit or the disappearance of polar data, our range statements
classify all limits arising from smooth nonnegative approximations within
the stated classes of inputs. The canonical reaction for the full input is
only one element of the range.

\subsection*{A common compactness question}
The link between these three parts is the localization of the relevant
observable before a limit is taken. For the local equation, restriction to
compact slabs separates escape at the ends from concentration on sets of
small capacity. The energy and the endpoint moments need separate controls.
For the fractional equation, cutting off the pole leaves a signed boundary
term in addition to the nonnegative remainder. For the heat equation, a
short-delay cutoff of the potential gives an upper bound on any limiting
reaction, and concentrated inputs may be made invisible. The three defects
are not identified with one another: they are, respectively, a norm defect,
a flux in a power identity, and a loss of effective input. The proofs of the
three parts are independent, but their common purpose is to replace an
unqualified assertion of weak stability by precise criteria or exact ranges.

\subsection{Guide to the results and notation}\label{sec-guide}
The three concentration phenomena have different compactness variables.
In the local Hardy problem the main unknown is an energy function on a
half-cylinder; in the fractional problem it is the profile after division
by a singular ground state; in the heat problem it is the \emph{reaction
measure}, which need not have the same limit as the forcing. The table
provides a guide; the full hypotheses are stated in the cited results.

\begin{center}
\small
\renewcommand{\arraystretch}{1.16}
\begin{tabularx}{\textwidth}{@{}p{0.31\textwidth}X@{}}
\toprule
\textbf{Result} & \textbf{Hypotheses and conclusion} \\
\midrule
Theorem~\ref{angular} (Section~\ref{sec-form}) & Bounded conditional angular densities; dyadic admissibility and strong-energy compactness. See Proposition~\ref{prop:angular} for the unrestricted one-sided results.\\
Theorem~\ref{endpoint} (Section~\ref{sec-form}) & Radial sources; sharp moments and separate endpoint-trace tightness criteria.\\
Theorems~\ref{condenser}--\ref{levels} (Section~\ref{sec-geometry}) & Diffuse sources; capacity discrepancies at $\gamma=1$ and minimum-action tests for other powers.\\
Theorems~\ref{angular-flux},~\ref{fractional-counter} (Section~\ref{sec-flux}) & Fractional Hardy operator; pole balance, angular control, and positive counterexamples. Corollary~\ref{cor:notrace} does not require a trace.\\
Theorem~\ref{atom-range} (Section~\ref{sec-impulse}) & A single time impulse; the exact interval of attainable reaction densities.\\
Theorem~\ref{evo2-range} and Corollary~\ref{mixed-range} (Section~\ref{sec-relaxation}) & Kato-exhaustible and mixed inputs; canonical reactions of submeasures and accumulating impulses.\\
Proposition~\ref{evo-modulus} (Section~\ref{sec-clock}) & Spatially absolutely continuous time measures; stretched-time stability estimates.\\
\bottomrule
\end{tabularx}
\end{center}

Here is the terminology needed to read the statements. On
$\Cyl=(0,\infty)\times\mathbb S^{N-1}$, the energy space $\HH$ is the
completion under the Dirichlet norm in Section~\ref{sec-form}. A positive
source $\sigma$ is \emph{admissible} if the equation
$-\Delta_{\Cyl}v=v^{-\gamma}\sigma$ has a positive energy solution whose
reaction $v^{-\gamma}\sigma$ defines a bounded energy functional. We write
$E(\sigma)=\|v_\sigma\|^2$. With $\zeta$ the radial marginal,
\[
 I_k=[2^k,2^{k+1}),\qquad
 D_\gamma(\zeta)=\sum_{k\in\mathbb Z}
 \bigl(2^{k(1-\gamma)/2}\zeta(I_k)\bigr)^{2/(1+\gamma)}.
\]
For compact $K\subset U\Subset\Cyl$, the \emph{condenser capacity}
$\ecap(K,U)$ is the least energy of a compactly supported function in $U$
that is at least one near $K$. The discrepancy $A_L(\sigma,\tau)$ is the
largest positive difference between the mass of one source on $K$ and
that of the other on $U$, with the roles reversed as well, over condensers
of capacity at most $L$; see Section~\ref{sec-geometry}. For
$\gamma\ne1$, the minimum action $m_\gamma$ is the minimum of the power
energy functional, and
\[
 G_\gamma(\sigma,\tau)
 =m_\gamma\bigl((\sigma+\tau)/2\bigr)
  -\tfrac12\bigl(m_\gamma(\sigma)+m_\gamma(\tau)\bigr).
\]
At $\gamma=1$ this difference is defined by a common exhaustion.
The nested-condenser cost $C_\gamma$ gives an equivalent capacity-only
formula for $m_\gamma$ (Theorem~\ref{levels}).
We use $L$ for the angular-density bound within
Theorem~\ref{angular}, for the capacity threshold in $A_L$ within
Theorem~\ref{condenser}, and as a subscript in the operator
$L_\lambda$; these quantities are unrelated. The energy-space symbol
$\HH$ is distinguished from the radial subspace $\Hrad$ and from
$h_\beta$, the fractional ground state.

In the heat problem, write $Q=[0,T]\times\mathbb T^d$ and
$S(t)=e^{t\Delta}$. With $F>0$ decreasing, the scalar jump and its inverse
change of time are
\[
 \Phi(r)=\int_0^r\frac{dz}{F(z)},\qquad
 K=\Phi^{-1},\qquad J_c(r)=K(\Phi(r)+c).
\]
A \emph{submeasure} $\Gamma\le\Lambda$ means that $\Lambda-\Gamma$ is
positive. A parabolic \emph{polar set} has zero capacity for the norm
$L_t^2H_x^1\cap H_t^1H_x^{-1}$; precise cutoffs appear in
\eqref{par-cap}. The one-sided \emph{Kato condition} is
\[
 k_\Lambda(h):=\sup_{t,x}
  \int_{t-h<s<t}p_{t-s}(x-y)\,d\Lambda(s,y)
 \longrightarrow0\qquad (h\downarrow0).
\]
A measure is \emph{exhausted by Kato measures} if it is the increasing
limit of finite Kato submeasures. We denote by $\mathcal R(\Gamma)$ its
\emph{canonical reaction}, constructed first for Kato measures and then
by total variation completion; $\mathcal A(\Lambda)$ is the set of
reaction limits obtained from smooth narrow approximations to $\Lambda$.
In particular, the canonical reaction of a submeasure is one \emph{choice}
of limit, not necessarily the reaction of the entire datum.

For quick reference, the main symbols are collected below.
\begin{center}
\small
\renewcommand{\arraystretch}{1.12}
\begin{tabularx}{\textwidth}{@{}p{0.31\textwidth}X@{}}
\toprule
\textbf{Symbol} & \textbf{Meaning} \\
\midrule
$\Cyl,\ \HH,\ \Hrad$ & Half-cylinder, energy space, and its radial subspace.\\
$\sigma,\ \zeta,\ E,\ D_\gamma$ & Transformed source, radial marginal, energy, and dyadic functional.\\
$\ecap,\ A_L,\ G_\gamma,\ C_\gamma$ & Capacity, condenser discrepancy, action gap, and layered cost.\\
$h_\beta,\ \mathcal R_{s,p}$ & Fractional Hardy ground state and Bregman remainder.\\
$\Phi,\ K,\ J_c$ & Scalar clock, its inverse, and pulse-jump map.\\
$\ecap_{\mathrm{par}},\ k_\Lambda$ & Parabolic capacity and one-sided Kato modulus.\\
$\mathcal R(\Gamma),\ \mathcal A(\Lambda)$ & Canonical reaction and all attainable reaction limits.\\
\bottomrule
\end{tabularx}
\end{center}

A \emph{vague} limit is tested against compactly supported continuous
functions, and a \emph{narrow} limit against all bounded continuous
functions. Loss of mass at either end of $\Cyl$ is compatible with vague
convergence, whereas the heat problem is formulated on the compact
space--time domain $Q$ and uses narrow convergence throughout. Thus energy
compactness, endpoint trace compactness, and reaction compactness require
different additional tightness hypotheses.

\section{Critical Hardy energy and radial traces}\label{sec-form}

\subsection{The cylinder and its energy solutions}
Set
\[
 \begin{gathered}
 b=(N-2)/2,\qquad \phi(x)=|x|^{-b},\qquad
 \mathcal T(x)=(\log(R/|x|),x/|x|),\\
 u=\phi(v\circ\mathcal T),\qquad
 \Cyl=(0,\infty)\times\mathbb S^{N-1}.
 \end{gathered}
\]
We take \(d\omega=dS/\sph\) as the angular probability measure and
write \(dm=dt\,d\omega\). Let \(\HH\) be the
completion of \(C_c^\infty(\Cyl)\) with respect to the norm given by
\[
 \|h\|^2=\B(h,h),\qquad
 \B(v,h)=\int_\Cyl(v_t h_t+\nabla_Sv\cdot\nabla_Sh)\,dm.
\]
For test functions supported away from the pole, integration by parts gives
\[
 \int_{B_R}(|\nabla u|^2-b^2|x|^{-2}u^2)\,dx
   =\int_{B_R}\phi^2|\nabla(u/\phi)|^2\,dx=\sph\|v\|^2 .
\]
The last expression defines the completed critical form~\cite{VZ}. In the
new variables we have
\begin{equation}\label{weight}
 \sigma=\sph^{-1}\mathcal T_\#(\phi^{1-\gamma}\mu),\qquad
      -\Delta_{\Cyl}v=v^{-\gamma}\sigma .
\end{equation}
Let \(\zeta\) denote the \(t\)-marginal of \(\sigma\). Then
\[
 \mu(B_R\setminus\{0\})=\sph R^{-b(\gamma-1)}
\int e^{b(\gamma-1)t}\,d\zeta(t),
\]
where both sides may be infinite.

Indeed, $\phi^{1-\gamma}(x)=|x|^{b(\gamma-1)}$, so for a
nonnegative Borel test $F(t)$, the definition of the pushforward gives
\[
 \int_\Cyl F(t)\,d\sigma
  =\sph^{-1}\int_{B_R\setminus\{0\}}
       F(\log(R/|x|))|x|^{b(\gamma-1)}\,d\mu(x).
\]
Substituting $|x|=Re^{-t}$ and then taking $F(t)=e^{b(\gamma-1)t}$
recovers the stated mass formula by monotone approximation, even when
one of the integrals is infinite.

For each \(h\in\HH\), the \(L^2(d\omega)\) trace is continuous in \(t\) and
satisfies
\begin{equation}\label{time-trace}
 \|h(t,\cdot)\|_2\le\sqrt{t}\,\|h\|,\qquad
 \|h(t,\cdot)-h(s,\cdot)\|_2\le\sqrt{|t-s|}\,\|h\|.
\end{equation}
For smooth $h$ the fundamental theorem of calculus holds in
$L^2(d\omega)$. Consequently
\[
 \begin{aligned}
 \|h(t,\cdot)\|_2^2
 &\le t\int_0^t\|\partial_rh(r,\cdot)\|_2^2\,dr
  \le t\|h\|^2,\\
 \|h(t,\cdot)-h(s,\cdot)\|_2^2
 &\le |t-s|\int_{\min(s,t)}^{\max(s,t)}
              \|\partial_rh(r,\cdot)\|_2^2\,dr.
 \end{aligned}
\]
These bounds persist in the completion, yield a continuous
$L^2(d\omega)$ trace, and show that its value at $t=0$ is zero. The radial
subspace is
\[
 \Hrad=\{V(t)=\int_0^t f(s)\,ds: f\in L^2(0,\infty)\}.
\]
Here the radial norm is $\|V\|^2=\int_0^\infty|V'|^2dt$.
To verify density of radial test functions, choose
$f_j\in C_c^\infty(0,\infty)$ with $f_j\to f$ in $L^2$ and set
$c_j=\int_0^\infty f_j$. For each $j$, select a smooth function
$\eta_j$ of integral one, supported beyond $\supp f_j$, with
$\|\eta_j\|_2\le (j(1+|c_j|))^{-1}$. This is possible by spreading a
fixed bump over a sufficiently long interval. Then
\[
 V_j(t)=\int_0^t(f_j(r)-c_j\eta_j(r))\,dr
\]
is smooth, vanishes both near $0$ and for large $t$, and satisfies
$\|V_j-V\|=\|f_j-c_j\eta_j-f\|_2\to0$.

For a compact set, the energy capacity \(\ecap(K)\) is the infimum of
\(\|h\|^2\) over smooth functions with compact support that are at least one
in a neighborhood of \(K\). We extend the capacity to general sets in the
usual way. A measure is called diffuse if it vanishes on sets of zero
capacity. A property that holds outside such a set is said to hold
quasi-everywhere (q.e.), and we use the quasi-continuous representatives of
energy functions. The form is transient since \eqref{time-trace} bounds
\(\int e^{-t}|h|\,dm\) by \(C\|h\|\). We use two facts about the extended
space, see \cite[Secs.~1.5,
2.1--2.2]{FOT}. A strongly convergent sequence
in \(\HH\) has a subsequence that converges q.e. Moreover, a diffuse Radon
measure can be exhausted by increasing restrictions to compact sets on which
it has finite energy integral. For a measure to have finite energy integral
means that
\[
 \int|h|\,d\rho\le C\|h\|
\]
for all \(h\in\HH\).

We say that \(\sigma\) is \emph{admissible} if there exists
\(v=v_\sigma\in\HH\), positive q.e., such that
\[
 \rho=v^{-\gamma}\sigma,
 \qquad \B(v,h)=\int h\,d\rho
 \quad\text{for all }h\in\HH,
\]
and \(\rho\) has finite energy integral. For the zero source, we take both the solution and the reaction
to be zero. We write
\[
 E(\sigma)=\|v_\sigma\|^2.
\]

\begin{proposition}[Variational facts]\label{form}
Let \(\sigma\ne0\) be a positive diffuse measure and set \(q=1-\gamma\).
Then the source is admissible if and only if the corresponding quantity
below is finite.
\[
 \begin{cases}
 S_q(\sigma)=\displaystyle\sup_{\|h\|=1}\int|h|^q\,d\sigma,
     & 0<\gamma<1,\\[3pt]
 \sigma(\Cyl), & \gamma=1,\\[3pt]
 S_q(\sigma)=\displaystyle\inf_{h\ge0,\ \|h\|=1}\int h^q\,d\sigma,
     & \gamma>1 ,
 \end{cases}
\]
Here we use the convention \(0^q=+\infty\) when \(q<0\). The solution is
unique and increases with the datum. Moreover,
\begin{equation}\label{extreme}
 E(\sigma)=\int v_\sigma^{\,1-\gamma}\,d\sigma,\qquad
 E(\sigma)=S_q(\sigma)^{2/(1+\gamma)}\quad(\gamma\ne1).
\end{equation}
If \(\sigma_j\uparrow\sigma\) and \(\sigma\) is admissible, then
\(v_{\sigma_j}\to v_\sigma\) strongly in \(\HH\).
\end{proposition}

\begin{proof}
Let \(q\ne0\). We minimize
\[
 J_\sigma(h)=\|h\|^2/2-q^{-1}\int h^q\,d\sigma
\]
over the nonnegative cone. If \(0<q<1\), the trace bound gives coercivity,
and norm continuity follows from \(|z^q-w^q|\le|z-w|^q\). If \(q<0\), there
is a competitor with finite energy functional, and norm lower semicontinuity
follows from Fatou's lemma along q.e.\ convergent subsequences. By
convexity, the functional is also weakly lower semicontinuous. The direct
method gives a minimizer. For \(q>0\), a sufficiently small multiple of a
nonnegative test function with positive source integral gives a negative
value of the functional. For \(q<0\), the functional is infinite at zero.
Thus the minimizer \(v\) is nonzero in both cases.

Using nonnegative variations and applying the monotone convergence theorem
to the difference quotients, we obtain
\[
 \int v^{-\gamma}h\,d\sigma\le\B(v,h)\qquad(h\ge0).
\]
Thus \(v\) is a nonnegative weak supersolution of the homogeneous Laplace
equation on the cylinder. The weak Harnack inequality on overlapping balls
gives a positive lower bound for \(v\) on each compact subcylinder. The
bound holds q.e. by taking the Sobolev representative of the positive part.
We can now use variations of both signs with smooth compact support. The
inequality for nonnegative variations bounds the reaction on \(\HH\), so the
weak equation follows by density. To prove comparison and uniqueness, we
subtract the two equations and test with the positive part of the
difference.

Taking \(v\) as a test function gives the first equality in \eqref{extreme}.
If \(0<q<1\), we apply H\"older's inequality with \(d\sigma=v^\gamma d\rho\)
to obtain
\[
 \int |h|^q\,d\sigma\le
 (\int |h|\,d\rho)^q(\int v\,d\rho)^\gamma
 \le \|h\|^q\|v\|^{1+\gamma}.
\]
Equality holds for the normalized solution. If \(q<0\), the solution itself
is a competitor with finite functional. Minimizing over \(ch\) with \(c>0\)
and \(\|h\|=1\) gives the minimum \((1/2-1/q)S_q(\sigma)^{2/(2-q)}\).
Evaluating the functional at \(v\) gives the second equality in
\eqref{extreme}.

For \(\gamma=1\), testing with \(v\) shows that the source must have finite
mass. To prove existence, we first consider a restriction of the source with
compact support and finite energy integral and minimize
\(\|h\|^2/2-\int\log h\,d\sigma\). The function \(\min(t,1)\) has finite
functional. We write the logarithmic term using \(\int(h-\log h)\,d\sigma\)
and subtract the continuous linear integral of \(h\). Fatou's lemma then
gives lower semicontinuity, and the linear bound gives coercivity. We use
the same variations as above. For ordered solutions \(v_i\le v_j\), we test
the equation for \(v_i\) with \(v_j-v_i\ge0\) and obtain
\begin{equation}\label{order-norm}
 \|v_j-v_i\|^2\le\|v_j\|^2-\|v_i\|^2 .
\end{equation}
At $\gamma=1$, use the finite-energy-integral exhaustion theorem for
diffuse measures recalled above to select compact sets
$F_1\subset F_2\subset\cdots\Subset\Cyl$ such that
$\sigma(F_j)\uparrow\sigma(\Cyl)$ and each restriction
$\sigma_j=\mathbf1_{F_j}\sigma$ has finite energy integral.
After deleting any initial zero restrictions, the compact-source
logarithmic minimization gives $v_j=v_{\sigma_j}$. Comparison for
the inverse-power equation implies $v_i\le v_j$ quasi-everywhere for
$i\le j$. Testing the equation for $v_j$ by $v_j$ and using
\eqref{order-norm} yields
\[
 \|v_j\|^2=\sigma_j(\Cyl),\qquad
 \|v_j-v_i\|^2\le\sigma_j(\Cyl)-\sigma_i(\Cyl)
       \quad(j\ge i).
\]
The right-hand side tends to zero uniformly for $j\ge i$ as
$i\to\infty$ because the masses increase to the finite mass of
$\sigma$. Thus $v_j$ converges strongly in $\HH$; the argument
requires only the compact logarithmic minimizers, not the existence of
a global logarithmic minimum. Indeed,
$\int(-\log h)_+\,d\sigma$ may be infinite for every
$h\in\HH$ if the source accumulates sufficiently rapidly near
$t=0$. Starting with the first
nonzero approximant, comparison gives a common positive lower bound on each
compact set. By dominated convergence with respect to the limiting source,
we can pass to the limit in the equation for smooth compactly supported
tests. To extend the equation to \(\HH\), we first approximate the absolute
values of smooth test functions to get
\[
 \int|h|\,d\rho\le\|v\|\|h\|
\]
and
then use density, a q.e.\ convergent subsequence and Fatou's lemma. For any
increasing exhaustion of an admissible source, the finiteness conditions
hold for each submeasure. Comparison bounds the energies by the energy of
the limiting solution. The same Cauchy estimate and passage to the limit on
compactly supported tests prove the last assertion.
\end{proof}

\begin{lemma}[Changing the source]\label{change}
Finite sums and submeasures of admissible sources are admissible. If
\(\delta=|\sigma-\tau|\), then
\begin{equation}\label{change-bound}
 \|v_\sigma-v_\tau\|^2\le
 \begin{cases}E(\delta),&0<\gamma\le1,\\2E(\delta),&\gamma>1.\end{cases}
\end{equation}
Moreover,
\[
 v_{\sigma+\tau}\le v_\sigma+v_\tau
 \quad\text{q.e.}
\]
\end{lemma}
\begin{proof}
We first apply the finiteness conditions in Proposition~\ref{form}. For a
negative power, the sum of two competitors with finite functional has finite
integral against both sources. The sum of the two solutions is a
supersolution for the sum of the sources. We obtain the bound by testing
with the positive part of the difference between the solution for the sum
and this supersolution. For \(\gamma\le1\), we test the difference equation
with \(D=v_\sigma-v_\tau\). The contribution from the common part of the
sources is nonpositive. For the excess of \(\sigma\), the positive part of
the integrand is bounded by \(|D|^{1-\gamma}\) since \(D_+\le v_\sigma\). We
use the same bound for the excess of \(\tau\) and then apply
\eqref{extreme}. When \(\gamma=1\), each integrand is bounded above by one.
See \cite[Theorem~4.7]{Hara23} for this estimate. For $\gamma>1$, use the lattice supremum
\[
 \sigma\vee\tau
   =\tfrac12\bigl(\sigma+\tau+|\sigma-\tau|\bigr).
\]
It is a positive diffuse measure, because it is dominated by
$\sigma+\tau$; by the first part of the lemma it is admissible. Moreover,
$0\le(\sigma\vee\tau)-\sigma\le\delta$ and likewise with $\tau$.
Set $U=v_{\sigma\vee\tau}$. Testing the difference equation for
the ordered solutions gives
\[
 \|U-v_\sigma\|^2\le
 \int U^{1-\gamma}\,d((\sigma\vee\tau)-\sigma).
\]
We add the corresponding estimate for \(\tau\). Since \(U\ge v_\delta\), the
sum of the two squared norms is bounded by \(E(\delta)\). The triangle
inequality then gives \eqref{change-bound}.
\end{proof}

\subsection{Dyadic energy and compactness}
Let \(\zeta\) be the radial marginal of \(\sigma\). We set
\begin{equation}\label{gauge}
 I_k=[2^k,2^{k+1}),\qquad
 c_k=2^{k(1-\gamma)/2}\zeta(I_k),\qquad
 D_\gamma(\zeta)=\sum_{k\in\mathbb Z} c_k^{2/(1+\gamma)} .
\end{equation}
The intervals are half-open, so each spherical atom is counted once. By
vague convergence we mean convergence against \(C_c\) test functions. For
finite measures, narrow convergence means convergence against bounded
continuous test functions.

\begin{theorem}[Dyadic compactness]\label{angular}
Suppose
\[
 \begin{gathered}
 d\sigma=a(t,\theta)\,d\zeta(t)d\omega(\theta),\qquad
 0\le a\le L,\\
 \int a(t,\theta)\,d\omega(\theta)=1
 \quad\text{for }\zeta\text{-a.e. }t.
 \end{gathered}
\]
Then \(\sigma\) is admissible if and only if
\(D_\gamma(\zeta)<\infty\). In this case
\[
 E(\sigma)\asymp_{\gamma,L}D_\gamma(\zeta).
\]
At \(\gamma=1\), the exact normalization is
\[
 E(\sigma)=D_1(\zeta)=\sigma(\Cyl).
\]

Let \(\sigma_n\) be admissible sources with the same density bound and
suppose that they converge vaguely to \(\sigma\). The limit satisfies the
same bound. It is admissible and \(v_{\sigma_n}\to v_\sigma\) strongly in
\(\HH\) if and only if
\begin{equation}\label{tight}
 \lim_{K\to\infty}\sup_n
 \sum_{|k|>K}\big(2^{k(1-\gamma)/2}\zeta_n(I_k)\big)^{2/(1+\gamma)}
 =0 .
\end{equation}
Taking \(a=1\) gives both criteria for the one-dimensional equation
\(-V''=V^{-\gamma}\zeta\) in \(\Hrad\).
\end{theorem}

\begin{proof}
By \eqref{time-trace}, the restriction of the source to a compact radial
interval has finite energy integral. Hence the source is diffuse. To prove
the upper bound for the energy, set \(q=1-\gamma\) and \(r=2/(1+\gamma)\).
For \(h\in\HH\), let
\[
 \begin{aligned}
 A_j&=\|\partial_t h\|_{L^2(I_j\times\mathbb S^{N-1})},\\
 B_k&=\sum_{j\le k}2^{-(k-j)/2}A_j.
 \end{aligned}
\]
Integrating the derivative in the
radial variable gives
\[
 \|h(t,\cdot)\|_2\le2^{k/2}B_k
\]
for \(t\in I_k\).
The one-sided convolution weights are
$a_m=2^{-m/2}$ for $m\ge0$, and $a_m=0$ otherwise. Their $\ell^1$ norm
is $(1-2^{-1/2})^{-1}$. Young's discrete convolution inequality therefore
gives the explicit bound
\[
 \sum_{k\in\mathbb Z}B_k^2
 \le (1-2^{-1/2})^{-2}\sum_{j\in\mathbb Z}A_j^2
 \le (1-2^{-1/2})^{-2}\|h\|^2.
\]
This convolution constant is independent of $\gamma$, $L$, and the
source.
If \(q>0\), we apply H\"older's inequality
first on the sphere and then to the sum over the dyadic intervals to obtain
\[
 \int |h|^q\,d\sigma\le L\sum_k c_k B_k^q
 \le L D_\gamma(\zeta)^{1/r}(\sum_k B_k^2)^{q/2}.
\]
The bound follows from Proposition~\ref{form}. For \(q<0\), choose
\(0\le\psi\in C_c^\infty((1/2,4))\) equal to one on \([1,2]\). Set
\[
 h(t)=\sum_k c_k^{1/(1+\gamma)}2^{k/2}\psi(2^{-k}t).
\]
The summands have
bounded overlap, so the energy is at most \(C D_\gamma\). On each interval
of positive mass, the corresponding summand bounds the profile below. Hence
\[
 \int h^q\,d\sigma\le\sum_{c_k>0}c_k^{1+q/(1+\gamma)}=D_\gamma.
\]
Minimizing the functional gives the energy bound. For \(q=0\), we use the
mass identity.

To prove the reverse inequality, let
\[
 M(t)=\int v_\sigma(t,\theta)\,d\omega.
\]
Then \(\|M'\|_2\le\|v_\sigma\|\).
By approximation in energy, \eqref{time-trace} and diffuseness, the traces
agree with the quasi-continuous representatives almost everywhere with
respect to \(\zeta\,d\omega\). We use radial test functions and apply
Jensen's inequality with respect to \(a\,d\omega\) to obtain
\[
 -M''\ge (LM)^{-\gamma}\zeta.
\]
The distribution $-M''$ is positive, hence $M$ has a concave
representative and its one-sided derivatives are nonincreasing. A negative
slope at a finite height would force $M$ eventually to become negative,
contrary to $M\ge0$ on the whole half-line. Thus $M'\ge0$ a.e. and,
because $M'\in L^2(0,\infty)$ and is nonincreasing, $M'(t)\to0$ as
$t\to\infty$. We define \(A_j,B_k\) as
above, with \(M'\) in place of \(\partial_t h\). For the marginal reaction $\nu=-M''\ge0$ and our half-open
intervals, monotonicity of $M'$ gives
\[
 \nu(I_k)\le M'_-(2^k)
 \le \frac{1}{2^{k-1}}\int_{2^{k-1}}^{2^k}M'(t)\,dt
 \le 2^{-(k-1)/2}A_{k-1}.
\]
The middle inequality uses that $M'$ is nonincreasing and nonnegative;
the last is Cauchy--Schwarz. Also $M(t)\le2^{k/2}B_k$ for $t\in I_k$.
Since $\nu\ge(LM)^{-\gamma}\zeta$, it follows that
\[
 \begin{aligned}
 c_k&=2^{k(1-\gamma)/2}\zeta(I_k)\\
 &\le L^\gamma 2^{k(1-\gamma)/2}
     (2^{k/2}B_k)^\gamma\nu(I_k)
 \le\sqrt2\,L^\gamma A_{k-1}B_k^\gamma.
 \end{aligned}
\]
Using
the upper bound for \(L^\gamma M^\gamma\), we obtain, for every set of
indices \(J\),
\begin{equation}\label{localized}
 \sum_{k\in J}c_k^r\le
 2^{r/2}L^{\gamma r}
 \left(\sum_{k\in J}A_{k-1}^2\right)^{1/(1+\gamma)}
 \left(\sum_{k\in J}B_k^2\right)^{\gamma/(1+\gamma)} .
\end{equation}
Apply H\"older's inequality to the last estimate for $c_k$ on an
arbitrary set of indices $J$; this gives exactly \eqref{localized}, since
$r=2/(1+\gamma)$. For $J=\mathbb Z$, the convolution bound and
$\sum A_j^2\le\|M'\|_2^2$ imply
\[
 D_\gamma(\zeta)
 \le C_{\gamma,L}\|M'\|_2^2
 \le C_{\gamma,L}E(\sigma),
\]
which is the reverse estimate.

For necessity of the uniform tail condition, strong energy convergence
implies $M_n'\to M'$ in $L^2(0,\infty)$. The reverse triangle inequality
on each $I_j$ gives
\[
 \bigl|\|M_n'\|_{L^2(I_j)}-\|M'\|_{L^2(I_j)}\bigr|
 \le\|M_n'-M'\|_{L^2(I_j)}.
\]
Consequently the sequences $(A_{j,n})_j$ converge in $\ell^2$ to
$(A_j)_j$; convolution by $(2^{-m/2})_{m\ge0}$ gives the same conclusion
for $(B_{k,n})_k$. A convergent family in $\ell^2$ has uniformly vanishing
tails. Taking $J=\{|k|>K\}$ in \eqref{localized} proves
\eqref{tight}.

For sufficiency we first consider sources restricted to a compact interval
\(J=[A,B]\Subset(0,\infty)\) on which the restrictions converge narrowly.
The inequality
\[
 \sigma_n\le L\zeta_n\otimes\omega
\]
passes to the limit by testing with continuous nonnegative functions. Locally, the same argument
applies under vague convergence. On $J\times\mathbb S^{N-1}$ the energy unit ball is bounded in
$H^1$, and the global trace estimate gives a uniform local $L^2$ bound.
Rellich compactness gives convergence in $L^2(J\times\mathbb S^{N-1})$
along a subsequence. To see explicitly why the convergence improves to
uniform convergence in $t$ with values in $L^2(d\omega)$, let
$h_n$ be a bounded sequence in the energy unit ball converging to
$h$ in this local $L^2$ space. Choose for each $t\in J$ a one-sided
or centered interval $J_{t,\delta}\subset J$ of length $\delta$,
whose endpoints may depend on $t$. Denoting the average over this
interval by $\langle\cdot\rangle_{t,\delta}$, the
trace estimate gives
\[
 \|h_n(t)-h(t)\|_2
 \le 2\sqrt{\delta}\,C
    +\|\langle h_n-h\rangle_{t,\delta}\|_2
 \le 2C\sqrt{\delta}
    +\delta^{-1/2}\|h_n-h\|_{L^2(J\times\mathbb S^{N-1})}.
\]
Here $C$ is independent of $n$, $t$, and $\delta$; the same trace
bound applies to the weak limit $h$. First choose $\delta$ small and
then $n$ large. This proves convergence in
$C(J;L^2(d\omega))$ and the asserted relative compactness.
Write $\mathcal L_n(F)=\int F\,d\sigma_n$ for
$F\in C(J;L^1(d\omega))$. The domination $\sigma_n\le L\zeta_n\otimes\omega$
and the uniform boundedness of $\zeta_n(J)$ give
\[
 |\mathcal L_n(F)|\le L\sup_n\zeta_n(J)
             \sup_{t\in J}\|F(t,\cdot)\|_{L^1(d\omega)}.
\]
Continuous functions of $(t,\theta)$ are dense in this Banach space.
Narrow convergence first proves convergence of $\mathcal L_n(F)$ on this
dense class, and the displayed uniform bound extends it to every $F$.
For clarity, if $\mathcal K\subset C(J;L^1(d\omega))$ is
compact, choose a finite $\epsilon$-net
$F_1,\ldots,F_M\in\mathcal K$. The uniform norm bound on
$\mathcal L_n$ and on $\mathcal L$ implies
\[
 \sup_{F\in\mathcal K}|\mathcal L_n(F)-\mathcal L(F)|
 \le \max_{1\le j\le M}
       |\mathcal L_n(F_j)-\mathcal L(F_j)|+2C\epsilon.
\]
For each fixed $j$ the first term tends to zero; sending
$\epsilon\downarrow0$ proves the claimed uniform convergence.

If the limit measure is nonzero, its solution \(v\) is bounded below by a
positive constant on \(J\times\mathbb S^{N-1}\). Extending the inverse power
below this bound as a Lipschitz function, we get
\(v^{-\gamma}\in C(J;L^2(d\omega))\). Hence the products \(v^{-\gamma}h\)
with \(\|h\|\le1\) form a relatively compact family of test functions in
\(C(J;L^1(d\omega))\). We subtract the two equations and use the
monotonicity of the inverse power to obtain
\[
 \|v_{\sigma_n}-v\|
 \le\sup_{\|h\|\le1}
      \left|\int v^{-\gamma}h\,d(\sigma_n-\sigma)\right|\longrightarrow0 .
\]
If the limit measure is zero, the masses tend to zero, and the energy bound
on the compact interval gives the conclusion.

Condition \eqref{tight} bounds the two dyadic tails uniformly. On
every fixed compact radial interval, vague convergence and the density
bound give a uniform mass bound; therefore the finitely many remaining
$c_k$ are bounded and the upper energy estimate yields
$\sup_n E(\sigma_n)<\infty$. To prove admissibility of the limit without
assuming continuity at dyadic endpoints, exhaust the cylinder by compact radial intervals whose
endpoints are continuity points of the limiting marginal. On each such
interval the restricted energies converge by the preceding argument.
Comparison bounds them by the global energy bound. The reverse dyadic
estimate and monotone convergence in \eqref{gauge} imply that the limit is
admissible. We choose the endpoints in \((2^{-K},2^{-K+1})\) and
\((2^K,2^{K+1})\). The restrictions to these intervals converge narrowly because neither
endpoint is an atom of the limiting radial marginal. Let
$\sigma_n^{(K)}$ and $\sigma^{(K)}$ denote the restrictions. The source
change estimate and the upper energy bound give, for all $K$,
\[
 \begin{aligned}
 \|v_{\sigma_n}-v_{\sigma_n^{(K)}}\|^2
 &\le C_\gamma E(\sigma_n-\sigma_n^{(K)})
 \le C_{\gamma,L}
       \sum_{|k|\ge K-1}c_{k,n}^{2/(1+\gamma)},\\
 \|v_\sigma-v_{\sigma^{(K)}}\|^2
 &\le C_{\gamma,L}
       \sum_{|k|\ge K-1}c_k^{2/(1+\gamma)}.
 \end{aligned}
\]
The first bound is uniform in $n$ by the full two-ended condition
\eqref{tight}; the second tends to zero because the limit is
admissible. This is precisely where the dyadic tail condition upgrades
local convergence to whole-cylinder convergence. For fixed $K$,
the compact-interval argument yields
$\|v_{\sigma_n^{(K)}}-v_{\sigma^{(K)}}\|\to0$. The triangle inequality,
followed by $n\to\infty$ and then $K\to\infty$, proves sufficiency.

The uniform bound on the angular conditional densities is used
explicitly in the reverse estimate \eqref{localized} and in the
compact-interval convergence argument. The unrestricted admissibility
implications are considered in Proposition~\ref{prop:angular}.

For radial data, invariance under rotations and uniqueness imply that the
cylinder solution is radial. Conversely, a positive distributional solution
in \(\Hrad\) is concave and nondecreasing, and its slope tends to zero at
infinity. With \(\rho=V^{-\gamma}\zeta\), integration gives
\begin{equation}\label{green}
 V'_+(t)=\rho((t,\infty)),\qquad V(t)=\int\min(t,s)\,d\rho(s).
\end{equation}
By Tonelli's theorem,
\[
 \int |h|\,d\rho\le\int |h'| V'\,dt
\]
for radial test
functions. Applying this estimate to the spherical mean of the absolute
value of a test function, we extend the radial equation to \(\HH\).
\end{proof}

\begin{suppprop}[One-sided dyadic admissibility without angular control]\label{prop:angular}
If $\gamma>1$ and $D_\gamma(\zeta)<\infty$, then $\sigma$ is admissible and
\[
 E(\sigma)\leq C_\gamma D_\gamma(\zeta).
\]
At $\gamma=1$, every finite positive diffuse measure is admissible, with $E(\sigma)=\sigma(\Cyl)$,
regardless of its angular distribution. For $0<\gamma<1$, finiteness of $D_\gamma(\zeta)$ alone does
not ensure admissibility. For $\gamma>1$, admissibility alone does not imply finiteness of
$D_\gamma(\zeta)$.
\end{suppprop}
\begin{proof}
For $\gamma>1$ put $q=1-\gamma<0$. Choose a nonnegative $\psi\in C_c^\infty((1/2,4))$ equal to one on $[1,2]$. Define the radial competitor
\[
 w(t)=\sum_{c_k>0}2^{k/2}c_k^{1/(1+\gamma)}\psi(2^{-k}t).
\]
The supports have uniformly bounded overlap. For a finite partial sum, the integral of the squared
derivative is at most a constant times the sum of $c_k^{2/(1+\gamma)}$ over its indices. Applying
the same bound to differences of partial sums proves convergence in $\HH$; local finiteness
identifies the limit with the displayed function. For $t\in I_k$ with $c_k>0$ we have $w(t)\geq
2^{k/2}c_k^{1/(1+\gamma)}$. Since $w$ is radial and $q<0$,
\begin{align*}
 \int_{\Cyl} w^q\,d\sigma
 &=\int_{0}^{\infty}w^q\,d\zeta\\
 &\leq\sum_{c_k>0} c_k^{1+q/(1+\gamma)}=D_\gamma(\zeta).
\end{align*}
Intervals with zero marginal mass are irrelevant to this integral. Proposition~\ref{form} gives admissibility and, for a unit-energy test obtained by normalizing $w$,
\[
 E(\sigma)\leq
 \left(\|w\|^{-q}\int_{\Cyl} w^q\,d\sigma\right)^{2/(1+\gamma)}
 \leq C_\gamma D_\gamma(\zeta).
\]
The zero-source case is immediate. At $\gamma=1$, the assertion is the mass identity and admissibility criterion in Proposition~\ref{form}.

For $0<\gamma<1$, put $q=1-\gamma>0$. Fix a relatively compact interior coordinate ball strictly
inside one dyadic slab. In it choose a sequence of pairwise disjoint smaller balls $B_j$ of radii
$r_j\downarrow0$. There are smooth nonnegative cutoffs $\eta_j$ equal to one on $B_j$, supported in
slightly larger balls and satisfying $\|\eta_j\|^2\leq C r_j^{N-2}$. Let $\sigma_j$ be a smooth
nonnegative density of mass $2^{-j}$ supported in $B_j$, and put $\sigma=\sum_j\sigma_j$. This is a
finite absolutely continuous measure, hence diffuse, supported in a single slab, and its dyadic sum
is finite. However,
\[
 \sup_{\|h\|=1}\int |h|^q\,d\sigma
 \geq 2^{-j}\|\eta_j\|^{-q}
 \geq C^{-q/2}2^{-j}r_j^{-(N-2)q/2}.
\]
Choose $r_j$ small enough that the last quantity exceeds $j$. Then the trace supremum in Proposition~\ref{form} is infinite.

For $\gamma>1$, select one coordinate ball strictly within each $I_k\times\mathbb S^{N-1}$,
$k\leq-1$, and prescribe a smooth source there of mass $m_k=2^{k(\gamma-1)/2}$. The resulting
measure has finite total mass but $c_k=1$ for every $k\leq-1$. Thus $D_\gamma=\infty$. Let
$A_k=(2^{-k}m_k)^{1/(\gamma-1)}$ and choose cutoffs $\eta_k$ equal to one on the supports of the
sources, with mutually disjoint supports in their respective slabs. Shrink the balls so that
$A_k^2\|\eta_k\|^2\leq 2^k$, possible since $N\geq3$. The function $w=\sum_{k\leq-1} A_k\eta_k$
belongs to $\HH$ and
\[
 \int w^{1-\gamma}\,d\sigma
 =\sum_{k\leq-1}m_k A_k^{1-\gamma}=\sum_{k\leq-1}2^k<\infty.
\]
The negative-power criterion in Proposition~\ref{form} proves admissibility, completing the proof.
\end{proof}

The uniform angular bound controls information not contained in the radial marginal.
Proposition~\ref{prop:angular} separates the one-sided statements that survive without it from the
two-sided and sequential statements of Theorem~\ref{angular}. In particular, no extension of the
dyadic strong-convergence criterion to arbitrary diffuse data is claimed. For $\gamma=1$ the
relevant general compactness criterion is the condenser condition of Theorem~\ref{condenser}.

\begin{remark}
Suppose the transformed radial marginal $\zeta$ has a density comparable
near $t=0$ to
\[
 t^\kappa(\log(e/t))^{-\eta}.
\]
The dyadic criterion holds if $\kappa>(\gamma-3)/2$, or, at the
borderline exponent, if $\eta>(1+\gamma)/2$. If the density is
comparable near infinity to
\[
 t^{-1-\delta}(\log t)^{-\eta},
\]
then the corresponding criterion is $\delta>(1-\gamma)/2$, with
the same logarithmic restriction at equality. These statements follow by integrating over the dyadic intervals.
If the density of \(\zeta\) (not of the original measure \(\mu\)) is
bounded above and below by positive constants near \(t=0\), the energy condition is \(\gamma<3\). This agrees with the boundary
threshold for the usual singular Dirichlet problem~\cite{SunZhang}. In
particular, a finite radial measure in the physical variables, supported
away from the outer boundary and with no mass at the pole, is admissible for
\(\gamma\ge1\). This follows from \eqref{weight} and the competitor
\(\min(t,1)\). For \(\gamma<1\), the transformed measure
\(\sum_{k\ge1}\delta_{2^k}\) has finite physical mass and infinite
\(D_\gamma\).
\end{remark}

The finiteness of $D_\gamma$ does not depend on selecting the dyadic
base $2$. For any fixed base $a>1$, define intervals $[a^k,a^{k+1})$
and the corresponding sum with weight $a^{k(1-\gamma)/2}$. Each
$a$-adic interval meets only a bounded number of dyadic intervals and
vice versa, with ratios of their endpoints bounded in terms of $a$.
Finite-dimensional inequalities for the exponent $2/(1+\gamma)>0$ compare
the two sums up to constants depending on $a$ and $\gamma$, even when the
marginal has atoms.

\subsection{Sharp endpoint moments}
For a radial solution, let \(L_\infty=V(\infty)\) be the pole trace and
\(L_0=V'(0+)\) the slope at the boundary. These quantities may be infinite.
In the original variables,
\[
 r^bu(r)=V(\log(R/r)),\qquad
 u_r(r)=-r^{-b-1}\bigl(bV(t)+V'(t)\bigr),\quad t=\log(R/r).
\]
Therefore $L_\infty=\lim_{r\downarrow0}r^bu(r)$ and the outer normal
derivative is $u_r(R^-)=-R^{-b-1}L_0$ (the positive boundary flux has
magnitude $R^{-b-1}L_0$).

\begin{theorem}[Endpoint traces]\label{endpoint}
Let \(V\) be the radial solution for an admissible source \(\zeta\) and set
\[
 p=1+\gamma,\qquad
 M_\infty=\int s\,d\zeta,\qquad
 M_0=\int s^{-\gamma}\,d\zeta.
\]
Then
\begin{equation}\label{two-moments}
 M_i\le L_i^p\le pM_i,\qquad i\in\{0,\infty\}.
\end{equation}
All quantities in \eqref{two-moments} take values in $[0,+\infty]$:
when $L_i=+\infty$, the upper inequality forces $M_i=+\infty$,
and the lower inequality is understood in the extended order.
The lower constant is attained by a single point source; the upper
constant is approached by the family below, but cannot be attained
by a nonzero finite trace. For every \(A>0\), the integrals
\(\int_{(0,A]}s\,d\zeta\) and \(\int_{[A,\infty)}s^{-\gamma}\,d\zeta\) are
finite.

Suppose that \(V_n\to V\) locally uniformly and that the pole traces of the
approximants are finite. The pole traces converge to a finite pole trace if
and only if
\begin{equation}\label{trace-tight}
 \lim_{A\to\infty}\sup_n\int_{(A,\infty)}s\,d\zeta_n(s)=0 .
\end{equation}
If the boundary slopes of the approximants are finite, the corresponding
necessary and sufficient condition for convergence to a finite boundary
slope is
\[
 \lim_{A\downarrow0}\sup_n\int_{(0,A)}s^{-\gamma}\,d\zeta_n(s)=0.
\]
\end{theorem}

\begin{proof}
By \eqref{green}, we have \(L_\infty=\int s\,d\rho\) and
\(d\zeta=V^\gamma d\rho\). Since \(V\le L_\infty\), we obtain the lower
bound. For the upper bound, integrate \(pV^\gamma V'\) from \(A\) to
infinity. By Tonelli's theorem and \(V(t)\le V(s)\) for \(t<s\), we get
\begin{equation}\label{trace-diff}
 L_\infty^p-V(A)^p\le p\int_{(A,\infty)}(s-A)\,d\zeta(s).
\end{equation}
Letting \(A\downarrow0\) gives the upper bound. The part of \(M_\infty\)
near zero is bounded by
\[
 V(A)^\gamma\int_{(0,A]}s\,d\rho\le V(A)^p.
\]
The inversion $\mathcal IV(t)=tV(1/t)$ is an isometry of $\Hrad$.
Indeed, with $W=\mathcal IV$ and $s=1/t$,
\[
 \int_0^\infty|W'(t)|^2dt
   =\int_0^\infty|V'(s)-V(s)/s|^2ds.
\]
Expanding the square gives
\[
 |V'-V/s|^2=|V'|^2-\frac{d}{ds}\bigl(V(s)^2/s\bigr).
\]
The boundary term vanishes at zero because
\[
 \frac{|V(s)|^2}{s}
 \le\int_0^s|V'(r)|^2dr\longrightarrow0\quad(s\downarrow0).
\]
At infinity, split $V(s)=V(A)+\int_A^sV'$ and use
\[
 \limsup_{s\to\infty}\frac{|V(s)|}{\sqrt s}
 \le\left(\int_A^\infty|V'(r)|^2dr\right)^{1/2};
\]
then let $A\to\infty$. It follows that $V(s)^2/s\to0$ there as well,
and the isometry follows. Approximation by smooth radial functions shows
that $W$ has the required zero trace at the left endpoint. Under
inversion the measures in \eqref{green} change from \(\rho\) to
\((s\mapsto1/s)_\#(s\rho)\) and from \(\zeta\) to
\((s\mapsto1/s)_\#(s^{1-\gamma}\zeta)\). If
$\widetilde\rho=(s\mapsto1/s)_\#(s\rho)$ and
$\widetilde\zeta=(s\mapsto1/s)_\#(s^{1-\gamma}\zeta)$, then
$-W''=\widetilde\rho$ and
$\widetilde\zeta=W^\gamma\widetilde\rho$. Moreover,
\[
 W(\infty)=\lim_{t\to\infty}tV(1/t)=V'(0+)=L_0,
 \qquad
 \int t\,d\widetilde\zeta(t)
   =\int s^{-\gamma}\,d\zeta(s)=M_0.
\]
Thus the pole moment inequality for the inverted pair is exactly the
boundary slope inequality for the original source.

The lower constant is attained by a single atom. For the upper constant,
choose \(0<\alpha<1\) and \(B>1\) and take a concave profile equal to \(t\)
on \([0,1]\), \(t^\alpha\) on \([1,B]\), and \(B^\alpha\) thereafter. The
source is given by \(d\zeta=V^\gamma d(-V'')\). Adding the contributions of
the two jumps in the slope and the density between them gives
\[
 M_\infty=(1-\alpha)+\frac{1-\alpha}{p}(B^{\alpha p}-1)
                       +\alpha B^{\alpha p}.
\]
The ratio \(L_\infty^p/M_\infty\) tends to \(p\) if we let first
\(B\to\infty\) and then \(\alpha\downarrow0\). The upper bound is strict for every nonzero source with finite pole
trace. Indeed, writing $d\zeta=V^\gamma d\rho$, Tonelli gives
\[
 L_\infty^p=p\int_{(0,\infty)}
           \left(\int_0^s V(t)^\gamma\,dt\right)d\rho(s),
\]
whereas $pM_\infty=p\int sV(s)^\gamma\,d\rho(s)$.
For $d\rho$-almost every $s>0$ carrying nonzero mass, the
monotonicity formula \eqref{green} implies
$V(t)<V(s)$ for $0<t<s$, so the first integral is strictly
smaller. This proves nonattainment without excluding point sources:
for a single atom it is the \emph{lower} bound that is attained.
Inversion gives the optimal constants for the boundary slope.

To prove convergence, we use \eqref{trace-diff} and the reverse estimate
\[
 \int_{(A,\infty)}(s-A)\,d\zeta_n
 \le L_{\infty,n}^{\gamma}(L_{\infty,n}-V_n(A)),
\]
To derive the reverse inequality, use $d\zeta=V^\gamma d\rho$,
$V(s)\le L_\infty$, and \eqref{green}:
\[
 \begin{aligned}
 \int_{(A,\infty)}(s-A)\,d\zeta(s)
 &\le L_\infty^\gamma
     \int_{(A,\infty)}(s-A)\,d\rho(s)\\
 &=L_\infty^\gamma\bigl(L_\infty-V(A)\bigr).
 \end{aligned}
\]
Assume first $L_{\infty,n}\to L_\infty<\infty$ and $V_n\to V$ locally
uniformly. Given $\epsilon>0$, choose $A$ large enough that
$L_\infty-V(A)<\epsilon$. For all sufficiently large $n$, the reverse
estimate makes $\int_{(A,\infty)}(s-A)d\zeta_n$ small uniformly in those
indices. Since $s\le2(s-A)$ on $(2A,\infty)$, their first-moment tails
are small. Each of the finitely many remaining $\zeta_n$ has a finite
first moment by \eqref{two-moments}, so enlarging $A$ also controls these
indices. This proves necessity of \eqref{trace-tight}.

Conversely, assume \eqref{trace-tight}. Fix $A$ large. The estimate
\eqref{trace-diff}, applied to every $V_n$, gives
\[
 0\le L_{\infty,n}^p-V_n(A)^p
 \le p\sup_n\int_{(A,\infty)}s\,d\zeta_n(s).
\]
The local convergence bounds $V_n(A)$ and the assumed tightness bounds
$L_{\infty,n}$. Passing to limits and then sending $A\to\infty$ shows
that the monotone limit $V(\infty)$ is finite and that
$L_{\infty,n}^p\to V(\infty)^p$, hence the traces converge. Applying the
same argument to $\mathcal IV_n$, using the transformed sources above,
gives the criterion for the boundary slopes.
\end{proof}

\begin{remark}
By \eqref{time-trace}, strong convergence in energy implies the local
convergence assumed in Theorem~\ref{endpoint}. The moment condition is still
needed when the limit source is nonzero. Take a radial profile \(V\) with
compactly supported reaction \(\rho\), and add \(c\min(t,s)/s\) with
\(c>0\). The new source is \(\rho+(c/s)\delta_s\) multiplied by the new
profile to the power \(\gamma\). As \(s\to\infty\), the sources converge
vaguely. The change in the energy norm is \(c/\sqrt{s}\) and the pole trace
increases by \(c\). Inversion gives an example for the boundary slope. The
weighted dilation \(\ell^{(\gamma-1)/2}(s\mapsto\ell s)_\#\zeta\) gives the
profile \(\sqrt\ell V(t/\ell)\), with the same energy. The pole trace and
the boundary slope are multiplied by \(\sqrt\ell\) and \(1/\sqrt\ell\),
respectively.
\end{remark}

\begin{remark}
Let \(\zeta,\eta\) be two nonzero admissible radial sources supported on
\([A,\infty)\) with \(A>0\), and set \(c=\min(V_\zeta(A),V_\eta(A))\) and
\(\Delta=|\zeta-\eta|\). We test the difference of the equations with the
difference of the profiles and use \(|h(s)|\le\sqrt{s}\|h'\|_2\) to obtain
\[
 \|(V_\zeta-V_\eta)'\|_2\le c^{-\gamma}\int\sqrt{s}\,d\Delta(s).
\]
We compare both solutions with the solution for \(\zeta\vee\eta\) and apply
\eqref{green}. This bounds the absolute value of their difference by
\[
 c^{-\gamma}\int\min(t,s)\,d\Delta(s).
\]
In both estimates the common part
of the sources has the favorable sign. For finite pole traces, the pointwise
bound gives local Lipschitz continuity with respect to the first moment of
\(\Delta\).
\end{remark}

\begin{remark}
Every strict cylinder solution yields a distributional solution of
\eqref{physical} on the whole ball, with zero mass at the pole. By
\eqref{time-trace}, the spherical mean of \(v\) is \(O(\sqrt t)\).
Consequently, \(|x|^{-2}u\) is locally integrable and
\[
 \eps^{-2}\int_{B_{2\eps}}u=O(\eps^b\sqrt{\log(R/\eps)}).
\]
The physical
reaction is finite near the pole since it is the pullback of
\(\sph\phi^{-1}\rho\) and a cutoff of \(e^{-bt}\) on the cylinder tail
belongs to \(\Hrad\). We first test away from the pole, using a cutoff equal to
zero on \(B_\eps\) and to one outside \(B_{2\eps}\). The preceding estimate
controls the terms containing derivatives of the cutoff. We pass to the
limit in the remaining terms by dominated convergence. For a nonzero radial
solution, monotonicity gives \(\sph\int V^2dt=\infty\) for the Hardy
integral.
\end{remark}
\section{Capacity tests for general sources}\label{sec-geometry}

\subsection{A condenser criterion for inverse power one}
Let \(\ecap(K,U)\) denote the infimum of \(\|h\|^2\) over
\(h\in C_c^\infty(U)\) that are at least one near \(K\). Here
\(K\subset U\Subset\Cyl\), with \(K\) compact and \(U\) open. This is the
usual condenser capacity~\cite{Mazya}. For finite diffuse measures
\(\lambda,\eta\), we define
\[
 A_L(\lambda,\eta)
 =\sup_{\ecap(K,U)\le L}
   \max\{0,\lambda(K)-\eta(U),\eta(K)-\lambda(U)\}.
\]
We compare the mass on the inner plate for one measure with the mass in the
outer open set for the other. The capacity bounds the energy needed to
separate the plates.

The next lemma isolates the metric comparison used in the theorem. It
also shows that the condition is pointwise in each fixed capacity
threshold $L$; no uniformity over all $L$ is needed.

\begin{condenserdistance}[Condenser discrepancy and test distance]
For finite diffuse measures $\lambda,\eta$ on $\Cyl$, set
\[
 d(\lambda,\eta)=
 \sup_{\substack{h\in C_c^\infty(\Cyl)\\
             \|h\|\le1,\;\|h\|_\infty\le1}}
  \left|\int h\,d(\lambda-\eta)\right|.
\]
For every $L>0$ and $0<\delta<1$,
\[
 \begin{split}
 A_L(\lambda,\eta)&\le\max(1,\sqrt L)\,d(\lambda,\eta),\\
 d(\lambda,\eta)&\le
   2A_{16/\delta^2}(\lambda,\eta)
       +2\delta\max\{\lambda(\Cyl),\eta(\Cyl)\}.
 \end{split}
\]
Thus on a bounded-mass family the condition
$A_L(\lambda_n,\lambda)\to0$ for each fixed $L>0$
is equivalent to $d(\lambda_n,\lambda)\to0$.
\end{condenserdistance}
\begin{proof}
The estimate extends to compactly supported energy test functions with the
same energy and supremum bounds. Indeed, we use smooth approximation with a
common supremum bound and take a q.e.\ convergent subsequence. If $\varphi$ is a nonnegative admissible condenser test, truncate it
to $h=\min(1,\varphi)$; then $h=1$ near $K$, $h=0$ outside $U$,
$0\le h\le1$, and $\|h\|^2\le L+o(1)$. Thus
\[
 \lambda(K)-\eta(U)
 \le\int h\,d(\lambda-\eta)
 \le\max(1,\sqrt{L+o(1)})\,d(\lambda,\eta).
\]
Interchanging $\lambda,\eta$ and sending the approximation error to
zero gives $A_L\le\max(1,\sqrt L)d$. The extension of $d$ from smooth
tests to $h$ follows by truncation, energy approximation, and a
quasi-everywhere convergent subsequence; diffuseness allows passage in
the measure integrals by bounded convergence. For the
converse, take a test function \(0\le h\le1\) and \(0<t<1-\delta\), and set
\(K=\{h\ge t+\delta\}\) and \(U=\{h>t\}\). The clipped level function
\[
 h_{t,\delta}(x)=\min\bigl(1,\max(0,(h(x)-t)/\delta)\bigr)
\]
is one on $K$, zero outside $U$, and has energy at most
$\delta^{-2}\|h\|^2\le\delta^{-2}$. Smooth approximation of its two
levels, keeping a safety margin, yields an admissible condenser test with
energy at most $16/\delta^2$. This gives
\[
 \lambda(\{h\ge t+\delta\})\le\eta(\{h>t\})+A_{16/\delta^2}.
\]
Integrating in $t\in(0,1-\delta)$ and using the layer-cake formula
compares $\int h\,d\lambda$ with $\int h\,d\eta$. The missing strips
of levels of total length at most $\delta$ contribute at most
$\delta\max\{\lambda(\Cyl),\eta(\Cyl)\}$. For signed tests we apply the
result to the positive levels of \(h\) and \(-h\). We obtain
\[
 d(\lambda,\eta)\le
 2A_{16/\delta^2}(\lambda,\eta)+
 2\delta\max\{\lambda(\Cyl),\eta(\Cyl)\}.
\]
Thus, on families of bounded mass, convergence in this distance is
equivalent to convergence of the condenser discrepancies.
\end{proof}

\begin{theorem}[Condenser criterion]\label{condenser}
Suppose \(\gamma=1\) and let \(\sigma_n,\sigma\) be finite positive diffuse
measures on \(\Cyl\). Then \(v_{\sigma_n}\to v_\sigma\) strongly in \(\HH\)
if and only if
\begin{equation}\label{two-plates}
 \sigma_n(\Cyl)\to\sigma(\Cyl),\qquad
 A_L(\sigma_n,\sigma)\to0\quad\text{for every }L>0 .
\end{equation}
Under vague convergence, the equivalent local formulation asks that,
for every $\epsilon>0$, some $0<a<b<\infty$ satisfy
\[
 \sup_n\sigma_n\bigl(((0,a)\cup(b,\infty))
           \times\mathbb S^{N-1}\bigr)<\epsilon,
\]
and that the condenser discrepancies tend to zero for condensers
contained in any fixed compact subcylinder. These are the precise
``tightness at both ends'' and local-discrepancy requirements.
\end{theorem}

\begin{proof}
By the condenser-discrepancy lemma, the mass and discrepancy conditions in
\eqref{two-plates} are equivalent, on bounded-mass families, to
$\sigma_n(\Cyl)\to\sigma(\Cyl)$ and
$d(\sigma_n,\sigma)\to0$.

Suppose that \(u_n=v_{\sigma_n}\to u=v_\sigma\) strongly in \(\HH\). Then
\(\rho_n\to\rho\) in the dual space \(\HH^*\) and \(\sigma_n=u_n\rho_n\).
For any \(f\in C_c^\infty(\Cyl)\) and any test function \(h\) in the
supremum, we have
\[
 \left|\int h\,d(\sigma_n-\sigma)\right|
 \le\|u_n-f\|\|u_n\|+\|fh\|\|u_n-u\|+\|f-u\|\|u\|.
\]
The product estimate
\[
 \|fh\|\le\|f\|_\infty+\|\nabla f\|_2
\]
is uniform
over the test functions. Approximating \(u\) by \(f\) gives convergence in
\(d\). The masses converge by \eqref{extreme}.

To prove sufficiency, fix a nonzero positive smooth compactly
supported measure $\beta_0$ and add $\varepsilon\beta_0$ to each
source and to the limit. This regularizing measure is unrelated to the
fractional exponent $\beta$ used later. The test distance is exactly unchanged by the common lift:
\[
 d(\sigma_n+\varepsilon\beta_0,
   \sigma+\varepsilon\beta_0)=d(\sigma_n,\sigma).
\]
The condenser discrepancy $A_L$ need not be exactly preserved:
$\beta_0(K)\le\beta_0(U)$ shows instead that it cannot increase
under a common positive lift. The mass bounds remain uniform in $n$.
At $\gamma=1$, the mass identity gives
\[
 \|v_{\sigma_n+\varepsilon\beta_0}\|^2
   =\sigma_n(\Cyl)+\varepsilon\beta_0(\Cyl),
\]
so the lifted sequence is bounded in energy. The mass identity bounds the new solutions \(w_n\) in
\(\HH\). On compact sets they have a positive lower bound by comparison with
\(v_{\varepsilon\beta_0}\). Up to a subsequence, \(w_n\rightharpoonup w\) in \(\HH\) and
\(w_n\to w\) locally in \(L^2\). Given \(\varphi\in C_c^\infty(\Cyl)\), we
extend the reciprocal below the local positive lower bound as a Lipschitz
function. For $K=\supp\varphi$ let $c_K>0$ be a common quasi-everywhere lower
bound for $w_n$ on a neighborhood of $K$. The chain rule gives
\[
 \nabla(\varphi/w_n)
   =\frac{\nabla\varphi}{w_n}
     -\frac{\varphi\nabla w_n}{w_n^2},
\]
whence $\|\varphi/w_n\|_{\HH}\le C_{K,\varepsilon}$ and
$\|\varphi/w_n\|_\infty\le\|\varphi\|_\infty/c_K$. As a result, after
normalization these functions belong to the test class defining $d$. Their integrals against \(\sigma_n-\sigma\)
therefore tend to zero. Local \(L^2\) convergence gives weak convergence in
\(\HH\) to \(\varphi/w\). To pass to the limit in the integrals against
\(\sigma+\varepsilon\beta_0\), we first restrict this measure to a compact set of finite
energy integral. The supremum bound controls the integral on the discarded
part. Thus
\[
 -\Delta_\Cyl w=w^{-1}(\sigma+\varepsilon\beta_0).
\]
The local lower bound
also holds q.e. for \(w\). As in Proposition~\ref{form}, the positive
reaction defines an energy functional. Uniqueness and convergence of the
squared norms give strong convergence. Finally, by \eqref{change-bound} and the exact $\gamma=1$ energy,
\[
 \|v_{\sigma_n+\varepsilon\beta_0}-v_{\sigma_n}\|^2
 \le E(\varepsilon\beta_0)
 =\varepsilon\beta_0(\Cyl),
\]
with the identical bound for the limiting source. It is independent of
$n$. First let $n\to\infty$ with $\varepsilon$ fixed, and then let
$\varepsilon\downarrow0$ to obtain strong convergence for the original
sources.

For measures with bounded mass, convergence in $d$ implies vague
convergence. Vague convergence and convergence of the masses imply the
stated two-ended tightness: choose continuity points $a,b$ of the
limiting radial marginal so that its mass outside $[a,b]$ is small;
weak convergence on the compact subcylinder and convergence of total
masses control the tails for all sufficiently large indices, while the
finitely many remaining measures can be controlled by enlarging
$[a,b]$. Conversely, for a given test $h$ defining $d$, multiply
it by a smooth cutoff supported in a large fixed compact
subcylinder and equal to one on a slightly smaller one. The energy
norms of the products are uniformly controlled using the trace and
product inequalities; the integrals over the deleted tails are at most
the total tail masses because $|h|\le1$. The local
condenser--distance estimate applies on this fixed compact set, so
first the approximation index and then the cutoff can be removed.
This yields convergence in $d$ and proves the local formulation.
\end{proof}

\subsection{Minimum actions and nested condensers}
For \(\gamma\ne1\), denote by \(m_\gamma(\sigma)\) the minimum of the
functional in Proposition~\ref{form}. Its value is
\((1/2-1/(1-\gamma))E(\sigma)\), and we take \(m_\gamma(0)=0\). For two
admissible sources, set \(\bar\sigma=(\sigma+\tau)/2\) and define
\[
 G_\gamma(\sigma,\tau)
 =m_\gamma(\bar\sigma)-\tfrac12(m_\gamma(\sigma)+m_\gamma(\tau)).
\]
For \(\gamma=1\), we first use the minima of the logarithmic functionals on
compact restrictions of finite energy integral. We then define \(G_1\) by
their common limit along an exhaustion as described below. We use this
difference of minima because the global logarithmic functional may be
infinite even for a finite diffuse source. For example, let
\[
 \zeta=\sum_{j\ge1}j^{-2}\delta_{\exp(-j^2)}.
\]
By the trace estimate and Jensen's inequality,
\[
 \int(-\log h)_+\,d(\zeta\otimes\omega)=\infty
\]
for every nonnegative energy function.

\begin{theorem}[Minimum-action criterion]\label{gap}
Let \(\sigma,\tau\) be admissible sources. The limit in the definition of
\(G_1\) exists and is independent of the exhaustion. For every \(\gamma>0\),
\begin{equation}\label{gap-bounds}
 \tfrac18\|v_\sigma-v_\tau\|^2
 \le G_\gamma(\sigma,\tau)
 \le\tfrac14(\|v_\sigma\|+\|v_\tau\|)\|v_\sigma-v_\tau\|.
\end{equation}
In particular, strong convergence in energy holds if and only if
\(G_\gamma(\sigma_n,\sigma)\to0\).
\end{theorem}

\begin{proof}
Set \(u=v_\sigma\), \(v=v_\tau\) and \(z=v_{\bar\sigma}\). For the power
functionals and for the logarithmic functionals on compact restrictions,
the functional $J_\sigma$ is strongly convex with quadratic modulus
one: its quadratic term is $\|\cdot\|^2/2$, while
$-q^{-1}\int h^q\,d\sigma$ (or $-\int\log h\,d\sigma$) is convex on the
positive cone. Hence for $0<\theta<1$,
\[
 J_\sigma((1-\theta)u+\theta z)
 \le (1-\theta)J_\sigma(u)+\theta J_\sigma(z)
     -\tfrac12\theta(1-\theta)\|u-z\|^2.
\]
The minimality of $u$ and division by $\theta$, followed by
$\theta\downarrow0$, yield
\[
 J_\sigma(z)-J_\sigma(u)\ge\tfrac12\|z-u\|^2.
\] We use the same inequality
for \(\tau\) and take the average, using the affine dependence on the
measure. This gives
\[
 G_\gamma\ge(\|z-u\|^2+\|z-v\|^2)/4
\]
and proves the
lower bound. For the upper bound, use \(w=\max(u,v)\) as a competitor for
the average measure. The primitive of the source term is decreasing, so replacing $u$ or
$v$ by their maximum cannot increase the corresponding source
contribution. For $w=u\vee v$, the local Dirichlet form has the lattice
identity
\[
 \|u\vee v\|^2+\|u\wedge v\|^2
       =\|u\|^2+\|v\|^2,
\]
which follows from
$\nabla(u\vee v)=\mathbf1_{\{u\ge v\}}\nabla u+
\mathbf1_{\{v>u\}}\nabla v$ almost everywhere. Using it gives
\[
 G_\gamma\le\tfrac12\|w\|^2-\tfrac14(\|u\|^2+\|v\|^2)
     =\tfrac14\B(u+v,|u-v|).
\]
The Cauchy--Schwarz inequality gives the upper bound.

For finite diffuse sources and \(\gamma=1\), exhaust both measures by
increasing compact sets \(F_j\) such that \((\sigma+\tau)|_{F_j}\) has
finite energy integral. Write \(u_j,v_j,z_j\) for the three solutions for
the restricted sources. For compact restrictions, the logarithmic functional is
$J_\sigma(h)=\|h\|^2/2-\int\log h\,d\sigma$. At its minimizer $u_j$,
$\|u_j\|^2=\sigma(F_j)$; similarly for $v_j$ and $z_j$. Since the mass
of the average source is the average of the masses, the quadratic terms
cancel in the three-minimum difference. Consequently
\[
 2G_{1,j}
   =\int_{F_j}\log(u_j/z_j)\,d\sigma
     +\int_{F_j}\log(v_j/z_j)\,d\tau.
\] For a nonzero component, once the
restriction of \(\sigma\) is nonzero, homogeneity and comparison give
\(z_j/u_j\ge
1/\sqrt2\). We also have
\[
 \int_{F_j}z_j/u_j\,d\sigma
   =\B(u_j,z_j)\le\|u_j\|\|z_j\|.
\]
The right-hand side is bounded independently of the restriction. Hence the
integral of \(\log(z_j/u_j)\) on \(\{z_j/u_j>A\}\) for \(A>e\) is bounded by
\(C(\log A)/A\). The same estimate holds for the other component. By
comparison and Proposition~\ref{form}, the profiles increase q.e. and
converge strongly in energy. We truncate the logarithms and use bounded
convergence. To pass from this formula to its limiting counterpart, note that
$z_j\ge 2^{-1/2}u_j$ by comparison and homogeneity. Thus the negative
part of $\log(z_j/u_j)$ is bounded by $\log\sqrt2$. Also
\[
 \int_{\{z_j/u_j>A\}}\log(z_j/u_j)\,d\sigma
 \le\frac{\log A}{A}\int z_j/u_j\,d\sigma
 \le C\frac{\log A}{A}\quad(A>e).
\]
These bounds supply uniform integrability with respect to the varying
restricted measures: the negative part of
$\log(z_j/u_j)$ is bounded by $\log\sqrt2$, while its
positive tail above $\log A$ has integral at most
$C(\log A)/A$. For fixed $A$ the truncated function
$T_A(\log(u_j/z_j))$ is uniformly bounded. On the increasing sets
$F_j$, quasi-everywhere convergence of $u_j,z_j$ and diffuseness of
$\sigma$ give pointwise convergence outside a set of
$\sigma$-measure zero; bounded convergence first applies on each
$F_k$, and the discarded part has mass at most
$\sigma(\Cyl\setminus F_k)$. Letting $j\to\infty$,
then $k\to\infty$, and finally $A\to\infty$ is justified
by the uniform tail estimate. Repeat the argument for $v_j$ and
$\tau$. We obtain the finite expression
\[
 2G_1(\sigma,\tau)=\int\log(u/z)\,d\sigma
                      +\int\log(v/z)\,d\tau.
\]
We take terms for zero measures to be zero. The right-hand side depends
only on the full solutions and measures, not on the compact sets
used to exhaust them; it therefore proves independence of the
exhaustion as well as passage of the energy bounds to the limit.
\end{proof}

We can express the minima in Theorem~\ref{gap} in terms of the masses of
layers. For \(\gamma\ne1\), put \(P(a)=-a^{1-\gamma}/(1-\gamma)\). For a
nonzero compactly supported admissible measure \(\xi\), define
\(C_\gamma(\xi)\) as the infimum of
\begin{equation}\label{level-cost}
 \tfrac12\sum_{i=1}^m(a_i-a_{i-1})^2\ecap(K_i,U_i)
       +\sum_{i=1}^m P(a_i)\,\xi(K_i\setminus K_{i+1})
\end{equation}
over levels \(0=a_0<a_1<\cdots<a_m\) and condensers
\(K_i\subset U_i\Subset\Cyl\) satisfying \(\supp\xi\subset K_1\),
\(\overline U_{i+1}\subset K_i\), and \(K_{m+1}=\varnothing\). We set
\(C_\gamma(0)=0\). For general admissible \(\xi\), we take the limit along a
compact radial exhaustion.

\begin{theorem}[Nested-condenser formula]\label{levels}
For \(\gamma\ne1\), the limit along the exhaustion exists and
\(C_\gamma(\xi)=m_\gamma(\xi)\) for every admissible \(\xi\). Thus we can
replace the three minima in \(G_\gamma\) by \(C_\gamma\) and test strong
convergence using only condenser capacities and masses of layers.
\end{theorem}

\begin{proof}
First assume that the support is compact. Choose nearly minimizing condenser
tests \(e_i\), truncated to \([0,1]\) and equal to one near \(K_i\). The nesting $\overline U_{i+1}\subset K_i$ implies that $e_i$ is
constant on the support of $\nabla e_j$ whenever $j>i$: the gradients
of the separate layers therefore have pairwise disjoint supports. In
particular,
\[
 \Bigl\|\sum_i(a_i-a_{i-1})e_i\Bigr\|^2
     =\sum_i(a_i-a_{i-1})^2\|e_i\|^2.
\]
On $K_i\setminus K_{i+1}$ the resulting sum is at least $a_i$; since
$P$ is decreasing, $P(h)\le P(a_i)$ there. The functional at \(\sum_i(a_i-a_{i-1})e_i\) is bounded above by
\eqref{level-cost} up to the errors in the capacities. For the reverse
bound, approximate the minimizing solution by smooth nonnegative functions
\(h\) of compact support with a positive minimum on \(\supp\xi\). For $0<\gamma<1$, put $q=1-\gamma$. Admissibility and
Proposition~\ref{form} give the trace inequality
\[
 \int |h_n-h|^q\,d\xi
    \le S_q(\xi)\|h_n-h\|^q.
\]
Together with $|a^q-b^q|\le|a-b|^q$ for $a,b\ge0$, this shows
continuity of $\int h^q\,d\xi$ under energy approximation, even for
diffuse compact sources that do not possess bounded densities. For $\gamma>1$, first add a
small fixed positive smooth cutoff $\varepsilon\psi$ on a neighborhood
of $\supp\xi$ to the energy approximants. The integrands are then
bounded by $(\varepsilon\min\psi)^{1-\gamma}$, so quasi-everywhere
convergence and diffuseness give dominated convergence for fixed
$\varepsilon$. As $\varepsilon\downarrow0$, monotone convergence
recovers the negative-power functional of the minimizer. Thus in both
cases a smooth compactly supported nonnegative $h$, strictly positive
on $\supp\xi$, approximates the minimum action. For such an \(h\) take the levels
\(a_i=i\delta\) with \(K_i=\{h\ge a_i\}\) and \(U_i=\{h>a_{i-1}\}\),
choosing \(\delta\) below the positive minimum. A smooth cutoff between the
levels gives
\[
 \delta^2\ecap(K_i,U_i)
     \le\int_{a_{i-1}<h<a_i}|\nabla h|^2\,dm .
\]
For a smooth $h$ we can select regular levels arbitrarily close
to the grid $i\delta$, so that the condenser sets are compactly nested.
The energy bounds on the strips $\{a_{i-1}<h<a_i\}$ sum to
$\int|\nabla h|^2dm$. On $\supp\xi$, the layer step function differs
from $h$ by at most $\delta$, and $h$ has a positive minimum there.
Consequently the potential costs converge to $\int P(h)d\xi$ by uniform
continuity of $P$ on this range. Sending $\delta\downarrow0$ proves
$C_\gamma(\xi)\le J_\xi(h)$, and then minimizing over $h$ gives the
reverse inequality to the preceding test-function construction. For a general admissible $\xi$, let $F_j$ be any increasing
compact exhaustion with $\xi(F_j)\uparrow\xi(\Cyl)$ and set
$\xi_j=\mathbf1_{F_j}\xi$. Proposition~\ref{form} gives
$\|v_{\xi_j}-v_\xi\|\to0$. Since the compact formula
already proved reads
\[
 C_\gamma(\xi_j)=m_\gamma(\xi_j)
   =\left(\frac12-\frac1{1-\gamma}\right)
           \|v_{\xi_j}\|^2,
\]
the sequence of condenser costs converges to $m_\gamma(\xi)$.
The same limit is obtained for every such compact exhaustion, because
the limit is the energy of the unique full solution. Thus the definition
for general sources is independent of the exhaustion.
\end{proof}

\subsection{Concentration on sets of small capacity}
Let \(K\Subset\Cyl\) be compact and have capacity \(C>0\). The equilibrium
potential \(e_K\) satisfies \(0\le e_K\le1\) and \(e_K=1\) q.e.\ on \(K\).
The equilibrium measure \(\nu_K=-\Delta_\Cyl e_K\) is supported on \(K\),
with \(\nu_K(K)=\|e_K\|^2=C\) by \cite[Sec.~2.1]{FOT}.

\begin{theorem}[Small supports and energy defects]\label{balls}
For an admissible measure of mass \(M>0\) supported on \(K\), define
\[
 Q_\gamma(M,C)=M^{2/(1+\gamma)}C^{(\gamma-1)/(1+\gamma)}.
\]
Then
\[
 \begin{cases}
 E(\sigma)\ge Q_\gamma(M,C),&\gamma<1,\\
 E(\sigma)\le Q_\gamma(M,C),&\gamma>1.
 \end{cases}
\]
In both cases equality holds
if and only if the source is \((M/C)\nu_K\). For \(\gamma=1\), the energy
equals \(M\). Every finite diffuse measure supported on \(K\) is admissible
for \(\gamma\ge1\).

Let \(\sigma_j\) be admissible measures with mass \(M\) supported in
interior balls of radii \(\eps_j\downarrow0\) centered at \(x_*\in\Cyl\).
Fix an admissible compactly supported source \(\tau\), possibly zero, and
set \(u_j=v_{\tau+\sigma_j}\) and \(w=v_\tau\). Then \(u_j\to w\) in
\(L^1_{\mathrm{loc}}\) and in the local essential supremum norm away from
\(x_*\). The reactions converge in total variation to the reaction of \(w\).
If \(\gamma<1\), then \(E(\tau+\sigma_j)\to\infty\). If \(\gamma>1\), then
\(u_j\to w\) strongly in \(\HH\). For \(\gamma=1\), we have
\(u_j\rightharpoonup w\) weakly in \(\HH\) and
\begin{equation}\label{defect}
 \|u_j-w\|^2\to M,\qquad
 |\nabla u_j|^2\,dm\ \longrightarrow\ |\nabla w|^2\,dm+M\delta_{x_*}
 \quad\text{narrowly}.
\end{equation}
\end{theorem}

\begin{proof}
Use \(e_K/\sqrt C\) in the supremum or infimum in Proposition~\ref{form}.
Its integral against the source in that expression is
\(MC^{-(1-\gamma)/2}\). If equality holds, minimization along the ray and
strict convexity imply that the minimizer is a multiple of \(e_K\). The
equation then gives \(\sigma=(M/C)\nu_K\), and this measure attains
equality. When $\gamma>1$, the competitor $e_K/\sqrt C$ is strictly positive
quasi-everywhere on $K$ and its negative-power integral equals
$MC^{(\gamma-1)/2}<\infty$; Proposition~\ref{form} gives admissibility
for every finite diffuse source supported on $K$. At $\gamma=1$ the
finite-mass criterion in the same proposition gives admissibility and
$E(\sigma)=M$. The capacity of an interior ball is
comparable to \(\eps^{N-2}\) by a scaled cutoff, the local Sobolev
inequality and \eqref{time-trace}. This gives the two limits for the
energies of \(v_j=v_{\sigma_j}\) when \(\gamma\ne1\).

On interior compact sets, the positive Green kernel of the cylinder with
zero boundary values satisfies \(G(x,y)\asymp d(x,y)^{2-N}\) near the
diagonal. It is smooth off the diagonal and bounded on subsets a positive
distance from it. Indeed, the Green kernel is the time integral of the
product of the reflected half-line heat kernel and the spherical heat
kernel. The local estimates follow from the short-time Gaussian
bounds~\cite{Grig}. For large heat times the half-line factor is
\(O(t^{-3/2})\) on compact sets. For compact $K,K'\Subset\Cyl$, the local singularity is integrable
uniformly in $y\in K$, since in dimension $N$,
\[
 \int_{B_r(y)}d(x,y)^{2-N}\,dm(x)\le Cr^2.
\]
Outside the diagonal the Green function is bounded on compact sets.
The heat-kernel representation and the $t^{-3/2}$ long-time bound for the
reflected half-line factor provide the remaining uniform estimate.
Thus $\sup_{y\in K}\int_{K'}G(x,y)dm(x)<\infty$. Set
\(\rho_j=v_j^{-\gamma}\sigma_j\) and \(H_j=\rho_j(\Cyl)\). This mass is
finite by the local positive lower bound for \(v_j\). The Green
representation on the supporting ball gives
\[
 v_j\ge c\eps_j^{2-N}H_j\quad\text{q.e.}
\]
Therefore
\[
 M=\int v_j^\gamma\,d\rho_j
 \ge(c\eps_j^{2-N}H_j)^\gamma H_j,\qquad
 H_j\le C'\eps_j^{(N-2)\gamma/(1+\gamma)} .
\]
The Green estimates imply that \(v_j\) converges to zero locally and
converges smoothly away from \(x_*\). By comparison and Lemma~\ref{change},
\(w,v_j\le u_j\le w+v_j\). The finite background reaction vanishes on the
polar point \(x_*\). By dominated convergence for this reaction and by
\(u_j^{-\gamma}\sigma_j\le\rho_j\), we obtain convergence of the reactions
in total variation. For \(\gamma<1\), comparison gives the divergence of the
energies. For \(\gamma>1\), strong convergence in energy follows from
\eqref{change-bound}.

When \(\gamma=1\), the mass identity gives \(\|u_j\|^2=\|w\|^2+M\). Local
convergence thus implies weak convergence in energy and the first limit in
\eqref{defect}. For the second limit, let \(\chi\) be a smooth cutoff equal
to zero near \(x_*\) and to one outside a larger interior ball. Put $z_j=u_j-w\ge0$, so $z_j\le v_j$. The product
$\chi^2z_j$ is an admissible energy test by the bounded multiplier
property; approximation by smooth tests justifies its use in the weak
equation. For large $j$ the source $\sigma_j$ is supported where
$\chi=0$, and therefore contributes nothing. On the common background
measure $\tau$, the difference of reactions is
$\bigl(u_j^{-\gamma}-w^{-\gamma}\bigr)\tau$, whose pairing with
$\chi^2z_j\ge0$ is nonpositive. Expanding $\nabla(\chi^2z_j)$, using
Cauchy--Schwarz and absorbing half the principal term, gives
\[
 \int\chi^2|\nabla(u_j-w)|^2\,dm
   \le4\int v_j^2|\nabla\chi|^2\,dm\longrightarrow0 .
\]
Choose a decreasing sequence of interior balls centered at $x_*$,
and take $\chi=1$ on the complement of a slightly larger ball.
The displayed estimate proves
$\nabla(u_j-w)\to0$ strongly on the complement of each such
ball, including both noncompact ends of the cylinder. The sequence
$z_j=u_j-w$ is weakly convergent to zero in $\HH$, so for any bounded
continuous test $\psi$,
\[
 \int\psi\,\nabla w\cdot\nabla z_j\,dm
       \longrightarrow0
\]
by weak convergence in $L^2$ and $\psi\nabla w\in L^2$.
Finally,
$\|z_j\|^2=\|u_j\|^2+\|w\|^2-2\B(u_j,w)\to M$.
The exterior strong convergence forces all this limiting
$\nabla z_j$-energy to lie in every ball about $x_*$. Thus it
converges to $M\delta_{x_*}$ as a finite measure: its total mass is
$M$, and no energy can escape to the ends or concentrate elsewhere.
Expanding $|\nabla u_j|^2=|\nabla w|^2+
2\nabla w\cdot\nabla z_j+|\nabla z_j|^2$ now proves
\eqref{defect}.
\end{proof}

\begin{remark}
The conditions in \eqref{two-plates} depend on the inverse power. For an
inverse power below one, take equilibrium sources of mass
\(C_j^{(1-\gamma)/2}\) on balls of capacity \(C_j\downarrow0\). Their masses
and their condenser discrepancies from zero tend to zero, and their energies
equal one. For an inverse power above one, take radial sources of mass one
at height \(s\to\infty\). Their energies are
\(s^{(1-\gamma)/(1+\gamma)}\to0\).
\end{remark}

\begin{remark}[Compactly supported sources with uniformly small Green singularities]
Let $\mu_n$ be positive measures on a fixed compact subset $K\Subset\Cyl$ converging narrowly to $\mu\ne0$. Suppose
\[
 \lim_{r\downarrow0}\sup_n\sup_{x\in K}
 \int_{d(x,y)<r}d(x,y)^{2-N}\,d\mu_n(y)=0,
\]
where $d$ is the cylinder distance. Then, for each $\gamma>0$, the corresponding solutions converge strongly in $\HH$. They also converge uniformly on $K$ in their continuous representatives.

To prove this, proceed as follows.
The same small-radius bound holds for $\mu$ by the lower semicontinuity inequality for narrow
convergence, since the integrand with the open-ball cutoff is nonnegative and lower semicontinuous.
The Green function $G(x,y)$ of the half-cylinder is continuous off the diagonal, positive, and on
$K$ has the local upper bound $C d(x,y)^{2-N}$ near the diagonal. It also has a positive lower bound
on $K\times K$ (using the local Green singularity near the diagonal and compactness away from it).
Since the masses converge to a positive number, for all sufficiently large $n$,
\[
 0<m\leq G\mu_n(x),\,G\mu(x)\leq M<\infty\qquad(x\in K).
\]
The same local bound and a continuous cutoff of the kernel near the diagonal show that these
potentials are continuous on $K$. They have finite energy: for a positive measure $\nu$ among these
measures, $\int G\nu\,d\nu\leq M\nu(K)<\infty$, which is the Green-energy criterion for finite
energy integral.

For $\nu=\mu_n$ or $\mu$, the functions $cG\nu$ and $CG\nu$ are sub- and supersolutions,
respectively, when the positive constants are chosen so that $c^{1+\gamma}M^\gamma\leq1\leq
C^{1+\gamma}m^\gamma$. Admissibility follows from Proposition~\ref{form}: for positive powers use
the finite-energy-integral bound and H\"older's inequality, at inverse power one use finite mass,
and for negative powers use $G\nu$ as a competitor. Comparison gives uniform positive upper and
lower bounds for $v_\nu$ on $K$ in the quasi-continuous sense. The reaction $v_\nu^{-\gamma}\nu$ is
therefore bounded by a common multiple of $\nu$ and has finite energy integral. The Green
representation $v_\nu=G(v_\nu^{-\gamma}\nu)$ and the same diagonal cutoff argument show that each
$v_\nu$ has a continuous representative on $K$. This step uses only the boundedness of the density
with respect to $\nu$, not its continuity.

In particular, for $v=v_\mu$, the function $v^{-\gamma}$ is now known to be continuous and bounded
on $K$. Set $z_n=G(v^{-\gamma}\mu_n)$. For fixed $r>0$ multiply $G$ by a continuous cutoff which
vanishes for $d(x,y)\leq r/2$ and equals one for $d(x,y)\geq r$. The resulting kernels, multiplied
by $v(y)^{-\gamma}$, form a compact family in $C(K)$ as $x$ ranges in $K$. Narrow convergence is
uniform on this family by a finite-net argument. The discarded portions tend to zero uniformly in
$n$ and $x$ by the assumed small-radius estimate. Consequently
\[
 \|z_n-v\|_{C(K)}\longrightarrow0.
\]
If $\eps_n=\|z_n/v-1\|_{C(K)}<1$, comparison in the nonlinear equation yields
\[
 (1+\eps_n)^{-\gamma/(1+\gamma)}z_n
 \leq v_{\mu_n}\leq
 (1-\eps_n)^{-\gamma/(1+\gamma)}z_n.
\]
Indeed, the displayed multiples satisfy the sub- and supersolution inequalities on the support of $\mu_n$; outside that support both are harmonic. Thus $v_{\mu_n}\to v$ uniformly on $K$.

Finally, the Green equation tested with $z_n-v$ gives
\[
 \|z_n-v\|^2=\int_K (z_n-v)v^{-\gamma}\,d(\mu_n-\mu)
 \leq C_0\|z_n-v\|_{C(K)}\longrightarrow0.
\]
Testing the difference of the equations for $v_{\mu_n}$ and $z_n$ with their difference gives
\begin{align*}
 \|v_{\mu_n}-z_n\|^2
 &=\int_K (v_{\mu_n}-z_n)(v_{\mu_n}^{-\gamma}-v^{-\gamma})\,d\mu_n\\
 &\leq C_1\|v_{\mu_n}-z_n\|_{C(K)}\longrightarrow0,
\end{align*}
where the coefficient is uniformly bounded by the positive lower bounds and the bounded masses. Both
tests are legitimate by finite energy integral. The triangle inequality proves strong convergence.
The hypothesis follows, for example, from a uniform growth estimate $\mu_n(B_r(x))\leq C r^a$ with
$a>N-2$, by summing the dyadic annuli around $x$.
\end{remark}

\begin{remark}[Radial transport]\label{transport}
For radial measures of equal finite mass, consider the linear interpolation
of their nondecreasing quantiles. Along this interpolation the energy is
concave if \(\gamma<1\), constant if \(\gamma=1\), and convex if
\(\gamma>1\). Indeed, for equal masses at ordered positions \(x_i\), the
minimum of the functional with prescribed values \(y_i\) contains the
quadratic term \(\tfrac12\sum (y_i-y_{i-1})^2/(x_i-x_{i-1})\), with
\(x_0=y_0=0\). The quotient is jointly convex when extended by lower
semicontinuity to repeated positions. The primitive is convex in \(y_i\), so
the minimum over the values is convex in the positions. The factor relating
the minimum to the energy changes sign at inverse power one. For compact
supports, sampling the quantiles gives the result by Theorem~\ref{angular}.
For general supports, restrict the quantiles to nested compact subintervals
of the mass parameter. When \(\gamma<1\), admissibility along the
interpolation follows by bounding \(|h(t)|\) by \(\int_0^t|h'|\) in the
trace inequality. When \(\gamma>1\), the initial profile is a competitor
with finite negative moment at a parameter \(\theta<1\) because
\[
 V_0((1-\theta)x+\theta y)\ge(1-\theta)V_0(x).
\]
We then pass to the limit
using Proposition~\ref{form}.

For equal-mass data on a fixed interior radial interval with uniform
positive lower and finite upper mass bounds, we have
\[
 \|V_1-V_2\|_\infty\le C W_1,
 \qquad \|(V_1-V_2)'\|_2\le C W_1^{1/2}.
\]
Indeed, the Green formula gives a uniform upper bound and a positive lower
bound for both profiles on the interval, and hence a uniform Lipschitz
bound. By Lipschitz duality the potential
\[
 z=\int\min(t,s)V_1(s)^{-\gamma}d(\zeta_2-\zeta_1)(s)
\]
is bounded in the
supremum norm by \(CW_1\). Subtract the equations for the profiles and use
the nonnegative zeroth-order coefficient given by the common inverse-power
law to compare \(V_2-V_1-z\) with \(\pm\|z\|_\infty\). Testing the
difference equation with the difference bounds the squared energy norm by
\(C\|V_2-V_1\|_\infty\). If a single atom is moved by a distance \(h\), the
pole trace changes by order \(h\) and the slope changes by order one on the
interval between the two positions. This gives the two rates.
\end{remark}
\section{The power balance at a Hardy pole}\label{sec-flux}

We now consider a linear Hardy equation. We test it with a power of the
solution and determine the term at the pole in the resulting identity.

\subsection{The fractional balance}
Let $N\ge2$, $0<s<1$ and $p>2$. Throughout this section,
\[
 \beta=\frac{N-2s}{p}>0,
 \qquad d=N-2s-2\beta=(N-2s)\left(1-\frac2p\right)>0.
\] We use the Fourier normalization of the fractional Laplacian, with
kernel constant
\[
 C_{N,s}=\frac{4^s\Gamma(N/2+s)}{\pi^{N/2}|\Gamma(-s)|}.
\]
The action of the fractional Laplacian on powers is
\[
 (-\Delta)^s|x|^{-b}=\kappa(b)|x|^{-b-2s},
\]
where
\[
 \kappa(b)=2^{2s}\frac{\Gamma((N-b)/2)\Gamma((b+2s)/2)}
             {\Gamma(b/2)\Gamma((N-b-2s)/2)}.
\]
Set \(h_\beta(x)=|x|^{-\beta}\) and define
\[
 L_\lambda=(-\Delta)^s-\lambda|x|^{-2s},
\]
where \(\lambda=\kappa(\beta)\) is the sharp \(L^p\) Hardy
coefficient~\cite{BogdanLp}. This coefficient is strictly below the
quadratic Hardy constant. Hence the quadratic form is coercive on
\(X_0^{s,2}(B_1)\), the closure in \(H^s(\R^N)\) of smooth test functions in
\(B_1\) extended by zero.

Write
\[
 u(x)=h_\beta(x)V(\log(1/|x|),x/|x|),
 \qquad V=0\quad(t<0),
\]
and set
\[
 D_p(y,x)=y^p-x^p-px^{p-1}(y-x).
\]
In these variables, the Sobolev--Bregman
remainder in \cite{BogdanLp} is
\begin{equation}\label{frac-rem}
 \mathcal R_{s,p}(V)=\frac{\sph}{p}\int_{\R^2}
   \iint_{(\mathbb S^{N-1})^2}
    D_p(V(t+h,\varphi),V(t,\theta))P(h,\theta,\varphi)
      \,d\omega(\varphi)\,d\omega(\theta)\,dh\,dt ,
\end{equation}
where
\[
 P(h,\theta,\varphi)
 =C_{N,s}\sph\frac{e^{-(N-\beta)h}}{|\theta-e^{-h}\varphi|^{N+2s}} .
\]
The kernel is obtained by setting
\[
 x=e^{-t}\theta,\qquad y=e^{-(t+h)}\varphi
\]
in the weighted ground-state form. The polar Jacobians contribute
$e^{-Nt}e^{-N(t+h)}$, the fractional kernel contributes
$e^{(N+2s)t}|\theta-e^{-h}\varphi|^{-N-2s}$, and the factor
$u(x)^{p-1}h_\beta(y)$ contributes
$e^{p\beta t+\beta h}$. Since $p\beta=N-2s$, the $t$-dependence
cancels. The remaining radial factor is $e^{-(N-\beta)h}$,
and $dS=\sph d\omega$ supplies the angular normalization in
\eqref{frac-rem}. This is the change of variables in the fractional
$L^p$ ground-state representation of \cite{BogdanLp}.

Both angular marginals of $P$ are given by
\[
 k(h)=C_{N,s}e^{-(N-\beta)h}
          \int_{\mathbb S^{N-1}}|e_1-e^{-h}\theta|^{-N-2s}dS(\theta)
\]
by rotation invariance. For radial profiles, we can therefore replace the
angular integral in \eqref{frac-rem} by integration against the kernel
\(k(h)\).

\begin{theorem}[Pole flux and angular traces]\label{angular-flux}
Let $u=h_\beta V\in X_0^{s,2}(B_1)$, where $V\ge0$ is bounded.
Assume that $L_\lambda u$ agrees, as a distribution on $B_1$, with an
$L^1(B_1,dx)$ \emph{function} $f$, rather than merely with a finite
measure; that is,
\[
 \langle L_\lambda u,\varphi\rangle=\int_{B_1}f\varphi\,dx
 \qquad(\varphi\in C_c^\infty(B_1)).
\]
Suppose that
\[
 g(t)=\int V(t,\theta)^p\,d\omega(\theta)
\]
has essential limit \(\ell\) at
infinity. Since $V$ is bounded and nonnegative, necessarily
$0\le\ell<\infty$; ``extended values'' refers to the
possible value $+\infty$ of the source pairing and the nonnegative
remainder, not to $\ell$. If \(f\ge0\), then the following identity holds with extended
values.
\begin{equation}\label{angular-balance}
 \int_{B_1}u^{p-1}f\,dx
    =\mathcal R_{s,p}(V)+\frac{\sph}{p}\kappa'(\beta)\ell .
\end{equation}
The same identity holds for signed \(f\) provided that
\(u^{p-1}f\in L^1(B_1)\). If \(\ell>0\), the pole term is positive and the
Hardy integral
\[
 \int_{B_1}|x|^{-2s}u^p\,dx=\sph\int_0^\infty g(t)\,dt
\]
is
infinite.

Let the nonnegative profile be locally \(L^p\) and have finite remainder,
and let
\[
 \bar V(t)=\int V(t,\theta)d\omega(\theta)
\]
be its spherical mean.
Then
\begin{equation}\label{angular-cost}
 \int_0^\infty\!\int_{\mathbb S^{N-1}}
       |V(t,\theta)-\bar V(t)|^p\,d\omega(\theta)\,dt
       \le C_{N,s,p}\mathcal R_{s,p}(V).
\end{equation}
In particular, if \(V(t,\cdot)\) has an essential limit in \(L^p(d\omega)\)
as \(t\to\infty\), the limit is constant.
\end{theorem}

\begin{lemma}[First moment of the radial kernel]\label{kernel-moment}
Set \(d=N-2s-2\beta>0\). For \(h>0\) we have \(k(-h)=e^{dh}k(h)\). Near
zero, \(k(h)\sim c_{N,s}|h|^{-1-2s}\) for some \(c_{N,s}>0\). At both ends
the kernel decays exponentially. In particular,
\begin{equation}\label{first-moment}
 m_\beta=\int_0^\infty h[k(-h)-k(h)]\,dh=\kappa'(\beta)>0 .
\end{equation}
\end{lemma}

\begin{proof}
We use
\[
 |e_1-e^{-h}\theta|^2=(1-e^{-h})^2+e^{-h}|e_1-\theta|^2.
\]
For $|h|\ll1$ only a spherical cap about $e_1$ produces a singular
contribution. In normal coordinates $\theta=(\sqrt{1-|z|^2},z)$,
$z\in\mathbb R^{N-1}$, the denominator is comparable to
$\bigl(h^2+|z|^2\bigr)^{(N+2s)/2}$. Rescale $z=|h|w$; the cap integral
then satisfies
\[
 |h|^{1+2s}k(h)\longrightarrow
 C_{N,s}\int_{\mathbb R^{N-1}}
       \frac{dw}{(1+|w|^2)^{(N+2s)/2}}=:c_{N,s}>0.
\]
Outside that cap the contribution is uniformly bounded; the rescaled
integrands are dominated by an integrable multiple of
$(1+|w|^2)^{-(N+2s)/2}$.

For the reflection identity, factor $e^h$ from the distance in $k(-h)$:
\[
 \begin{aligned}
 k(-h)&=C_{N,s}e^{(N-\beta)h}
      \int_{\mathbb S^{N-1}}|e_1-e^h\theta|^{-N-2s}dS(\theta)\\
 &=C_{N,s}e^{-(\beta+2s)h}
      \int_{\mathbb S^{N-1}}|e_1-e^{-h}\theta|^{-N-2s}dS(\theta)\\
 &=e^{(N-2s-2\beta)h}k(h).
 \end{aligned}
\]
For $h\to+\infty$, the integral in $k(h)$ stays bounded, so
$k(h)=O(e^{-(N-\beta)h})$; the reflection formula gives
$k(-h)=O(e^{-(\beta+2s)h})$. In particular, the first moments of both
tails are finite.

Fix a small real \(\xi\) and a test function supported on a compact annulus.
Subtract the ground-state equation for \(h_\beta\) from that for
\(|x|^{-\beta-\xi}\), using the test function multiplied by \(|x|^{-\xi}\)
in the equation for the first ground state. In both bilinear integrals we
first impose the symmetric exclusion \(|\log(|x|/|y|)|>\epsilon\). The two powers are smooth on the test annulus. For each form,
the product of differences is absolutely integrable on a larger annulus
by the local $H^s$ seminorm bounds; outside it, the test support is separated
from the other coordinate and the power is integrable against the kernel.
Thus dominated convergence removes the exclusion in the bilinear forms. Changing variables by \(|y|=|x|e^{-h}\) gives
\[
 \kappa(\beta+\xi)-\kappa(\beta)
    =\int_{\R}(1-e^{\xi h})k(h)\,dh ,
\]
where we pair the two signs of \(h\). Pairing $h$ and $-h$ rewrites the integrand as
\[
 (1-e^{\xi h})k(h)+(1-e^{-\xi h})k(-h).
\]
For $|\xi|\le\xi_0$ sufficiently small, Taylor expansion and
$k(-h)-k(h)=O(h^{-2s})$ show that this expression and its
$\xi$-derivative are $O(h^{1-2s})$ as $h\downarrow0$. This is
integrable for every $s<1$. The exponential bounds just established
provide an integrable majorant at infinity, possibly after reducing
$\xi_0$. Dominated differentiation at $\xi=0$ therefore yields
\[
 \kappa'(\beta)
   =\int_0^\infty h\bigl(k(-h)-k(h)\bigr)\,dh=m_\beta.
\]
Strict positivity follows from $d>0$ and the reflection identity.
\end{proof}

The moment gives the constant in the linear weighted pole identity of Chen
and Weth~\cite[Theorem~1.3]{CW}. For \(\mu=-\lambda\) their two powers are
\(|x|^{-\beta}\) and \(|x|^{-(p-1)\beta}\). Using the change of variables
\(z=e^{-h}\theta\) in \eqref{first-moment}, we obtain
\begin{equation}\label{linear-normalization}
 \sph m_\beta=C_{N,s}\sph\int_0^1\frac{da}{a}
   \int_{B_a}\frac{|z|^{-(p-1)\beta}-|z|^{-\beta}}
                      {|e_1-z|^{N+2s}}\,dz .
\end{equation}
This is the constant in their proof, with the radial measure \(da/a\) from
the change to polar coordinates, see \cite[Eqs.~(3.19)--(3.24)]{CW}. In
\eqref{angular-balance} this constant multiplies the limiting spherical mean
of \(V^p/p\) from the power test.

\begin{proof}[Proof of Theorem~\ref{angular-flux}]
Write $h=h_\beta$ and $M=\|V\|_{L^\infty}$. The quadratic form in the
Fourier normalization is
\[
 \mathcal E_s(a,b)=\frac{C_{N,s}}2\iint_{\R^N\times\R^N}
 \frac{(a(x)-a(y))(b(x)-b(y))}{|x-y|^{N+2s}}\,dx\,dy.
\]
\emph{Annular testing.}
Let $w\in H^s(\R^N)\cap L^\infty(\R^N)$ have zero exterior values on $B_1^c$ and vanish in a neighborhood of zero. Then
\[
 \mathcal E_s(u,w)-\lambda\int |x|^{-2s}uw\,d x=\int_{B_1}fw\,d x.
\]
The same testing procedure gives
\[
 \mathcal E_s(h,w)=\lambda\int |x|^{-2s}hw\,d x,
\]
where the mixed integral with $h$ is absolutely convergent as an integral of the product of differences.

Indeed, take inward dilations $w(x/r_j)$, $r_j\uparrow1$, followed by mollification at radii smaller
than both the gap to $\partial B_1$ and half the distance of their supports to the pole. These
operations converge in $H^s(\R^N)$ and preserve a common $L^\infty$ bound; the functions are smooth
with compact support in a fixed annulus inside the punctured ball, apart from the outer radius
approaching one. Pass to a subsequence converging almost everywhere. The form for $u$ converges by
Cauchy--Schwarz, the Hardy term by $L^2$ convergence on the annulus, and the right-hand side by
bounded convergence against $|f|\,d x$.

For the mixed integral with $h$, let $A$ be a fixed closed annulus containing the supports of all
tests and let $A'$ be a larger annulus with $A\Subset A'$. Choose a compactly supported smooth
cutoff equal to one on $A'$; its product with $h$ is a smooth compactly supported function. For
pairs in $A'\times A'$ the absolute integral of the product of differences is bounded by the product
of the two $H^s$ seminorms. For a pair with one coordinate in $A$ and the other outside $A'$, the
distance between the coordinates is bounded below. The integral of $(|h(x)|+|h(y)|)|x-y|^{-N-2s}$
over $y\notin A'$ is uniformly finite for $x\in A$: near zero this follows from $\beta<N$, and at
infinity from decay of the kernel. Since the test vanishes outside $A$, these separated-pair
contributions are bounded by a constant times its $L^1(A)$ norm. The distributional identity for $h$
therefore extends to $w$ by the same approximation in $H^s$. No classical smoothness of $V$ is used.

Choose $\chi\in C^\infty(\R)$, $0\leq\chi\leq1$, nonincreasing, equal to one on $(-\infty,0]$ and zero on $[1,\infty)$. Put $\chi_T(t)=\chi(t-T)$, and define tests (extended by zero outside the ball)
\[
 \phi_T(x)=\chi_T(\log(1/|x|))u(x)^{p-1},
 \qquad \psi_T(x)=h(x)^{-1}u(x)\phi_T(x).
\]
Both belong to $H^s\cap L^\infty$ and vanish near the pole. To see this explicitly, on the support
of the cutoff $u$ is bounded by $M e^{\beta(T+1)}$. Truncate $u$ at that level before composing with
the power functions $z^{p-1}$ and $z^p$, which are Lipschitz on the bounded nonnegative range and
vanish at zero. Multiplication by the smooth cutoff, and by $h^{-1}$ on its support in the second
test, preserves $H^s$. Thus the annular testing argument applies. The two testing identities, with
tests $\phi_T$ and $\psi_T$, have identical Hardy terms since $h\psi_T=u\phi_T$.

For $\eps>0$ impose in both forms the symmetric restriction $|\log(|x|/|y|)|>\eps$. The absolute
convergence established in the annular testing argument and Cauchy--Schwarz for the form of $u$
imply convergence of these restricted forms as $\eps\downarrow0$. The restriction exhausts almost
every pair, so there is no need to invoke a pointwise principal value or Taylor regularity of $V$.
On the restricted domain, the exact algebraic identity (with subscripts denoting evaluation at $x$
and $y$, not derivatives, and $\phi=\phi_T$)
\begin{align*}
 &(u_x-u_y)(\phi_x-\phi_y)
 -(h_x-h_y)(V_x\phi_x-V_y\phi_y)\\
 &\hspace{8mm}=h_y\phi_x(V_x-V_y)+h_x\phi_y(V_y-V_x)
\end{align*}
and symmetry of the restriction yield the integral of the first term on the right, multiplied by
$C_{N,s}$. For each fixed $\eps$ this integral is absolutely convergent: the $x$ support is in an
annulus, the logarithmic restriction separates $x$ and $y$, $V$ is bounded, and $h$ is integrable
against the separated kernel. In logarithmic coordinates $x=e^{-t}\theta$, $y=e^{-(t+a)}\varphi$,
the radial weights cancel because $p\beta=N-2s$. In the notation of \eqref{frac-rem} we obtain
\begin{align*}
 \int_{B_1}\phi_T f\,d x
 =\sph\lim_{\eps\downarrow0}
 \int_{|a|>\eps}\!\int_{\R}\!\iint
 &\chi_T(t)V(t,\theta)^{p-1}\big(V(t,\theta)-V(t+a,\varphi)\big)\\[-2pt]
 &\cdot P(a,\theta,\varphi)\,d\omega(\varphi)\,d\omega(\theta)\,d t\,d a.
\end{align*}
Apply $pX^{p-1}(X-Y)=D_p(Y,X)+X^p-Y^p$. The nonnegative first term converges monotonically to the
remainder $\mathcal R_{s,p}^{\chi_T}(V)$ with the extra weight $\chi_T(t)$. The angular marginals of
$P$ are both $k(a)$ by rotational invariance. Translating $t$ in the second term and pairing
positive and negative increments yields
\begin{align*}
 B_T&=\int_0^\infty\int_{\R}g(t)\,b_{T,a}(t)\,d t\,d a,\\
 b_{T,a}(t)&=k(a)\big(\chi_T(t)-\chi_T(t-a)\big)
      +k(-a)\big(\chi_T(t)-\chi_T(t+a)\big).
\end{align*}
For $0<a<1$, split the coefficient into its symmetric second difference with factor $k(a)$ and its first difference with factor $k(-a)-k(a)\geq0$. The $L^1(\R)$ norm is bounded by
\[
 a^2 k(a)\|\chi''\|_1+a\big(k(-a)-k(a)\big)\|\chi'\|_1.
\]
For $a\geq1$ the bound $a(k(a)+k(-a))\|\chi'\|_1$ suffices. Both bounds are independent of $T$ and
integrable in $a$: Lemma~\ref{kernel-moment} gives $k(a)=O(a^{-1-2s})$, $k(-a)-k(a)=O(a^{-2s})$ at
zero and exponential decay at infinity. As $0\leq g\leq M^p$, the limit in the signed second term
exists by absolute integration. In particular, for a finite constant $C_*$ depending only on the
kernel and $\chi$,
\begin{equation}\label{eq:fluxbound}
 |B_T|\leq C_* M^p\qquad\text{for all }T>0.
\end{equation}
The localized identity is now justified without additional smoothness hypotheses:
\begin{equation}\label{eq:localbalance}
 \int_{B_1}\chi_T u^{p-1}f\,d x
 =\mathcal R_{s,p}^{\chi_T}(V)+\frac{\sph}{p}B_T.
\end{equation}
Each localized source pairing is finite because the test is bounded and $f\in L^1$; hence the localized remainder is finite as well.

If $g$ has essential limit $\ell$ at infinity, translate $t=T+z$. For each fixed increment $a$, the
function $b_{T,a}(T+z)=b_{0,a}(z)$ has bounded support in $z$. On this support, $g(T+z)$ approaches
$\ell$ in essential supremum as $T\to\infty$. Its integral against $b_{0,a}$ therefore tends to
$\ell\int b_{0,a}$. The preceding integrable bounds justify convergence also in the increment
variable. Since
\[
 \int_{\R}(\chi(z)-\chi(z-a))\,d z=-a,
 \qquad \int_{\R}(\chi(z)-\chi(z+a))\,d z=a,
\]
we conclude from Lemma~\ref{kernel-moment} that
\[
 \lim_{T\to\infty}B_T=\ell\int_0^\infty a\big(k(-a)-k(a)\big)\,d a
 =\kappa'(\beta)\ell.
\]
The cutoff increases to one. Monotone convergence for the remainder and for the source pairing when $f\geq0$ (or dominated convergence when $u^{p-1}f\in L^1$) yields
\[
 \int_{B_1}u^{p-1}f\,d x=\mathcal R_{s,p}(V)+\frac{\sph}{p}\kappa'(\beta)\ell.
\]
The coefficient is positive since $k(-a)>k(a)$ for $a>0$. If $\ell>0$, the Hardy integral diverges because it equals $\sph\int_0^\infty g(t)\,d t$. This proves the balance identity.

As a direct normalization check, take a smooth radial cutoff
$V_*(t)$ that vanishes for $t$ near $0$ and equals $1$ for
$t\ge t_1$. Then $g_*(t)=V_*(t)^p$ is $1$ throughout the
support of each translated cutoff difference near $t=T$ for all
large $T$, apart from the exponentially small region of increments
$|h|\gtrsim T-t_1$. Therefore the radial boundary contribution is
\[
 \begin{aligned}
 \lim_{T\to\infty}B_T
  &=\int_0^\infty
      \bigl[-h k(h)+h k(-h)\bigr]\,dh\\
  &=\kappa'(\beta)>0.
 \end{aligned}
\]
Thus the coefficient is $+\sph\kappa'(\beta)/p$, not its
negative, for a profile with unit pole trace. The compact-transition
part contributes to the Bregman remainder and need not vanish.
This check is also compatible with the density criterion in
Proposition~\ref{fractional-flux}. The distributional assumption
on $B_1$, rather than just on the punctured ball, excludes an
unaccounted point source at the pole; allowing an arbitrary singular
measure on the right would require a separate analysis of the
pairing with unbounded $u^{p-1}$ and is not asserted here.

To prove \eqref{angular-cost}, observe that for $x,y\ge0$,
\[
 D_p(y,x)\ge c_p|y-x|^p,
\]
because the ratio of the two sides is bounded below after normalizing
$x+y=1$ and using strict convexity (the cases $x=0$ or $y=0$ are
included by continuity). The angular kernel
is bounded below by a positive constant for \(1\le h\le2\). To estimate the
difference between \(V(t,\theta)\) and \(V(t,\theta')\), insert the
intermediate value \(V(t+h,\varphi)\). Raise the triangle inequality to the
power \(p\) and integrate in all three angular variables. More precisely, for each $h\in[1,2]$,
\[
 \begin{aligned}
 &\iint|V(t,\theta)-V(t,\theta')|^p
        d\omega(\theta)d\omega(\theta')\\
 &\quad\le 2^p\iint
       |V(t,\theta)-V(t+h,\varphi)|^p
         d\omega(\theta)d\omega(\varphi).
 \end{aligned}
\]
By Jensen's inequality, the left-hand side dominates
$\int|V(t,\theta)-\bar V(t)|^pd\omega(\theta)$.
Since $P(h,\theta,\varphi)$ is uniformly bounded below on
$[1,2]\times(\mathbb S^{N-1})^2$, integration in $h\in[1,2]$ and
$t>0$ proves the stated estimate. Any limiting trace with positive
angular oscillation would give an infinite integral on the left.
\end{proof}

\begin{suppcor}[Finiteness without an endpoint limit]\label{cor:notrace}
Assume the hypotheses on $u,V,f$ of Theorem~\ref{angular-flux}, except for the existence of the essential limit of $g$. If $f\geq0$, then
\[
 \mathcal R_{s,p}(V)<\infty\quad\Longleftrightarrow\quad
 \int_{B_1}u^{p-1}f\,d x<\infty.
\]
When these quantities are finite,
\[
 \left|\int_{B_1}u^{p-1}f\,d x-\mathcal R_{s,p}(V)\right|
 \leq (\sph/p)C_*\|V\|_\infty^p.
\]
For signed $f$, the assumption $u^{p-1}f\in L^1$ still implies that the remainder is finite. No absolute-integrability converse for signed $f$ is claimed.
\end{suppcor}
\begin{proof}
Use \eqref{eq:localbalance} and the uniform estimate \eqref{eq:fluxbound} for every $T$. The two
nonnegative cutoff integrals increase to their full values when $f\geq0$. If either full value is
finite, the uniform difference bound gives finiteness of the other and passes to their limits. In
the signed case, absolute integrability bounds the cutoff source pairing and hence the increasing
remainder.
\end{proof}
\noindent The corollary requires boundedness of the profile but not a limiting trace. In contrast,
the positive-data counterexamples in Theorem~\ref{fractional-counter} use unbounded profiles. The
existence of the trace is required in Theorem~\ref{angular-flux} to identify the numerical value of
the flux, not merely to compare finiteness of the two terms. The corollary does not assert that
finite remainder forces existence of a trace.

\subsection{Radial criteria and positive examples}
Let the radial profile be extended by zero to \(t<0\), and set
\[
 \Delta_h^2V(t)=V(t+h)+V(t-h)-2V(t).
\] Using $x=r\theta$ and $y=re^{-h}\varphi$ in the fractional
principal-value integral, subtract the ground-state equation
$(-\Delta)^sh_\beta=\kappa(\beta)r^{-2s}h_\beta$.
For a radial $V$, integration in the angular variable then gives
\[
 L_\lambda(h_\beta V)(r)
 =r^{-\beta-2s}\int_{\mathbb R}
    \bigl(V(t)-V(t+h)\bigr)k(h)\,dh,
 \qquad t=\log(1/r),
\]
where the integral is understood by symmetric exclusion around $h=0$.
Pairing $h$ with $-h$ separates this expression into a symmetric
second difference and an asymmetric first difference, giving
\begin{equation}\label{frac-paired}
 AV(t)=\int_0^\infty
 \big[-k(h)\Delta_h^2 V(t)
   +(k(-h)-k(h))(V(t)-V(t-h))\big]\,dh .
\end{equation}

\begin{proposition}[Radial density criterion]\label{fractional-flux}
Let \(V\ge0\) be of bounded variation on \(\R\). Assume that it vanishes for
\(t<0\) and that
\begin{equation}\label{fractional-second-diff}
 \int_0^1 h^{-1-2s}\|\Delta_h^2V\|_{L^1(\R)}\,dh<\infty.
\end{equation}
Then \(u=h_\beta V\in X_0^{s,2}(B_1)\), and the density
\[
 f=L_\lambda u=r^{-\beta-2s}AV(\log(1/r))
\]
is integrable.
Moreover, \(u^{p-1}f\) is absolutely integrable and \eqref{angular-balance}
holds with \(\ell=V(\infty)^p\). If \(s<1/2\), condition
\eqref{fractional-second-diff} holds automatically.
\end{proposition}

\begin{proof}
Since the profile has bounded variation,
\[
 \|V(\cdot+h)-V\|_1\le |h|\operatorname{Var}(V).
\]
By
Lemma~\ref{kernel-moment}, the second term in \eqref{frac-paired} is
integrable in \((t,h)\). The first term is integrable by
\eqref{fractional-second-diff} and the exponential decay of the kernel. Also
\[
 \|\Delta_h^2V\|_1
 \le 2\|V(\cdot+h)-V\|_1
 \le2h\operatorname{Var}(V),
\]
so the integrand of \eqref{fractional-second-diff} is bounded by
$2\operatorname{Var}(V)h^{-2s}$. This is integrable at zero when
$s<1/2$, proving the last assertion.

To estimate the energy, order the radii so that \(r=|x|\le |y|\). The pairs
that contribute then satisfy \(r<1\). Split the difference of the products
into a difference of \(h_\beta\) and a difference of \(V\). The first gives a
squared contribution bounded by
\[
 C\|V\|_\infty^2\int_0^1 r^{d-1}dr.
\]
Here
we use the gradient of \(h_\beta\) when \(|x-y|\le r/2\) and the bound \(h_\beta(x)\) on
the difference otherwise. The term involving the profile has radial factor
\(e^{-dt}\le1\) and angular kernel \(e^{-\beta h}k(h)\) for \(h<0\). This
kernel decays exponentially at infinity and is of order \(|h|^{-1-2s}\) at
zero. By changing variables in the integrable first differences, we get
\[
 \|V(\cdot+h)-V\|_2^2=-\int V\Delta_h^2V
        \le\|V\|_\infty\|\Delta_h^2V\|_1 .
\]
For large \(|h|\), we use boundedness of the profile and the translation
estimate in \(L^1\) to bound the squared difference. Hence
\(u\in H^s(\R^N)\), and dilation and mollification give membership in the
space with zero exterior values.

We subtract the two ground-state forms on annular tests, with the same
logarithmic exclusion as above. In radial coordinates, the truncated
integrals in \eqref{frac-paired} converge in \(L^1(dt)\). Averaging the
tests over rotations gives the density formula on the punctured ball. The
density is integrable since \(N-\beta-2s>0\). As \(N>2s\), removing the pole
by cutoffs gives convergence in energy to any smooth compactly supported
test. Dominated convergence on the density side then gives the equation on
the whole ball. Finally, \(p\beta=N-2s\) gives
\[
 \int |u^{p-1}f|=\sph\int_0^\infty V^{p-1}|AV|\,dt<\infty.
\]
We apply
Theorem~\ref{angular-flux} to obtain the identity.
\end{proof}

\begin{theorem}[Finite fractional remainder and divergent pairing]\label{fractional-counter}
Let \(0\le\alpha<s\), \(\max(\alpha,2s-1)<\eta<s\) and
\(C\ge\max(1,\alpha/d)\), where \(d\) is as in Lemma~\ref{kernel-moment}.
For sufficiently small \(\delta>0\), define
\[
 V(t)=
 \begin{cases}
  0,&t\le0,\\
  (C+\delta)^\alpha(t/\delta)^\eta,&0<t<\delta,\\
  (C+t)^\alpha,&t\ge\delta .
 \end{cases}
\]
Then \(u=h_\beta V\) is the variational solution for the finite positive density
\(f=L_\lambda u\). If \(\alpha=0\), the pole flux is
\(\sph\kappa'(\beta)/p\). If \(\alpha>0\), the source pairing is infinite
and the remainder is finite if and only if \(p\alpha<1\). Thus, for every
\(s\in(0,1)\) there is a positive source density for which the remainder is
finite and the pairing is infinite.
\end{theorem}

For each $0<s<1$ and $p>2$, the parameter restrictions are
nonempty. In fact choose
\[
 0<\alpha<\min\{s,1/p\},\qquad
 \max\{\alpha,2s-1\}<\eta<s,
\]
which is possible because $2s-1<s$ for $s<1$; then choose
$C\ge\max(1,\alpha/d)$ and finally sufficiently small
$\delta>0$. This makes $p\alpha<1$ while preserving
$\alpha>0$, as required for the divergent pairing.

\begin{proof}
We first prove positivity near the boundary using the function
\(B(t)=t_+^\eta\). The asymmetric term in \eqref{frac-paired} is nonnegative
for \(B\). After multiplication by \(t^{2s-\eta}\), the symmetric term
converges as \(t\downarrow0\) to \(c_{N,s}I_{\eta,s}\), where
\begin{align*}
 I_{\eta,s}
 &=\int_0^\infty
   [2-(1+z)^\eta-(1-z)_+^\eta]z^{-1-2s}\,dz\\
 &=\frac{\eta}{2s}\int_0^1 z^{-2s}
   [(1-z)^{\eta-1}-(1-z)^{2s-\eta-1}]\,dz>0 .
\end{align*}
To take the limit, use \(h=tz\), the bound \(k(h)\le Ch^{-1-2s}\) and the
cancellation of the terms of first order as \(z\downarrow0\). Integration by
parts and the change of variables \(z\mapsto z/(1+z)\) in the forward
increment give the second line. The two lower cutoffs differ by an error
\(O(\epsilon^{2-2s})\). For $0<z<1$, the integrand in the second line is positive since
$\eta-1<2s-\eta-1$ and $1-z<1$; thus $I_{\eta,s}>0$.
The rescaled integral is dominated at $z=0$ by $z^{1-2s}$ thanks to
second-order cancellation, and at infinity by
$z^{\eta-1-2s}$, which is integrable because $\eta<s<2s$.
Consequently the rescaled expression tends to the strictly positive
limit $c_{N,s}I_{\eta,s}$, and there exists $\delta_0>0$ such that
$AB(t)>0$ for every $0<t<\delta_0$. Choose $0<\delta<\delta_0$. The choice of $C$ and $\eta>\alpha$ imply
$\alpha/(C+t)\le\eta/t$. Therefore the function
\[
 t\longmapsto
 \log\left(\frac{(C+\delta)^\alpha(t/\delta)^\eta}
                  {(C+t)^\alpha}\right)
\]
is nondecreasing on $(0,\delta)$ and equals zero at $t=\delta$;
this also gives $V(t)\le W(t)$ there. In particular, we have
\[
 V\le (C+\delta)^\alpha\delta^{-\eta}B
\]
everywhere with equality on \((0,\delta)\). It follows that \(AV>0\) there.

For \(t>\delta\), compare with the upper bound \(W(t)=(C+t)^\alpha\) for
\(t>0\), extended by zero on the negative half-line. At such \(t\) we have
\(AV\ge AW\) because \(W=V\) there. If \(0<h<t\), the backward increment of
\(W\) is at least the forward increment by concavity, and \(k(-h)>k(h)\). If $h>t$, then $W(t-h)=0$. As $C\ge\alpha/d$ and
$\log(1+z)\le z$,
\[
 \frac{W(t+h)}{W(t)}
  =\left(1+\frac{h}{C+t}\right)^\alpha
  \le e^{\alpha h/(C+t)}
  \le e^{\alpha h/C}\le e^{dh}
  =\frac{k(-h)}{k(h)}.
\]
Thus $k(-h)W(t)-k(h)W(t+h)\ge0$ also in this case.
The integrand of the paired
integral is nonnegative in the first case and positive in the second case.
Thus $AV>0$ on the two open intervals $(0,\delta)$ and
$(\delta,\infty)$ wherever the paired integral exists. At the
single junction $t=\delta$, $V$ is continuous but its derivative
may have a downward jump. For a fractional operator the value of the
density at that individual height is immaterial; the estimate on
$\|\Delta_h^2V\|_1$ below shows that the operator has an
$L^1_{\mathrm{loc}}$ density there, with no negative point mass
concentrated at the junction. Consequently the source is nonnegative
as a measure and positive almost everywhere.

The derivative \(V'\) is nonincreasing on \((0,\infty)\), with a possible
downward jump at \(\delta\). Extend the derivative by zero. Its \(L^1\)
difference under translation is at most \(2\int_0^h V'(t)dt=O(h^\eta)\) as
\(h\downarrow0\). Integrating gives
\[
 \|\Delta_h^2V\|_{L^1}=O(h^{1+\eta}).
\]
Since \(\eta>2s-1\), this gives local \(L^1\) convergence in
\eqref{frac-paired}. The first differences are locally integrable with an
\(O(h)\) bound. The exponential tails of the kernel control the polynomial
growth. If $\alpha=0$, the profile is bounded, nondecreasing, and eventually
equal to one; its total variation is one. The condition
$\eta>2s-1$ gives \eqref{fractional-second-diff}, so
Proposition~\ref{fractional-flux} applies. The pole trace is
$V(\infty)=1$, and \eqref{angular-balance} gives the pole term
$\sph\kappa'(\beta)/p$. The Bregman remainder is \emph{not} zero:
$V$ varies on $(0,\delta)$ and vanishes for $t<0$, so strictly
positive pairs occur in \eqref{frac-rem}. This distinction matters even
though the tail profile is constant.

If \(\alpha>0\), the estimates near the pole are
\[
 |u(r)|\le C r^{-\beta}(\log(e/r))^\alpha,\qquad
 |u'(r)|\le C r^{-1-\beta}(\log(e/r))^\alpha.
\]
For pairs ordered by radius, split the fractional energy integral at
separation \(r/2\). The estimates give a finite bound by integration of
\(r^{d-1}(\log(e/r))^{2\alpha}\). Away from the pole, truncate the profile
at a large height and use the preceding energy estimate for the bounded
profile. Thus \(u\in X_0^{s,2}(B_1)\). At infinity, the paired operator has the asymptotic
\[
 AV(t)=\alpha m_\beta(C+t)^{\alpha-1}
             +O((C+t)^{\alpha-2})+o((C+t)^{\alpha-1}).
\]
In particular, $AV(t)\sim\alpha m_\beta(C+t)^{\alpha-1}$.
Indeed, for \(0<h<t/2\) the second difference is \(O(t^{\alpha-2}h^2)\). For
each fixed increment, the backward difference divided by \(V'(t)\) tends to
\(h\) and is bounded by a constant times \(h\) in this range. For $h>t/2$, the kernel bounds
$k(h)=O(e^{-(N-\beta)h})$ and
$k(-h)=O(e^{-(\beta+2s)h})$ dominate the polynomial growth of
$V(t\pm h)$, so both the second-difference and first-difference
tails are $O(e^{-ct}t^C)$ for suitable $c>0,C<\infty$.
They are therefore $o(V'(t))$, and the first-moment term gives
the stated asymptotic. The density formula then gives \(f\in L^1(B_1)\), and the same
cutoffs give the distributional equation across the pole. It follows that
\(f\) belongs to the dual of the energy space, and by coercivity \(u\) is
its unique variational solution. The pairing diverges like
\[
 \int^\infty t^{p\alpha-1}dt.
\]

For the remainder we use
\[
 D_p(y,x)\le C_p(x+y)^{p-2}|x-y|^2.
\]
On bounded
\(t\)-intervals, the bounds at the joins give integrability for small
increments. Exponential decay controls the region \(|h|>1\) with
\(|t|\le2|h|+2\). For large \(t\) and \(|h|\le t/2\), the factor multiplying
the kernel is bounded by \(C t^{p\alpha-2}h^2\). For a lower bound, restrict $h$ to $[1,2]$ and use the Taylor
expansion, uniform there for large $t$,
\[
 \begin{aligned}
 V(t+h)-V(t)&=\alpha h(C+t)^{\alpha-1}
                     +O(t^{\alpha-2}),\\
 D_p(V(t+h),V(t))
  &=\tfrac12p(p-1)\alpha^2h^2(C+t)^{p\alpha-2}
                     (1+o(1)).
 \end{aligned}
\]
The kernel is positive and integrable on $[1,2]$, so this yields a
positive multiple of $t^{p\alpha-2}$ for the tail remainder. The upper
bound follows by the quadratic Bregman estimate for $|h|\le t/2$ and
the exponential tails for $|h|>t/2$. Consequently, for $\alpha>0$,
\[
 \mathcal R_{s,p}(V)<\infty
 \quad\Longleftrightarrow\quad
 \int^{\infty}t^{p\alpha-2}\,dt<\infty
 \quad\Longleftrightarrow\quad p\alpha<1.
\]
The local contribution, including the join at $\delta$, is finite by
$\eta>2s-1$ and the preceding second-difference estimates.
\end{proof}

\subsection{The local balance and a tail criterion}
For the local equation let \(N\ge3\) and set \(a=(N-2)/p\), \(a_*=(p-1)a\),
\(h_0=a_*-a>0\) and \(c_p=4(p-1)/p^2\). Consider \(L=-\Delta-aa_*|x|^{-2}\)
on \(B_1\). Let \(\mu\) be a positive radial measure with zero mass at the
pole and support in \(0<|x|\le e^{-A}\), and put
\[
 \begin{aligned}
 \nu&=\sph^{-1}(|x|)_\#\mu,\\
 S&=(r\mapsto\log(1/r))_\#(r^{-a_*}\nu).
 \end{aligned}
\]
Assume
\[
 \int e^{-h_0s}dS(s)<\infty.
\]
We take the Green solution with zero
boundary values and the less singular homogeneous behavior at the pole. Its
profile \(Z(t)=r^a u(r)\) is given by
\begin{equation}\label{drift-kernel}
 Z(t)=\int\!\int_0^{\min(t,s)} e^{-h_0(s-v)}\,dv\,dS(s).
\end{equation}
To see the reduction, write $r=e^{-t}$ and
$u(r)=r^{-a}Z(t)$ in the radial Laplacian. The indicial roots of
$-\Delta-aa_*r^{-2}$ are $a$ and $a_*$; cancellation of the
zeroth-order term leaves the distributional ODE
\[
 -Z''+h_0Z'=S\quad\text{on }(0,\infty),
 \qquad Z(0)=0.
\]
For a unit atom at $s$, the kernel in \eqref{drift-kernel} has derivative
$e^{-h_0(s-t)}$ for $t<s$ and zero for $t>s$; its derivative has the
required jump of size $-1$ at $s$. The kernel is constant for $t>s$,
which selects the less singular homogeneous behavior at the pole. Equivalently, the same kernel follows from matching the two
homogeneous powers $r^{-a}$ and $r^{-a_*}$ with the outer boundary
condition. For \(t\ge A\), set \(C(t)=S([A,t])\) and
\[
 T(t)=\int_{(t,\infty)}e^{-h_0(s-t)}\,dS(s).
\]

\begin{proposition}[Local pole balance]\label{flux}
The local remainder
\[
 \mathcal D_p=\int_{B_1}|x|^{-2b_H}
       |\nabla(|x|^{b_H} u^{p/2})|^2\,dx,
       \qquad b_H=(N-2)/2 ,
\]
is finite if and only if
\[
 \int_A^\infty(C(t)+T(t))^{p-2}T(t)^2dt<\infty.
\]
The source pairing is finite if and only if \(S\) has finite total mass.
This is also equivalent to \(Z(\infty)<\infty\). In this case
\begin{equation}\label{flux-eq}
 \int_{B_1}u^{p-1}\,d\mu
  =c_p\mathcal D_p+\frac{h_0}{p}\sph Z(\infty)^p,\qquad
 Z(\infty)=h_0^{-1}\int(1-e^{-h_0s})\,dS(s).
\end{equation}
If \(S\ne0\), finiteness of the remainder implies \(u\in H_0^1(B_1)\) and
\(\mu\in H^{-1}(B_1)\).
\end{proposition}

\begin{proof}
It follows from \eqref{drift-kernel} that \(Z'=T\) almost everywhere and
that
\[
 h_0 Z(t)=\int_{[A,t]}(1-e^{-h_0s})\,dS(s)+(1-e^{-h_0t})T(t)
 \quad(t\ge A).
\]
Hence \(Z\asymp C+T\) on this interval. By polar integration,
\[
 \mathcal D_p=(p^2\sph/4)\int_0^\infty Z^{p-2}|Z'|^2dt.
\]
The explicit
kernel gives finiteness of the integral below \(A\). The source pairing
divided by \(\sph\) equals \(\int Z^{p-1}dS\). If \(S\) is nonzero, then
\(Z\ge Z(A)>0\) on \([A,\infty)\), and the kernel increases to
\((1-e^{-h_0s})/h_0\). This gives the conditions on the mass and the trace.
For $L>A$ with $S(\{L\})=0$, the derivative $Z'$ is of locally
bounded variation and $-dZ'+h_0Z'\,dt=dS$. As $Z$ is continuous and
locally absolutely continuous, Stieltjes integration by parts gives
\[
 \begin{aligned}
 \int_{[0,L]}Z^{p-1}\,dS
 &= -Z(L)^{p-1}Z'(L)+Z(0)^{p-1}Z'(0+)\\
 &\quad +(p-1)\int_0^LZ^{p-2}(Z')^2dt
       +h_0\int_0^LZ^{p-1}Z'\,dt.
 \end{aligned}
\]
Since $Z(0)=0$, the initial boundary term vanishes, and the last integral
is $Z(L)^p/p$. Thus
\[
 \int_{[A,L]} Z^{p-1}dS
 = (p-1)\int_0^L Z^{p-2}(Z')^2dt
      +\frac{h_0}{p}Z(L)^p-Z(L)^{p-1}Z'(L).
\]
If $S$ has finite mass, then
\[
 Z'(L)=\int_{(L,\infty)}e^{-h_0(s-L)}\,dS(s)
  \le S((L,\infty))\longrightarrow0.
\]
The boundedness of $Z(L)$ then eliminates the remaining boundary term.
Monotone convergence applies to the nonnegative energy and source
integrals, giving \eqref{flux-eq}. When the remainder is finite, we have \(Z'\in L^2\) and
\(Z(t)=O(\sqrt t)\). Thus
\[
 \int_0^\infty e^{-h_0t}|aZ+Z'|^2dt<\infty,
\]
which implies \(u\in H_0^1(B_1)\). The pole has zero capacity, and the
radial flux in the original variables \(-e^{-a_*t}(aZ+Z')\) tends to zero by
\[
 Z'(t)\le e^{h_0t}\int e^{-h_0s}dS
\]
off the atoms. Thus the equation holds
across the pole with zero mass there. By the usual Hardy inequality,
\(\mu=Lu\in H^{-1}\).
\end{proof}

\begin{remark}
Suppose the tail of \(S\) has density asymptotic to
\(t^{\alpha-1}(\log t)^{-\eta}\) with \(\alpha>0\). The exponential convolution gives
\[
 \begin{aligned}
 T(t)&\sim\frac{t^{\alpha-1}(\log t)^{-\eta}}{h_0},\\
 Z(t)&\sim\frac{t^\alpha(\log t)^{-\eta}}{h_0\alpha}.
 \end{aligned}
\] For this family, the
pairing is infinite. The remainder is finite for \(\alpha<1/p\) or, at
equality, for \(\eta>1/p\). When the pairing is finite and nonzero, the pole
trace is positive and the usual Hardy term for \(u^{p/2}\) is infinite. By
the gamma recurrence we have \(\kappa_1(a)=a(N-2-a)\), and
\(\kappa_s'((N-2s)/p)\) tends to \(h_0\) as \(s\uparrow1\).

In \eqref{drift-kernel}, fix a finite measure \(S\) and let
\(h_0\downarrow0\). The kernels increase to \(\min(t,s)\), and the
derivatives increase to \(S((t,\infty))\). By monotone and dominated
convergence, convergence in \(\Hrad\) holds if and only if this tail belongs to
\(L^2(dt)\). The pole traces increase to \(\int s\,dS(s)\), which may be
infinite even when the energy is finite.
\end{remark}
\section{Heat reactions at a time impulse}\label{sec-impulse}

We now consider concentration of the input at one time. The jump from the
incoming to the outgoing state is the reaction on that time slice. We show
how this jump depends on the approximation of the input.

Let \(D=\T^d\) and \(Q=[0,T]\times D\) with \(T>0\). Let
\(S(t)=e^{t\Delta}\) be the periodic heat semigroup with kernel \(p_t\). In
Sections~\ref{sec-impulse} and~\ref{sec-relaxation}, we assume that
\(F\in C^2(0,\infty)\) is positive and nonincreasing, with \(F(r)\to0\) as
\(r\to\infty\). For a smooth nonnegative coefficient \(a\) and bounded
nonnegative initial data \(u_0\), write \(R(a,u_0)\) for the reaction
measure \(aF(u)\,dt\,dx\) associated with
\begin{equation}\label{heat}
 \partial_tu-\Delta u=aF(u),\qquad u(0)=u_0 .
\end{equation}
If $u_0\ge e>0$, the reaction law is locally Lipschitz on the solution
range. The maximum principle gives
\[
 e\le u(t,x)\le\|u_0\|_\infty+F(e)\int_0^T\|a(s)\|_\infty\,ds,
\]
so local fixed-point iteration extends to $[0,T]$ and yields the mild
identity
\[
 u(t)=S(t)u_0+\int_0^t S(t-s)\bigl(a(s)F(u(s))\bigr)\,ds.
\]
If $u_0\le v_0$ and $u,v$ solve with the same input, the positive-part
test for $u-v$ and monotonicity of $F$ give $u\le v$. More generally,
comparison with $u+c$, using $F(u+c)\le F(u)$, shows that replacing
$u_0$ by $u_0+c$ changes the solution by at most $c$ everywhere.

For arbitrary nonnegative $u_0$, solve from $u_0+e$ and decrease
$e\downarrow0$. The solutions decrease uniformly because of this
comparison estimate. Their reactions $aF(u^{(e)})$ increase pointwise
and their total masses are bounded: integrating the mild equation over
$D$ yields
\[
 \int_QaF(u^{(e)})\,dt\,dx
 =\int_Du^{(e)}(T)\,dx-\int_D(u_0+e)\,dx.
\]
For $0<e'<e$, the reaction measures are ordered in the opposite
direction from the solutions, and the final-state comparison bounds
the difference of their masses by $(e-e')V_D$. Therefore the
reactions converge as finite measures, and monotone convergence in
Duhamel's formula proves the mild equation for the limiting state.
We set $R(a,u_0)=aF(u)\,dt\,dx$ with the convention described below. We take the reaction to be zero at points
where the coefficient is zero, and \(F(0)\) denotes the extended right
limit.

Set
\[
 \Phi(r)=\int_0^r\frac{dz}{F(z)},\qquad K=\Phi^{-1},\qquad
 J_c(r)=K(\Phi(r)+c),\quad c,r\ge0 .
\]
As $F>0$ and $F'\le0$,
\[
 \Phi'(r)=\frac1{F(r)}>0,\qquad
 \Phi''(r)=-\frac{F'(r)}{F(r)^2}\ge0.
\]
Thus $\Phi$ is increasing and convex, and
$\Phi(r)/r\to\infty$ since $F(r)\to0$ as $r\to\infty$.
Its inverse $K$ is increasing and concave with $K(0)=0$. Inverting
this superlinear growth gives $K(c)/c\to0$ as $c\to\infty$.
Moreover $J_c$ is increasing and nonexpansive in its initial state:
$J_c'(r)=F(J_c(r))/F(r)\le1$ wherever the derivative exists. For a
pulse with a given spatial density, the scalar map \(J_c\) gives the
outgoing value after input mass \(c\) when diffusion is omitted.

\begin{theorem}[Exact range at one time]\label{atom-range}
Fix \(t_*\in[0,T]\) and let \(\sigma=f\,dx+\sigma_{\rm s}\) be the Lebesgue
decomposition of a finite positive spatial measure. Set \(b=S(t_*)u_0\). For
smooth nonnegative inputs converging narrowly on \(Q\) to
\(\delta_{t_*}\otimes\sigma\), the possible limiting reactions are precisely
\begin{equation}\label{range-bound}
 \{\,\delta_{t_*}\otimes r\,dx:
       0\le r(x)\le J_{f(x)}(b(x))-b(x)\quad\text{a.e.}\,\}.
\end{equation}
Every such reaction is attained along a full sequence. For
\(F(z)=z^{-\gamma}\) the upper bound is
\((b^{1+\gamma}+(1+\gamma)f)^{1/(1+\gamma)}-b\). The range consists of a
single measure if and only if \(f=0\) almost everywhere.
\end{theorem}

\subsection{Reaction estimates and sets of zero capacity}
Write \(V_D=|D|\) and define
\[
 \omega_F(\delta)=\inf_{e>0}\big(2eV_D+2F(e)\delta\big),
 \qquad \lim_{\delta\downarrow0}\omega_F(\delta)=0 .
\]
The infimum is over positive lifts $e>0$. Given $\epsilon>0$,
first choose $e$ with $2eV_D<\epsilon/2$ and then choose
$\delta$ so that $2F(e)\delta<\epsilon/2$; this proves
the asserted modulus at zero.
We first prove two comparison estimates for the reaction.

\begin{lemma}[Total variation of the reaction]\label{reaction-tv}
For smooth nonnegative inputs and bounded nonnegative initial data, we have
\begin{equation}\label{tv-pair}
 \|R(a,u_0)-R(b,v_0)\|_{\TV}
 \le\|u_0-v_0\|_1+\omega_F(\|a-b\|_{L^1(Q)}) .
\end{equation}
If the initial data are the same, the reaction with an added input satisfies
\begin{equation}\label{tv-add}
 \|R(b+c,u_0)-R(b,u_0)\|_{\TV}\le 2 R(c,u_0)(Q).
\end{equation}
\end{lemma}

\begin{proof}
First fix the input $a$ and assume $u_0\le v_0$. Write $u,v$ for
the resulting solutions. Comparison gives $u\le v$, while $F(u)\ge F(v)$,
so $R(a,u_0)-R(a,v_0)$ is a nonnegative measure. Subtracting their
mass balances gives
\[
 \begin{aligned}
 0\le \bigl(R(a,u_0)-R(a,v_0)\bigr)(Q)
 &=\int_D(v_0-u_0)\,dx-\int_D(v(T)-u(T))\,dx\\
 &\le \|v_0-u_0\|_1.
 \end{aligned}
\]
For unordered data, compare each solution with that for
$u_0\vee v_0$, then add the two resulting estimates. Their total
initial-data differences sum to $\|u_0-v_0\|_1$. For the input variation, first lift the common initial datum by
$e>0$, and assume $a\le b$. Write $u_a\le u_b$ for the solutions.
The signed measure difference has density
\[
 bF(u_b)-aF(u_a)
 =(b-a)F(u_b)-a\bigl(F(u_a)-F(u_b)\bigr).
\]
Its positive part is at most $F(e)(b-a)$, while its total mass is
nonnegative by comparison of the final states and equality of the
initial masses. For any signed measure $\nu$, the identity
$\|\nu\|_{\TV}=2\nu^+(Q)-\nu(Q)$ therefore yields
\[
 \|R(b,u_0+e)-R(a,u_0+e)\|_{\TV}
 \le2F(e)\|b-a\|_{L^1(Q)}.
\] For general inputs, compare both with \(a\vee b\). This input
is bounded, and the mild construction and comparison also apply to it.
Removing the two lifts changes the reactions by at most \(2eV_D\) in total
variation. Taking the infimum over the lift gives \eqref{tv-pair}.

For \eqref{tv-add}, denote by $U,v,w$ the solutions for $b+c$, $b$,
and $c$ with the same initial data. Comparison gives $U\ge v,w$.
Decompose the reaction difference into a nonnegative added-source term
and a nonpositive background term,
\[
 R(b+c)-R(b)=cF(U)\,dt\,dx
                 -b\bigl(F(v)-F(U)\bigr)\,dt\,dx.
\]
The first term has mass $A\le R(c,u_0)(Q)$ by $U\ge w$. The second
has mass $B\ge0$. The equation for $U-v$ gives
$\int_D(U(T)-v(T))\,dx=A-B\ge0$, hence $B\le A$. It follows that
$\|R(b+c)-R(b)\|_{\TV}\le A+B\le 2R(c,u_0)(Q)$.
These estimates initially hold for positive lifts; the initial-data
bound proved above permits $e\downarrow0$.
\end{proof}

A Borel set $B\subset(0,T)\times D$ is called
\emph{parabolically polar} if every compact set $E\subset B$ has
zero capacity in the sense below; the notion is local to the open
time cylinder, so endpoint slices are dealt with separately. For a
compact set \(E\subset(0,T)\times D\), define the parabolic capacity
by
\begin{equation}\label{par-cap}
 \ecap_{\mathrm{par}}(E)=\inf_\varphi
  \big(\|\varphi\|_{L^2_tH^1_x}^2+
              \|\partial_t\varphi\|_{L^2_tH^{-1}_x}^2\big),
\end{equation}
where the infimum is over nonnegative
\(\varphi\in C_c^\infty((0,T)\times D)\) that are at least one in a
neighborhood of \(E\). This is the cutoff definition of
capacity for the \(L^2\) parabolic energy space~\cite{DPP,PPP}. In the
regularized local inverse-power problem with zero initial data and
\(\gamma\ge1\), Das~\cite[Theorem~2.15]{Das25} proved that the corresponding
power of the solution converges weakly to zero in energy for polar inputs.
Here we estimate the reaction mass for the decreasing factors specified
above.

\begin{proposition}[Parts of the input that always disappear]
\label{time-reduction}
Suppose that \(a_n\,dt\,dx\to\Lambda\) narrowly. The reactions are
relatively compact as finite measures on \(Q\). If \(R\) is a limit, then
\(R\ll\Lambda\) and for every time the spatial slice \(R(\{t\}\times dx)\)
is absolutely continuous with respect to \(dx\). Moreover, \(R\) vanishes on
polar Borel sets in the open time cylinder. If \(\Lambda\) is supported on
the union of such a polar set and countably many sets \(\{t_j\}\times N_j\)
with \(|N_j|=0\), then \(R(a_n,u_0)(Q)\to0\) and \(u_n(t)\to S(t)u_0\)
uniformly in time in \(L^1(D)\).
\end{proposition}

\begin{proof}
First assume that \(u_0\ge e>0\). Then \(0\le R_n\le F(e)a_n\,dt\,dx\), so
the masses are bounded and every limit satisfies \(R\le F(e)\Lambda\). Indeed, for positive $u_n$ the chain rule gives the pointwise
inequality
\[
 \begin{aligned}
 (\partial_t-\Delta)\Phi(u_n)
 &=\Phi'(u_n)a_nF(u_n)-\Phi''(u_n)|\nabla u_n|^2\\
 &=a_n-\Phi''(u_n)|\nabla u_n|^2\le a_n.
 \end{aligned}
\]
The positive-lift approximation justifies its integrated form when
the initial data vanish.
Integrating yields
\begin{equation}\label{entropy}
 \sup_{n,t}\int_D\Phi(u_n(t))\,dx
 \le\int_D\Phi(u_0)\,dx+\sup_n\int_Qa_n\,dt\,dx .
\end{equation}
For bounded nonsmooth data, smooth in space with the same positive lower
bound and use contraction for the fixed input. Monotonicity and Fatou's
lemma then give the estimate for data that may vanish.

To prove absolute continuity on a time slice, extend the solutions past
\(T\) by the free heat flow and extend \(R_n\) by zero. Jensen's inequality
preserves \eqref{entropy}. For \(t'>t\), take the positive Duhamel integral
at \(t'\) with a continuous time cutoff equal to one at \(t\) and supported
strictly before \(t'\). Along the chosen subsequence these integrals
converge uniformly in space since the kernels form a compact family of
continuous tests. The limit dominates \(S(t'-t)R(\{t\}\times\cdot)\).
Therefore its entropy is bounded by \eqref{entropy}. As $t'\downarrow t$, the heat potentials converge narrowly to the
spatial slice $R(\{t\}\times\cdot)$. A family with uniformly bounded
$\int_D\Phi$ is uniformly integrable by the de la Vallee--Poussin
criterion, because $\Phi(r)/r\to\infty$. Its narrow measure limit
must therefore be absolutely continuous with respect to $dx$.
For $t=0$ use $t'\downarrow0$ with only forward times; for $t=T$
use the extension by free heat flow. Thus the same argument includes
both endpoints.

For the assertion about capacity, let \(H_k\) be smooth, nondecreasing and
concave, with \(H_k(0)=0\), \(H_k'=1\) on \([0,k]\), \(0\le H_k'\le1\), and
\(H_k'=0\) on \([2k,\infty)\). Set \(L_n=H_k'(u_n)R_n\). Then
\[
 0\le L_n\le(\partial_t-\Delta)H_k(u_n),\qquad
 0\le R_n-L_n\le F(k)a_n\,dt\,dx .
\]
For the capacity estimate we test \eqref{heat} with
$T_{2k}(u_n)=\min(u_n,2k)$. If
$G_{2k}(r)=\int_0^rT_{2k}(z)\,dz$, the time derivative
integrates to
$\int_DG_{2k}(u_n(T))-
  \int_DG_{2k}(u_0)$. Since $0\le G_{2k}(r)\le2kr$
and $T_{2k}(u_n)\le2k$, integration by parts yields
\[
 \int_Q|\nabla T_{2k}(u_n)|^2\,dt\,dx
 \le2k\bigl(\|u_0\|_1+R_n(Q)\bigr).
\]
The chain rule gives
$\nabla H_k(u_n)=H_k'(u_n)\nabla u_n$ with
$H_k'$ supported in $[0,2k]$ and $|H_k'|\le1$.
Hence $|\nabla H_k(u_n)|\le|\nabla T_{2k}(u_n)|$
a.e. Since $0\le H_k\le2k$, these estimates yield a uniform
$L^2_tH^1_x$ bound for $H_k(u_n)$. Test the
distributional inequality with a nonnegative cutoff and integrate by parts.
The integral \(\int\varphi\,dL_n\) is bounded by
\(C_{k,e}(\|\partial_t\varphi\|_{L^2_tH^{-1}_x}
   +\|\varphi\|_{L^2_tH^1_x})\).
For fixed \(k\), pass to limits of the two parts of the reaction. On a
compact set of zero capacity the first part vanishes, and the remaining part
is bounded by \(F(k)\Lambda\). Let \(k\to\infty\) and use inner regularity.

For initial data that may vanish, apply the preceding estimates
first to $u_0+e$, where all chain-rule manipulations are classical.
For each fixed smooth input the comparison and mass-balance argument
in Lemma~\ref{reaction-tv} gives
$\|R(a,u_0+e)-R(a,u_0)\|_{\TV}\le eV_D$.
This is uniform in the approximation index and can therefore be
passed through narrow subsequences. Thus the mass bound and
the three absolute-continuity properties pass to the limit. Take a countable
sequence \(e\downarrow0\), extract the limits for the positive lifts by a
diagonal argument, and use the total variation bound. Under the assumption
on the support of the input this yields \(R=0\) for every subsequential
limit. Since the reactions are positive and the cylinder is compact, their
masses tend to zero. The formula of Duhamel and contraction of the heat
semigroup give the convergence of the solutions.
\end{proof}

\subsection{Proof of the impulse range}
\begin{proof}[Proof of Theorem~\ref{atom-range}]
By Proposition~\ref{time-reduction}, every limit has the form
\(\delta_{t_*}\otimes r\,dx\). Choose smooth time cutoffs supported on
relative intervals \(I_n\) shrinking to \(t_*\). By narrow convergence and a
diagonal choice, the cutoffs remove \(o(1)\) input mass. The reaction has
the same limit by \eqref{tv-pair}. Denote by \(u_n^-,u_n^+\) the incoming
and outgoing states at the endpoints of \(I_n\). The incoming states are
given by free heat evolution. They have a common supremum bound and converge
strongly in \(L^1\) to \(b\). Test
\[
 (\partial_t-\Delta)\Phi(u_n)\le a_n
\]
with a smooth spatial function \(\varphi\ge0\) and integrate over the
interval. We obtain
\[
 \int_D\varphi[\Phi(u_n^+)-\Phi(u_n^-)]\,dx
 \le\int_{I_n\times D}\varphi a_n\,dt\,dx
       +C\|\Delta\varphi\|_\infty |I_n| .
\]
Indeed, integration by parts of the diffusion term produces
$\int_{I_n\times D}\Phi(u_n)\Delta\varphi$, whose absolute value
is bounded by $\|\Delta\varphi\|_\infty|I_n|$ times the uniform
entropy mass from \eqref{entropy}. The chain rule is first applied
with an $e$-lift and then passed to the limit using monotone
convergence and the bound \eqref{entropy}. The outgoing states are
uniformly integrable. Integrating \eqref{heat} on the interval gives
\(u_n^+\rightharpoonup b+r\) weakly in \(L^1\). Indeed, for any smooth test
function the diffusion term is bounded by the length of the interval times a
uniform bound on the mass, with a constant depending on the test function.
Since the incoming states lie in a bounded range, their entropies converge
strongly in \(L^1\). Weak lower semicontinuity of the convex integral with
weight given by the test function yields
\[
 [\Phi(b+r)-\Phi(b)]\,dx\le\sigma.
\]
The absolutely continuous part gives
the upper bound in \eqref{range-bound}.

Conversely, set \(g=\Phi(b+r)-\Phi(b)\), so that \(0\le g\le f\). Take
smooth \(g_n\ge0\) converging to \(g\) in \(L^1\). Also choose smooth
\(b_n\ge0\) converging narrowly to \(\sigma-g\,dx\) with the same mass \(M\)
and supported on sets \(E_n\) with \(|E_n|\to0\). We obtain the second
sequence by replacing the remaining measure with discrete measures on fine
partitions and smoothing each atom in a sufficiently small ball. To quantify the scalar approximation, concavity and $K(0)=0$
imply subadditivity $K(x+y)\le K(x)+K(y)$ and hence
$|K(x+y)-K(x)|\le K(y)$ for $x,y\ge0$. Applying this with
$x=\Phi(b)+\min(c,d)$ and $y=|c-d|$ yields
\[
 |J_c(b)-J_d(b)|\le K(|c-d|).
\]
The scalar flow is also contractive in the initial state, since
$J_c'(b)=F(J_c(b))/F(b)\le1$ for $b>0$ and extends by continuity
to $b=0$.
Split the difference into the contribution of $g_n-g$ and that
of $b_n$, and apply the preceding bound twice. Jensen's inequality
for the concave function $K$ then gives
\begin{align*}
 \|J_{g_n+b_n}(b)-(b+r)\|_1
 &\le V_D K(\|g_n-g\|_1/V_D)
       +|E_n|K(M/|E_n|)\longrightarrow0 .
\end{align*}
The first term tends to zero by $L^1$ convergence of $g_n$, and
the second tends to zero because $M$ is fixed and
\[
 |E_n|K(M/|E_n|)=M\frac{K(M/|E_n|)}{M/|E_n|}\longrightarrow0.
\]
When $M=0$ this term is defined as zero. Thus no reaction survives
from the spatially concentrated part of the approximation.

Choose smooth positive functions \(b_n^*\to b\) in \(L^1\). For each \(n\),
the scalar paths \(J_{\theta(g_n+b_n)}(b_n^*)\), \(0\le\theta\le1\), have a
common bound \(C_n\) on their spatial Laplacians. Use the spatial density
\(g_n+b_n\) with a smooth time pulse of unit mass and width \(\ell_n\) at
\(t_*\), supported inside the time interval. Choose $\ell_n\downarrow0$ so that $C_n\ell_n\to0$. Let
$\theta(t)$ be the cumulative mass of the time cutoff, so that
$\theta'=q_n\ge0$ and $\int q_n=1$. The spatial scalar path
\[
 W_n(t,x)=J_{\theta(t)(g_n(x)+b_n(x))}(b_n^*(x))
\]
obeys $\partial_tW_n=q_n(g_n+b_n)F(W_n)$ and
$|\Delta W_n|\le C_n$. Let $t_n^-$ be the start of the pulse and set
$e_n(t)=C_n(t-t_n^-)$. Then $W_n+e_n$ is a supersolution, since
\[
 (\partial_t-\Delta)(W_n+e_n)
 =q_n(g_n+b_n)F(W_n)-\Delta W_n+C_n
 \ge q_n(g_n+b_n)F(W_n+e_n).
\]
Similarly, $(W_n-e_n)_+$ is a subsolution: where $W_n>e_n$ the
inequality is reversed because $F$ is nonincreasing, and the positive
part preserves the subsolution inequality. Both barriers agree with
$b_n^*$ at the initial endpoint. Comparison gives
\[
 \|u_n(t_n^+)-J_{g_n+b_n}(b_n^*)\|_\infty
 \le C_n\ell_n=o(1)
\]
for the pulse problem started from $b_n^*$. The
actual incoming states converge to \(b\) in \(L^1\). The heat equation for a
fixed input and the scalar map \(J_c\) are both contractive with respect to
the incoming data. Hence the actual outgoing states converge to \(b+r\) in
\(L^1\). For a single full realizing sequence, at stage $n$ make the three
errors (spatial approximation of $g$, concentrated-part reaction,
and diffusion during the pulse) smaller than $1/n$ by the successive
choices above. The input then converges narrowly to
$\delta_{t_*}\otimes\sigma$, while the outgoing state tends
to $b+r$ in $L^1(D)$. Integrating the equation across the shrinking
pulse identifies the reaction limit as
$\delta_{t_*}\otimes r\,dx$ without passing to any
further subsequence.
\end{proof}

\begin{corollary}[Complete disappearance can be attained]\label{evo2-diagonal}
For every finite positive measure \(\Lambda\) on \(Q\), there is a smooth
nonnegative approximation for which the reaction mass tends to zero.
\end{corollary}

\begin{proof}
Approximate \(\Lambda\) narrowly by finite positive sums of point masses at
interior times. For each fixed sum, smooth the point masses. The resulting
reactions tend to zero by Proposition~\ref{time-reduction}. Choose one
smoothing for each sum so that the reaction mass and the bounded-Lipschitz
distance from that sum both tend to zero.
\end{proof}
\section{Atomless and mixed heat inputs}\label{sec-relaxation}

We now consider inputs with atomless time marginal. We first use short-time
heat potentials to construct the reaction for Kato measures. This
construction will be used to characterize the reaction limits for data
exhausted by Kato measures.

\subsection{Canonical reactions for Kato measures}
We extend measures on \(Q\) by zero in time and set
\[
 k_\Lambda(h)=\sup_{t\in\R,\,x\in D}
    \int_{t-h<s<t}p_{t-s}(x-y)\,d\Lambda(s,y).
\]
The one-sided parabolic Kato condition is \(k_\Lambda(h)\to0\), see
\cite{evo2-BHJ} for the linear heat potential. A Kato measure cannot carry a time atom. Indeed, if
$\Lambda(\{t_0\}\times D)=m>0$, then for $0<\varepsilon<h$
we have, by positivity of the heat kernel and $\int_Dp_\varepsilon(x-y)dx=1$,
\[
 \sup_{x\in D}\int_Dp_\varepsilon(x-y)\,
               d\Lambda(\{t_0\},y)\ge m/V_D.
\]
Taking $t=t_0+\varepsilon$ in the definition of $k_\Lambda(h)$
forces $k_\Lambda(h)\ge m/V_D$ for every $h>0$, a contradiction.
The one-sided condition only measures potentials from preceding
times; it is precisely the direction needed in the mild equation. For any bounded integrand \(b\), define
\[
 H_\Lambda b(t,x)=\int_{s<t}p_{t-s}(x-y)b(s,y)d\Lambda(s,y).
\]

\begin{proposition}[Kato approximation]\label{evo2-kato}
Let \(\Lambda\) be a Kato measure and \(u_0\) continuous and nonnegative.
Suppose \(\Lambda_n\to\Lambda\) narrowly for smooth nonnegative approximants
with
\[
 \lim_{h\downarrow0}\sup_n k_{\Lambda_n}(h)=0.
\]
Then the solutions
converge uniformly on \(Q\) to the unique nonnegative continuous mild
solution with finite reaction given by
\[
 u(t)=S(t)u_0+H_\Lambda(F(u))(t),\qquad
 \mathcal R(\Lambda)=F(u)\Lambda.
\]
The reactions converge narrowly to \(\mathcal R(\Lambda)\), which is
equivalent to \(\Lambda\). Such approximations can be obtained by mollifying
in space and time with nonnegative kernels.
\end{proposition}

\begin{proof}
Mollification after extension by zero preserves the bound on the modulus,
since translations preserve it and restriction decreases it. The endpoint
time slices have zero mass. Thus restricting the smooth approximants to
\(Q\) preserves narrow convergence. For a continuous integrand \(b\),
multiply the heat kernel by a continuous cutoff for delays below \(h\). The
error is bounded uniformly by \(\|b\|_\infty k_\Lambda(h)\), with the
corresponding bound for \(\Lambda_n\). The kernels with cutoff form a
compact family of continuous test functions indexed by \((t,x)\). For each fixed $h>0$, the cutoff kernels form an equicontinuous,
uniformly bounded family indexed by $(t,x)\in Q$; a finite-net
argument and narrow convergence give uniform convergence of the
cutoff potentials. On the discarded region $0<t-s<h$ the error is at
most
\[
 \|b\|_\infty\bigl(k_\Lambda(h)+\sup_n k_{\Lambda_n}(h)\bigr).
\]
Letting $n\to\infty$ and then $h\downarrow0$ proves
$H_{\Lambda_n}b\to H_\Lambda b$ uniformly and shows that
$H_\Lambda b$ is continuous. The modulus at a fixed positive delay
and the total masses give a uniform bound for the operator norms.

First add a positive constant \(e>0\) to the initial data. Truncate the
nonlinearity outside the interval from \(e\) to
\(\|u_0+e\|_\infty+F(e)\sup_\sigma\|H_\sigma1\|_\infty\), where \(\sigma\)
ranges over the limit and the approximants. The resulting extension is
Lipschitz. Put $L_e=\operatorname{Lip}(F)$ for the truncated reaction factor.
Choose $h>0$ such that
\[
 q:=L_e\sup_\sigma k_\sigma(2h)<1,
\]
where $\sigma$ ranges over the fixed limit and the approximants.
For two candidate paths on one interval of length $h$, with the past
held fixed, the mild map is a contraction in the supremum norm with
factor at most $q$; the short-time integral has exactly the above
operator bound. At the first interval the fixed-inhomogeneity error
is $\|H_{\Lambda_n}(F(u))-H_\Lambda(F(u))\|_\infty\to0$ by the
preceding paragraph. The difference between the solutions is bounded
by this error divided by $1-q$. On each successive interval, the
already controlled past contributes at most its supremum error times
the uniform full-potential operator bound. Induction over finitely
many intervals covers $[0,T]$ and proves uniform convergence of the
lifted solutions.
Their reactions therefore converge narrowly.

For smooth inputs, comparison and mass balance show that changing the
positive lift by \(e\) changes the state by at most \(e\) in the supremum
norm and the reaction by at most \(e|D|\) in total variation. The same
estimates hold in the limits with fixed lift. As the lift decreases to zero,
the solutions converge uniformly. For fixed input, the reaction densities
increase and their masses stay bounded. Their limit is \(F(u)\Lambda\), and
monotone convergence gives the mild equation. The total variation estimate
also gives narrow convergence of the approximating reactions for the
original initial data. Since \(F\) is positive and \(u\) is bounded, the two
measures are equivalent. If \(F(0)=\infty\), finiteness implies that the
input gives zero mass to the zero set of the solution.

For uniqueness, let \(d=u-v\) for two continuous solutions and set
\(M=(F(u)-F(v))\Lambda\). This difference of reactions is atomless in time.
Let $P_j$ be convolution with an even, nonnegative smooth
approximate identity on the torus. After this symmetric spatial
smoothing, the difference satisfies
\[
 \partial_tP_jd=\Delta P_jd+P_jM.
\]
At each spatial point the time
dependence is continuous and of bounded variation. The chain rule for the
square and integration in space give
\[
 \|P_jd(t)\|_2^2/2\le\int_{[0,t]\times D}P_j^2d\,dM.
\]
As $j\to\infty$, the continuous solutions $u,v$ are uniformly
approximated by their spatial mollifications; finite reaction mass
allows passage in the integrals by dominated convergence after a
bounded truncation of the integrands. The remaining terms follow
from the $L^2$ convergence of the spatial smoothing. Monotonicity of
$F$ then gives
\[
 \|d(t)\|_2^2/2\le\int d\,dM\le0.
\]
\end{proof}

We say that \(\Xi\) is \emph{exhausted by Kato measures} if there are finite
positive Kato measures \(\Xi_j\) increasing to \(\Xi\). The same holds for
every submeasure, by multiplying each \(\Xi_j\) by the corresponding
Radon--Nikodym density. In particular, this includes every measure of the
form \(m(t,x)\,dx\,d\nu(t)\) with \(m\) integrable and \(\nu\) finite and
atomless. Truncating \(m\) at \(j\) bounds the modulus by
\(j\sup_t\nu((t-h,t))\to0\). We define \(\mathcal R(\Xi)\) for this class by
completion in total variation. For Kato data we mollify the measures with
the same kernels. For two Kato measures $\Gamma_1,\Gamma_2$, mollify both with the
same nonnegative kernel. Convolution does not increase the total
variation norm of their difference, and the corresponding reaction
measures converge narrowly by Proposition~\ref{evo2-kato}. Apply
\eqref{tv-pair} to the smooth inputs, take lower semicontinuity of
total variation for the narrow limits, and then let the mollifier
shrink. This yields
\begin{equation}\label{complete-reaction}
 \|\mathcal R(\Gamma_1)-\mathcal R(\Gamma_2)\|_{\TV}
       \le\omega_F(\|\Gamma_1-\Gamma_2\|_{\TV}) .
\end{equation}
The modulus $\omega_F$ tends to zero at zero and is nondecreasing.
Thus for any increasing Kato exhaustion $\Xi_j\uparrow\Xi$,
\[
 \|\mathcal R(\Xi_j)-\mathcal R(\Xi_k)\|_{\TV}
 \le\omega_F(\|\Xi_j-\Xi_k\|_{\TV})\longrightarrow0.
\]
The limit is independent of the exhaustion: compare two exhaustions
against each other by the same estimate and use their total-variation
convergence to $\Xi$. Define this limit to be $\mathcal R(\Xi)$.
Duhamel's formula and $L^1$ contraction of the heat semigroup show
that the corresponding states converge uniformly in time in $L^1(D)$.
Notice that exhaustion is genuinely weaker than the global Kato
property: truncations of an arbitrary $L^1$ space--time density are
Kato measures, although the original density need not be Kato.

\begin{suppprop}[Necessary properties of Kato exhaustions]\label{prop:kato}
Every finite measure exhausted by Kato measures has an atomless time marginal and vanishes on
parabolic polar Borel sets in the open cylinder. Every positive submeasure is again exhausted by
Kato measures. In a decomposition $\Lambda=\Xi+\Pi$ as in Theorem~\ref{evo2-range}, the two
components are mutually singular, and $\Xi$ is uniquely determined among decompositions of this
type.
\end{suppprop}
\begin{proof}
First let $\nu$ be Kato. As shown at the beginning of this subsection, its time marginal has no atoms.

The full forward potential of $\nu$ is bounded on $Q$. Indeed, choose a positive delay $h_0$ for
which the short-delay modulus is finite, and use boundedness of the heat kernel for delays at least
$h_0$, together with the finite mass, for the other part. The same argument gives a uniform bound
for smooth nonnegative mollifications $\nu_\eps$ of the zero extension, restricted to $Q$: their
masses are bounded, and convolution by a nonnegative probability kernel does not increase the
supremum of the short-delay modulus. There is no mass on the endpoint slices, so these restricted
mollifications converge narrowly to $\nu$.

Let $w_\eps$ be the smooth heat potential of $\nu_\eps$ with zero initial data. There is a uniform
bound $0\leq w_\eps\leq M$. Multiplication of its heat equation by $w_\eps$ and integration in space
and time give
\[
 \frac12\|w_\eps(T)\|_2^2+
 \int_0^T\|\nabla w_\eps\|_2^2\,d t
 =\int_Q w_\eps\,d\nu_\eps
 \leq M\nu_\eps(Q).
\]
The potential is therefore uniformly bounded in $L^2(0,T;H^1(D))$, including its $L^2$ term by the
supremum bound. For every nonnegative $\varphi\in C_c^\infty((0,T)\times D)$, the heat equation and
duality imply
\begin{align*}
 \int_Q\varphi\,d\nu
 &=\lim_{\eps\downarrow0}\left[-\int_0^T
 \langle\partial_t\varphi,w_\eps\rangle\,d t+
 \int_Q\nabla w_\eps\cdot\nabla\varphi\,d t\,d x\right]\\
 &\leq C\left(\|\partial_t\varphi\|_{L^2_tH^{-1}_x}
            +\|\varphi\|_{L^2_tH^1_x}\right).
\end{align*}
For a compact polar set $E$, apply this bound to the nonnegative cutoffs equal to at least one near
$E$ with arbitrarily small capacity norm. It follows that $\nu(E)=0$. Inner regularity gives the
conclusion for polar Borel sets.

If $\nu_j\uparrow\nu$ is a Kato exhaustion, the assertions about time atoms and polar sets follow by
monotone convergence. If $\eta\leq\nu$, write $\eta=b\nu$ with $0\leq b\leq1$; then
$b\nu_j\uparrow\eta$ and $k_{b\nu_j}\leq k_{\nu_j}$. This proves stability under submeasures.

Finally, let $B$ be the union of a polar Borel set and the countably many spatially null time slices
carrying $\Pi$ in Theorem~\ref{evo2-range}. By the preceding properties, $\Xi(B)=0$, whereas
$\Pi(Q\setminus B)=0$. Thus the components are mutually singular. If a second such decomposition has
carrier $B'$, both Kato-exhaustible components vanish on $B\cup B'$, and both residual components
vanish outside this union. Restriction of $\Lambda$ to the complement yields the same $\Xi$ in both
decompositions.
\end{proof}

\subsection{Reaction limits for atomless data}
For a finite positive measure \(R\) on \(Q\), define the heat potential with
initial datum by
\begin{equation}\label{heat-potential}
 v_R(t,x)=S(t)u_0(x)+
 \int_{[0,t)\times D}p_{t-s}(x-y)\,dR(s,y).
\end{equation}
The pointwise value may be \(+\infty\). Let \(\mathcal A(\Lambda)\) denote
the set of all narrow limits of the reactions with the same initial data and
with smooth nonnegative approximations to \(\Lambda\).

\begin{theorem}[Range for Kato exhaustions]\label{evo2-range}
Let \(u_0\) be continuous and nonnegative and suppose that
\(\Lambda=\Xi+\Pi\) is finite. Assume that \(\Xi\) is exhausted by Kato
measures and \(\Pi\) is concentrated on the union of a polar Borel set in
the open cylinder and countably many sets \(\{t_j\}\times N_j\) with
\(|N_j|=0\). Then
\begin{equation}\label{range-exhausted}
 \mathcal A(\Lambda)=\{\mathcal R(\Gamma):0\le\Gamma\le\Xi\}.
\end{equation}
Every reaction in this set is attained along a full sequence. For any
sequence whose reactions converge to \(R\), the solutions converge to
\(v_R\) uniformly in time in \(L^1(D)\).
\end{theorem}

To prove necessity, we first obtain a pointwise lower bound for the
solutions from heat potentials with a cutoff at short times.

\begin{lemma}[An inequality for limiting reactions]\label{potential-limit}
Suppose the reactions \(R_n\) converge narrowly to \(R\) and the inputs
converge narrowly to \(\Lambda\). If \(R\) has atomless time marginal, the
solutions converge to \(v_R\) uniformly in time in \(L^1(D)\). Moreover,
\begin{equation}\label{reaction-obstacle}
 R\le F(v_R)\Lambda\ \text{on }\{v_R>0\},\qquad
 R(\{v_R=0\})=0 ,
\end{equation}
where we set \(F(+\infty)=0\).
\end{lemma}

\begin{proof}
In \eqref{heat-potential} multiply the kernel by a continuous nondecreasing
function of \((t-s)/h\), equal to zero up to \(1/2\) and equal to one from
\(1\) on. Write \(v_{R,h}\) and \(v_{R_n,h}\) for the resulting potentials
with cutoff, including the initial term. For fixed \(h\), these converge
uniformly by compactness of the family of kernels with cutoff. Their
\(L^1(D)\) differences from the full potentials are bounded by the masses of
intervals of length \(h\) for the corresponding time marginals. More precisely, let $m_n$ be the time marginal of $R_n$ and
$m$ that of $R$. If $m$ has no atoms, then
\[
 \lim_{h\downarrow0}\limsup_{n\to\infty}
       \sup_t m_n([t-h,t+h])=0.
\]
Otherwise choose $h_k\downarrow0$, indices $n_k\to\infty$ and centers
$t_k\to t_*$. The masses of the intervals would stay bounded below,
so the closed-set inequality for narrow convergence on slightly
larger fixed intervals would give $m(\{t_*\})>0$, contradicting
atomlessness. Since the heat semigroup is $L^1$-contractive, this
controls the $L^1$ norms of the discarded short-time potentials,
uniformly in their observation times. This proves convergence of the solutions.

For fixed \(h\), positivity gives \(u_n\ge v_{R,h}-o(1)\) uniformly. On
compact subsets of \(\{v_{R,h}>0\}\), monotonicity gives
\(F(u_n)\le F(v_{R,h})+o(1)\). Testing against continuous functions
supported there and passing to the limit gives \(R\le F(v_{R,h})\Lambda\).
By Proposition~\ref{time-reduction}, we have \(R\ll\Lambda\). Let
\(h\downarrow0\) along a decreasing sequence. Then \(v_{R,h}\uparrow v_R\),
which gives the inequality on \(\{v_R>0\}\). If \(R\ne0\), the heat kernel
is positive and gives \(v_R>0\) at every time strictly after the first point
of the support of the time marginal of \(R\). This first time has zero mass
for \(R\), proving the assertion on the zero set.
\end{proof}

\begin{proof}[Proof of Theorem~\ref{evo2-range}]
Let \(R\) be a reaction limit. By Proposition~\ref{time-reduction},
\(R\ll\Xi\), and the reaction vanishes on the set carrying \(\Pi\). Its time
marginal is atomless. Lemma~\ref{potential-limit} gives convergence of the solutions
and $R\le F(v_R)\Xi$ on $\{v_R>0\}$ (after discarding the
$\Pi$-part by Proposition~\ref{time-reduction}). With the
convention $F(+\infty)=0$, this inequality immediately yields
$R(\{v_R=+\infty\})=0$; the same lemma gives
$R(\{v_R=0\})=0$. Therefore the measure
\[
 d\Gamma
  =\mathbf1_{\{0<v_R<\infty\}}
        F(v_R)^{-1}\,dR
\]
is well defined and finite. On the set $0<v_R<\infty$,
$F(v_R)^{-1}dR\le d\Xi$, while on its complement the measure
$\Gamma$ vanishes. Hence $0\le\Gamma\le\Xi$.
This also verifies explicitly that no reaction is hidden on the
infinite-potential set.

We show that \(R\) is the reaction obtained by completion for \(\Gamma\).
Write \(b=dR/d\Xi\), take a Kato exhaustion \(\Xi_j\uparrow\Xi\), and set
\(R_j=(b\wedge j)\Xi_j\). If \(\Xi=0\), the reaction is zero. Otherwise,
\(R_j\) are Kato measures and increase to \(R\). The associated heat
potentials \(v_j=v_{R_j}\) are continuous and bounded and increase pointwise
to \(v_R\). Define
\[
 d\Gamma_j=\frac{1}{F(v_j)}\,dR_j,
\]
using the right limit of $1/F$ at zero. Each $v_j$ is bounded, so
$1/F(v_j)$ is a bounded nonnegative multiplier. Thus $\Gamma_j$ is
Kato and finite. Since $R_j\uparrow R$, $v_j\uparrow v_R$ and
$1/F$ is nondecreasing, the measures $\Gamma_j$ increase.
Monotone convergence gives $\Gamma_j\uparrow\Gamma$ because
$R(\{v_R=\infty\})=R(\{v_R=0\})=0$ and $F(v_R)>0$ on the
remaining set. Positivity of the heat kernel and the absence of time atoms
give \(R_j(\{v_j=0\})=0\). Hence \(F(v_j)\Gamma_j=R_j\). By
Proposition~\ref{evo2-kato}, we obtain \(\mathcal R(\Gamma_j)=R_j\). Since each sequence increases to a finite measure, convergence
holds in total variation:
\[
 \|\Gamma-\Gamma_j\|_{\TV}=\Gamma(Q)-\Gamma_j(Q)\to0,
 \qquad \|R-R_j\|_{\TV}=R(Q)-R_j(Q)\to0.
\]
By the completion estimate \eqref{complete-reaction},
$\mathcal R(\Gamma_j)\to\mathcal R(\Gamma)$ in total variation.
As $\mathcal R(\Gamma_j)=R_j$, it follows that
$\mathcal R(\Gamma)=R$.

To prove attainment, fix \(\Gamma\le\Xi\) and choose a Kato exhaustion
\(\Gamma_j\uparrow\Gamma\). For fixed \(j\), mollify \(\Gamma_j\) as in
Proposition~\ref{evo2-kato}. By Corollary~\ref{evo2-diagonal}, we can
approximate \(\Lambda-\Gamma_j\) by smooth coefficients whose reactions have
vanishing mass. Adding these coefficients changes the reaction by \(o(1)\)
in total variation by \eqref{tv-add}. The sums of the coefficients converge
narrowly to \(\Lambda\) and the reactions converge to
\(\mathcal R(\Gamma_j)\). For each $j$ choose an approximation index $n(j)$ making both
bounded--Lipschitz errors less than $1/j$ relative to
$\Lambda$ and $\mathcal R(\Gamma_j)$, respectively. This is a
standard diagonal selection on the compact space $Q$, where the
bounded--Lipschitz metric metrizes narrow convergence on bounded
sets of measures. Since $\|\Gamma-\Gamma_j\|_{\TV}\to0$,
\eqref{complete-reaction} implies
$\mathcal R(\Gamma_j)\to\mathcal R(\Gamma)$ in total variation.
The selected smooth inputs therefore converge narrowly to $\Lambda$
and their reactions converge narrowly to the desired
$\mathcal R(\Gamma)$ along a full sequence.
\end{proof}

\noindent\textit{Stability under changes of initial data.}
For any finite input $\Lambda$ and bounded nonnegative initial data
$u_0,v_0$, denote the corresponding attainable-reaction sets by
$\mathcal A_{u_0}(\Lambda)$ and $\mathcal A_{v_0}(\Lambda)$.
The initial-data estimate in Lemma~\ref{reaction-tv} implies the
following two-sided stability: for every
$R\in\mathcal A_{u_0}(\Lambda)$ there exists
$R'\in\mathcal A_{v_0}(\Lambda)$ such that
\[
 \|R-R'\|_{\TV}\le\|u_0-v_0\|_{L^1(D)},
\]
and conversely with the initial data interchanged. To prove this,
choose a realizing sequence of smooth inputs for $R$ and solve
both initial-value problems with those same inputs. Extract a
narrowly convergent subsequence of the second reactions; the
first reactions still converge to $R$. Lower semicontinuity of
total variation for the difference of two narrowly converging
measure sequences gives the displayed bound from
\eqref{tv-pair}. This statement applies in particular to the
range in Theorem~\ref{evo2-range}, without identifying new
reaction limits or changing the underlying input measure.

\begin{remark}[Examples of Kato and polar inputs]
Let \(\Lambda=d\sigma_s\,d\nu(s)\), with
\[
 \sigma_s(B_r(x))\le C r^{d-k},\qquad
 \nu((t-h,t+h))\le C h^\alpha
\]
uniformly for small radii. If \(k<2\alpha\), then
\[
 k_\Lambda(h)\le C' h^{\alpha-k/2}.
\] Indeed, the
spatial Gaussian bound gives \(C(t-s)^{-k/2}\) for each convolution. Summing
over dyadic intervals of time delays gives the result by the bound on the
time marginal. Thus, with Lebesgue measure in time, one can take arbitrary
measurable translates of a fixed smooth hypersurface, or moving points in
dimension one. If \(k\ge2\alpha\), take a spatial density comparable to
\(|y|^{-k}\) near zero (a point mass if \(k=d\)) and a time density
\((t_*-s)^{\alpha-1}\) just before \(t_*\). The short-time potential at
\((t_*,0)\) is infinite.

The product of a fixed codimension-two subtorus and a compact interior time
interval is polar. In the two normal dimensions, take a logarithmic cutoff
between radii \(\eps^2\) and \(\eps\). Its squared gradient norm is
\(O(1/|\log\eps|)\) and its \(L^2\) norm tends to zero. Multiplication by a
fixed smooth time cutoff gives a test function whose norm in \eqref{par-cap}
tends to zero, including the term with the time derivative. Measures
concentrated on this product thus give zero limiting reaction for every
approximation.

Exhaustion by Kato measures is strictly weaker than the Kato condition. For
example, the density \(s^{-1}\mathbf1_{\{0<s<\delta,\,|y|<\sqrt{s}\}}\) on a
coordinate ball is integrable in space and time. At \((h,0)\), the
contribution from \(h/2<s<3h/4\) gives a positive lower bound for the
short-time heat potential, uniform as \(h\downarrow0\) by the Gaussian lower
bound.
\end{remark}

\paragraph{Scope of the reaction range.}
As noted in the canonical construction above, the class in Theorem~\ref{evo2-range} includes all
integrable space--time densities and, more generally, integrable spatial densities against any
finite atomless time measure. Proposition~\ref{prop:kato} shows that Kato exhaustibility also
implies diffuseness in the open cylinder. The converse, even after excluding time atoms, is not
established here. Time atoms themselves cannot be Kato-exhaustible, although spatially absolutely
continuous time impulses need not disappear; they are handled by Theorem~\ref{atom-range}. Thus
Theorem~\ref{evo2-range} is a range theorem on its specified class, not on all diffuse parabolic
measures. Whether its range description extends to every finite diffuse input with atomless time
marginal, or to more general diffusion operators, is a separate question. Any such extension must
justify the analogue of the canonical reaction construction and the total variation completion,
rather than only the elimination of polar sets.

\subsection{Accumulating impulses}
The set \(\mathcal A(\Lambda)\) is compact for narrow convergence. The
estimate for a positive lift and Lemma~\ref{reaction-tv} give a uniform
bound on the reaction masses. Closedness follows by a diagonal choice of
smooth approximants, since narrow convergence on bounded sets of measures on
\(Q\) is metrizable.

\begin{lemma}[Removing input mass]\label{thinning}
If \(0\le\lambda\le\Lambda\), then
\(\mathcal A(\lambda)\subset\mathcal A(\Lambda)\). For each
\(R\in\mathcal A(\Lambda)\) there exists \(R'\in\mathcal A(\lambda)\) such
that
\[
 \|R-R'\|_{\TV}\le\omega_F((\Lambda-\lambda)(Q)).
\]
The two sets are equal if \(\Lambda-\lambda\) is concentrated on a set of
the form in Proposition~\ref{time-reduction}.
\end{lemma}

\begin{proof}
To prove the inclusion, approximate \(\Lambda-\lambda\) with coefficients
whose reaction mass tends to zero and apply \eqref{tv-add}. Conversely,
choose \(a_n\) attaining \(R\). We can choose smooth \(0\le b_n\le a_n\)
whose input measures converge narrowly to \(\lambda\). Indeed, approximate
\(d\lambda/d\Lambda\) in \(L^1(\Lambda)\) by smooth functions with values
between zero and one. Multiply \(a_n\) by these functions, choosing indices
by a diagonal argument over continuous test functions. Let \(R'\) be a
subsequential reaction limit for \(b_n\). The discarded masses converge to
\((\Lambda-\lambda)(Q)\). Apply \eqref{tv-pair} for each positive lift, use
lower semicontinuity, and take the infimum over the lifts. In the last case
of the lemma, the reaction mass for the discarded input tends to zero for
every approximation, by Proposition~\ref{time-reduction}. Hence
\eqref{tv-add} gives \(R'=R\).
\end{proof}

\begin{corollary}[Mixed time measures]\label{mixed-range}
Suppose
\[
 \Lambda=\Xi+\Pi+\sum_i\delta_{t_i}\otimes\sigma_i
\]
is finite,
with \(\Xi,\Pi\) as in Theorem~\ref{evo2-range} and distinct times
\(t_i\in[0,T]\). Keep the first \(j\) impulses. Remove \(\Pi\) and restrict
\(\Xi\) to the complement of disjoint relative open time intervals about the
retained times, chosen to have total \(\Xi\)-mass less than \(1/j\). Write
\(\lambda_j\) for the new input and \(\delta_j\) for the removed mass
excluding \(\Pi\).

Define \(\mathcal B_j\) by choosing the reactions in order of time on the
remaining atomless intervals and at the retained impulses. Between
consecutive parts use the free heat flow to pass from the outgoing data to
the next incoming data. On an atomless interval take \(\mathcal R(\Gamma)\)
for any submeasure of the input on that interval. At an impulse with incoming state \(b\) and spatial density
\(f_i=d(\sigma_i)_{\rm ac}/dx\), take any jump density satisfying
\[
 0\le r_i\le J_{f_i}(b)-b.
\]
Then
\[
 \mathcal B_j=\mathcal A(\lambda_j)\subset\mathcal A(\Lambda),
 \qquad
 \sup_{R\in\mathcal A(\Lambda)}\inf_{S\in\mathcal B_j}
            \|R-S\|_{\TV}\le\omega_F(\delta_j).
\]
Precisely, for each $R\in\mathcal A(\Lambda)$ one can select
$S_j\in\mathcal B_j$ for all $j$ with
$\|S_j-R\|_{\TV}\le\omega_F(\delta_j)\to0$.
Conversely, any total-variation limit of a sequence
$S_j\in\mathcal B_j$ belongs to $\mathcal A(\Lambda)$.
This remains true when the impulse times accumulate.
\end{corollary}

\begin{proof}
Since the time marginal of $\Xi$ is atomless, we may choose pairwise
disjoint relative intervals around the finitely many retained impulse
times with total $\Xi$-mass less than $1/j$; their endpoints can be
chosen outside the atoms (there are none). Since the impulse masses
are summable, the discarded mass other than $\Pi$ satisfies
\[
 0\le\delta_j\le \frac1j+\sum_{i>j}\sigma_i(D)
 \longrightarrow0.
\] Fix $j$. Only finitely many impulses remain, and the deleted
neighborhoods isolate them from the remaining atomless pieces. Choose
one time in each resulting gap. For approximating coefficients,
construct smooth time cutoffs that equal one on the corresponding
closed pieces and vanish outside slightly larger intervals. Narrow
convergence lets the cutoffs be chosen so that the total input mass
removed from their union tends to zero. By \eqref{tv-pair}, this
modification changes the reaction in total variation by $o(1)$.
There are now finitely many consecutive independent input windows,
separated by intervals of free heat flow.

Pass to a subsequence on which the reactions in each window converge.
Between the windows, Duhamel's formula expresses the state as free heat
flow of the outgoing state plus earlier reaction potentials. Since the
observation time is separated by a positive gap from the previous
reaction support, the heat kernels and all their spatial derivatives
are uniformly bounded and continuous in the measure variables.
Consequently narrow convergence of the preceding reaction, together
with $L^1$ convergence of the preceding state, yields uniform spatial
convergence at the chosen gap times. At these times the heat kernels
are separated from zero delay. In particular, the next incoming data are
continuous and bounded even if the preceding jump is only in \(L^1\). The
initial-data estimate in \eqref{tv-pair} allows us to use the limits of the
incoming data. Theorems~\ref{atom-range} and~\ref{evo2-range} give the possible
reactions on each separated input window for this limiting incoming
state. In the converse direction, choose smooth realizing sequences
for each of the finitely many desired window reactions. Apply the
initial-data part of \eqref{tv-pair} to replace the approximating
incoming state by its limit without changing the reaction limit.
Concatenate these finitely many sequences across the gaps, allowing
free heat evolution in between, and choose a common approximation
index by a finite diagonal argument. The resulting smooth inputs
converge to $\lambda_j$ and realize precisely the selected reaction;
this proves $\mathcal B_j=\mathcal A(\lambda_j)$. To quantify this passage, first apply the final assertion of
Lemma~\ref{thinning} to remove the polar and spatially singular input
$\Pi$, which does not change the reaction range. Apply the second
assertion of that lemma to the remaining deletion of total mass
$\delta_j$. For each $R\in\mathcal A(\Lambda)$ this produces
$S_j\in\mathcal A(\lambda_j)=\mathcal B_j$ with
\[
 \|R-S_j\|_{\TV}\le\omega_F(\delta_j).
\]
Taking the infimum and supremum gives the stated bound. Since
$\delta_j\to0$, $S_j\to R$ in total variation. Conversely, any
total-variation limit of elements of $\mathcal B_j$ belongs to
$\mathcal A(\Lambda)$ because every $\mathcal B_j$ is contained in
$\mathcal A(\Lambda)$ and the latter is narrowly closed. This also
explains why accumulating impulse times require no extra condition.
\end{proof}
\section{Selection by a stretched time variable}\label{sec-clock}

When the spatial density of an impulse is given, a change of time realizes
the upper jump in Theorem~\ref{atom-range}. We also use this change of time
to obtain a stability estimate. Throughout this section, we take
\(F\in C^2(0,\infty)\) positive and nonincreasing.

Let $\nu$ be a finite positive Borel measure on $[0,T]$ with
$\nu(\{T\})=0$, and assume
$\Lambda=m(t,x)\,dx\,d\nu(t)$ where $m\ge0$ is
integrable with respect to $d\nu(t)\,dx$. Thus the conditional
spatial measures of $\Lambda$ have densities. Write the unique
decomposition of its time measure into an atomless and an atomic part as
\[
 \nu=\nu_{\mathrm{cont}}+
       \sum_{t\in[0,T]}\nu(\{t\})\delta_t,
 \qquad \nu_{\mathrm{cont}}(\{t\})=0
           \quad\text{for every }t.
\]
No absolute continuity of $\nu_{\mathrm{cont}}$ with respect
to $dt$ is assumed; singular-continuous time measures are allowed. Use \(\Phi,K,J_c\) from
Section~\ref{sec-impulse} and write \(V_D=|D|\).

\begin{proposition}[Stability on the stretched time interval]\label{evo-modulus}
For bounded nonnegative initial data \(u_{i,\mathrm{in}}\) and inputs of the
stated form, the construction gives right-continuous solutions in \(L^1(D)\)
with left limits. At a time atom of density \(c=m(t,\cdot)\nu(\{t\})\), the
outgoing state is \(J_c(u(t-))\). At zero we take the given incoming state
for \(u(0-)\). With \(\Delta_t=|\Lambda_1-\Lambda_2|([0,t]\times D)\), the
following estimates hold.
\begin{align}
 \int_D\Phi(|u_1(t)-u_2(t)|)\,dx
 &\le\int_D\Phi(|u_{1,\mathrm{in}}-u_{2,\mathrm{in}}|)\,dx+\Delta_t,
 \label{evo-modular-bound}\\
 \|u_1(t)-u_2(t)\|_1
 &\le\|u_{1,\mathrm{in}}-u_{2,\mathrm{in}}\|_1+V_DK(\Delta_t/V_D).
 \label{evo-modular-L1}
\end{align}
For identical initial data, the modulus in \eqref{evo-modular-L1} is sharp
at every positive value of \(\Delta_t\). On the atomless part, the reaction given by the construction is
$F(u(t,x))m(t,x)\,dx\,d\nu_{\mathrm{cont}}(t)$.
At each time atom its reaction measure is
$\delta_t\otimes(u(t)-u(t-))\,dx$.
\end{proposition}

\begin{proof}
Take the common time measure \(d\kappa=dt+d\nu_1+d\nu_2\). Its cumulative
function has continuous nondecreasing generalized inverse
\(\tau:[0,\kappa([0,T])]\to[0,T]\). This map sends Lebesgue measure to
\(\kappa\) and is constant on an interval for each atom. Set
\[
 b_i(s,x)=(d\Lambda_i/d(\kappa\otimes dx))(\tau(s),x).
\]
The generalized inverse is a nondecreasing continuous map: on the
interval of length $\kappa(\{t\})$ corresponding to an atom at time
$t$, it is constant. Because $\kappa\ge dt$, increasing the clock by
$\Delta s$ can advance physical time by at most $\Delta s$, so
$\tau$ is $1$-Lipschitz. The pushforward identity
$\tau_\#(ds)=\kappa(dt)$ is understood in the usual Stieltjes sense;
it is also obtained first for step measures and then by monotone
approximation. For bounded \(b_i\)
and initial data bounded below by \(e>0\), solve
\[
 U_i(s)=S(\tau(s))u_{i,\mathrm{in}}
       +\int_0^s S(\tau(s)-\tau(r))b_i(r)F(U_i(r))\,dr .
\]
Use a Lipschitz extension on the range from \(e\) to
\(\|u_{i,\mathrm{in}}\|_\infty+F(e)s\|b_i\|_\infty\) and iterate, using the
time integral of the \(L^1\) difference to bound its convolution. An
integrating factor gives comparison. The diffusion coefficient in stretched
time is \(\tau'\), and the state at physical time \(t\) is
\(U_i(\kappa([0,t]))\). On an interval of the stretched clock with $\tau(s)=t_*$,
there is no diffusion and $\partial_sU=b(s,x)F(U)$.
Integrating $\partial_s\Phi(U)=b$ gives, pointwise in $x$,
\[
 \Phi(U(s_+,x))-\Phi(U(s_-,x))
   =\int_{s_-}^{s_+}b(s,x)\,ds=c(x),
\]
so the outgoing state is exactly $J_c(U(s_-))$. To see independence of
the dominating measure, take a third finite measure dominating two
candidate clocks. Each old clock embeds monotonically into the new
one; by the substitution formula for their cumulative functions,
both embedded paths satisfy the same enlarged-clock mild equation.
Uniqueness for its bounded, positively lifted problem makes the paths
identical, and truncation removes these restrictions.

Set \(w=U_1-U_2\) and test the difference equation with the derivative of
\(\Phi(|w|)\). By convexity, the diffusion term is nonpositive. On \(w>0\),
the reaction term satisfies
\[
 \frac{b_2(F(U_1)-F(U_2))+(b_1-b_2)F(U_1)}{F(w)}
      \le(b_1-b_2)_+
\]
since \(w\le U_1\). Interchanging the indices gives the estimate on \(w<0\).
For rigor, use the smooth monotone approximation
\[
 \psi_\delta(w)=\operatorname{sgn}(w)
     \frac{|w|}{(|w|+\delta)F(|w|)},\qquad\psi_\delta(0)=0.
\]
On $w>0$, it is the product of the two nonnegative nondecreasing
functions $w/(w+\delta)$ and $1/F(w)$; therefore
$\psi_\delta'\ge0$ away from zero, and likewise on $w<0$.
Spatial convolution of the difference equation permits the chain
rule and the periodic integration by parts
\[
 \int_D\psi_\delta(w)\Delta w\,dx
      =-\int_D\psi_\delta'(w)|\nabla w|^2\,dx\le0.
\]
The reaction contribution on $w>0$ is bounded by $(b_1-b_2)_+$
using $w\le U_1$ and $F(U_1)\le F(w)$; the same argument with the
indices reversed applies on $w<0$. The positive lower bound for
$U_i$, bounded reaction coefficients and truncated $F$ permit passage
from smoothed paths to the original ones. Then monotone convergence
in $\delta\downarrow0$ identifies the primitive of $\psi_\delta$
with $\Phi(|w|)$ and gives the modular inequality. Integration up to a
time in the stretched variable gives \eqref{evo-modular-bound}, since the
change of time preserves the total variation of the input difference. For equal initial data, Jensen's inequality and the convexity of
$\Phi$ yield
\[
 \Phi\left(\frac{\|w(t)\|_1}{V_D}\right)
 \le\frac1{V_D}\int_D\Phi(|w(t)|)\,dx
 \le\frac{\Delta_t}{V_D}.
\]
Inverting gives $\|w(t)\|_1\le V_DK(\Delta_t/V_D)$.
For different initial data, compare $U_1$ with the solution for
$\Lambda_1$ starting from the second initial datum. The
positive-part contraction for equal inputs bounds this intermediate
state difference by $\|u_{1,\mathrm{in}}-u_{2,\mathrm{in}}\|_1$.
Apply the equal-initial-data estimate to the two remaining states,
and use the triangle inequality. This proves \eqref{evo-modular-L1}.

The estimates hold at each intermediate stretched time. For an integrable density $b_i$, take the increasing truncations
$b_i^{(j)}=\min(b_i,j)$. Since $\|b_i^{(j)}-b_i^{(k)}\|_1\to0$,
\eqref{evo-modular-L1} makes the associated states Cauchy in
$L^\infty_sL^1_x$; its modulus tends to zero by continuity of $K$
at zero. The same inequality permits the lift of the initial data
to decrease to zero without affecting the limit. With a positive initial lower bound, the solutions increase with the
truncations, and domination by \(b_i F(e)\) gives the mild equation. As the
lift decreases to zero, contraction gives a uniform \(L^1\) limit and the
reactions increase. Integrating the mild equation in space and comparing
with the unforced solution by \eqref{evo-modular-L1} bounds the reaction
masses by \(V_DK(\Lambda_i(Q)/V_D)\). Monotone convergence identifies the
reaction at the limit, with the conventions in \eqref{heat}. Pushforward by
\(\tau\) gives the reaction in physical time. Integration on each interval
of constant physical time gives the jump. Uniform \(L^1\) convergence and
Fatou's lemma give the time traces and both estimates. From zero, a
spatially constant impulse of density \(c>0\) gives the state \(K(c)\) and
an input difference of \(cV_D\) in total variation. This attains the
modulus.
\end{proof}

\begin{remark}[Spatially resolved regularizations]
A narrow approximation converges to the solution given by the change of time
if the smooth densities \(a_n\) have uniformly bounded
\(\int_0^T\|a_n(t)\|_\infty\,dt\) and if, for a smooth spatial approximate
identity \(P_j\),
\[
 \lim_{j\to\infty}\sup_n\int_0^T\|P_j a_n-a_n\|_1\,dt=0 .
\]
First add a positive constant to the initial data and fix the spatial
smoothing. Set
\[
 \begin{aligned}
 H_n(t)&=\int_0^t P_j a_n(s)\,ds,\\
 H(t)&=\int_{[0,t]}P_jm(s)\,d\nu(s).
 \end{aligned}
\] By narrow convergence and compactness
of the smoothed kernel tests, \(H_n\to H\) in \(L^1_t C^2_x\) and at
continuity times. For both the smooth evolution and the evolution with the
change of time, consider \(Y=\Phi(u)-H\). This quantity is continuous across
the intervals of constant physical time. In the physical time variable it
satisfies
\[
 Y_t-\Delta Y=\Delta H+
      F'(K(Y+H))|\nabla(Y+H)|^2 .
\]
By the chain rule, \(\Phi(U)-H\) has derivative \(\tau'\Delta U/F(U)\) in
stretched time. With smoothed initial data, differentiate the equation once
and twice in space. Comparison gives uniform derivative bounds, using the
integrated \(C^2_x\) bound for the smoothed input and the sign of the
zeroth-order term. The estimates also apply with diffusion coefficient
\(\tau'\ge0\) in stretched time. Including the time marginal of \(\Lambda\)
in \(\kappa\) makes the smoothed density bounded in stretched time. On the
positive range, the difference equation for the compensated unknowns has
bounded coefficients. The maximum principle gives
\[
 \|Y_n-Y\|_\infty\le C_j\int\|H_n-H\|_{C^2}dt.
\]
We remove spatial
smoothing by \eqref{evo-modular-L1} and dominated convergence of
\(\|P_jm(t)-m(t)\|_1\) against \(\nu\) for the limiting input. Contraction
for fixed input removes the lift and the initial smoothing. This proves
convergence in space--time \(L^1\) and in \(L^1(D)\) at continuity times of
the cumulative spatial input.

For a fixed \(C^2\) multiplier \(m\), positive \(C^2\) initial data and
equal-mass time measures \(\nu_1,\nu_2\), take \(H_i=m A_i\) and
\(A_i(t)=\nu_i([0,t])\). The same estimate gives
\[
 \|Y_1-Y_2\|_\infty
 +\int_0^T\|u_1(t)-u_2(t)\|_\infty\,dt
 \le C W_1(\nu_1,\nu_2),
\]
since
\[
 \|A_1-A_2\|_{L^1_t}=W_1(\nu_1,\nu_2).
\]
Moving a spatially constant impulse attains the linear rate in the
time-integrated state norm.
\end{remark}

\section{Discussion and open questions}\label{sec-discussion}
The preceding results distinguish three types of information lost under weak approximation. In the
critical local Hardy equation, the dyadic energy criterion uses both radial tails and uniform
control of angular distributions; the general capacity tests use different data. The endpoint
moments are additional observations not determined by energy convergence alone. At inverse power
one, the condenser discrepancy replaces radial information by a test of the relative masses on sets
separated by a capacitary barrier.

In the fractional problem the defect appears in a different form. The power pairing is the sum of a
nonnegative remainder and an explicit pole term when a trace exists. Boundedness of the profile
already controls the difference of the localized quantities; the trace determines its limiting
value. The failure of boundedness permits positive sources whose remainder is finite but whose
pairing diverges.

For the parabolic problem, the limit of the input is insufficient to specify the reaction, even for
a single absolutely continuous impulse or an integrable density. The range theorem separates the
input components which always disappear from those for which every canonical reaction of a
submeasure can be attained. Additional restrictions on the approximating sequence may select a
particular element of that range.

These conclusions suggest three natural directions: geometric conditions beyond the bounded
angular-density hypothesis that lead to necessary and sufficient energy tests; conditions on
fractional profiles that guarantee a trace rather than merely finite remainder; and extension of the
parabolic range theorem beyond Kato-exhaustible inputs. The present results do not require these
extensions.

\section*{Acknowledgments}
The author thanks Sharif University of Technology for its support.

\section*{Declarations}
\noindent\textit{Funding.}
This research was supported by the Iran National Science Foundation
(INSF), project No.~40505687.

\smallskip
\noindent\textit{Competing interests.}
The author declares no relevant financial or non-financial competing interests.

\smallskip
\noindent\textit{Data availability.}
No datasets were generated or analysed in this theoretical study.

\smallskip
\noindent\textit{Use of generative AI.}
Generative-AI tools, including OpenAI ChatGPT, were used in preparing
and revising this manuscript for language and \LaTeX{} editing,
bibliographic cross-checks, critical review, and drafting and checking
selected mathematical arguments, including supplementary results and
proof refinements. These tools are not authors. The author retains full
responsibility for the mathematical claims, computations, citations,
and the submitted version.

\end{document}